\documentclass[a4wide,10pt]{article}%
\usepackage{geometry}                % See geometry.pdf to learn the layout options. There are lots.
\usepackage{tikz} 
\usetikzlibrary{positioning}
\usepackage{graphicx}
\usepackage{authblk}
\usepackage{amssymb}
\usepackage{epstopdf}
\usepackage{subfigure}  % use for side-by-side figures
\usepackage[sans]{dsfont}
\usepackage[applemac]{inputenc}
\usepackage[english]{babel}
\usepackage{latexsym}
\usepackage{mathrsfs}
\usepackage{amscd}
\usepackage{graphicx}
\usepackage{color}
\usepackage{float}
\usepackage{amsmath}
\usepackage{amsfonts}
\numberwithin{equation}{section}
\usepackage{enumerate}
\usepackage{amsthm}
\usepackage{algorithm}
\usepackage{algorithmic}
\usepackage[numbers,sort]{natbib}
\usepackage[bookmarks=true,colorlinks=true,linkcolor={blue},urlcolor={blue}, citecolor={blue},pdfstartview={XYZ null null 1.22}]{hyperref}%

\newcommand{\R}{\mathbb R}

\def\be#1\ee{\begin{equation}#1\end{equation}}

\theoremstyle{definition}

\newtheorem{remark}{Remark}

\newcommand{\xx}{{\bf x}}

\def\RR{\mathbb R}

\def\p{\rho}

\def\T{\mathcal T}

\def\be{\begin{equation}}
	\def\ee{\end{equation}}
\def\bea{\begin{eqnarray}}
	\def\eea{\end{eqnarray}}

\title{Kinetic description of swarming dynamics with topological interaction and transient leaders}

\author{Giacomo Albi\footnote{	Dipartimento di Informatica, Universit\`a di Verona, Verona, Italy, e-mail: giacomo.albi@univr.it}, and Federica Ferrarese\footnote{		Dipartimento di Matematica, Universit\`a di Trento, e-mail: federica.ferrarese@unitn.it}}

\begin{document}
	\date{}
	\maketitle
	
	%% Enter the first author's name and address:
	%\centerline{\scshape Giacomo Albi}
	%\medskip
	%{\footnotesize
	%% please put the address of the first author
	% \centerline{TU M\"unchen, Faculty of Mathematics}
	%   \centerline{Boltzmannstra\ss e 3, D-85748, Garching (M\"u√ºnchen), Germany}
	%%   \centerline{ Springfield, MO 65801-2604, USA}
	%} % Do not forget to end the {\footnotesize by the sign }
	%
	%\medskip
	%
	%\centerline{\scshape Lorenzo Pareschi and Mattia Zanella}
	%\medskip
	%{\footnotesize
	% % please put the address of the second  and third author
	% \centerline{ University of
	%Ferrara, Department of Mathematics and Computer Science}
	%   \centerline{Via N. Machiavelli 35, 44121, Ferrara, Italy}
	%%   \centerline{Springfield, MO 65810, USA}
	%}
	%
	%\bigskip
	%
	%% The name of the associate editor will be entered by an editorial staff
	%% "Communicated by the associate editor name" is not needed for special issue.
	%% \centerline{(Communicated by the associate editor name)}
	%
	%
	%%The abstract of your paper
	\begin{abstract}

In this paper, we present a model describing the collective motion of birds.  The model introduces spontaneous changes in direction which are initialized by few agents, here referred as leaders, whose influence act on their nearest neighbors, in the following referred as followers. Starting at the microscopic level, we develop a kinetic model that characterizes the behaviour of large flocks with transient leadership. One significant challenge lies in managing topological interactions, as identifying nearest neighbors in extensive systems can be computationally expensive. To address this, we propose a novel stochastic particle method to simulate the mesoscopic dynamics and reduce the computational cost of identifying closer agents from quadratic to logarithmic complexity using a $k$-nearest neighbours search algorithm with a binary tree.
Lastly, we conduct various numerical experiments for different scenarios to validate the algorithm's effectiveness and investigate collective dynamics in both two and three dimensions.	\end{abstract}
	{\bf Keywords:} mean-field models, kinetic equations, Monte-Carlo methods, topological interactions, transient leadership
	\\
	{\bf AMS classification}: 65C05, 65Y20, 82C80, 92B05
	
	\tableofcontents
\section{Introduction}
\label{sec:intro}
In the last decades, there has been a notable surge in interest regarding the study of mathematical models describing collective behaviour of animals such as bacterial swarm \cite{koch1998social}, self-organization in insects \cite{dussutour2006collective,bonabeau1997self},  bird flocking \cite{cucker2007emergent,parrish1999complexity,ballerini2008empirical,lukeman2010inferring}, and  fish schooling \cite{hemelrijk2008self,d2006self}.
This captivating area of investigation has garnered substantial interest, with researchers increasingly delving into the complexities of emergent behaviours exhibited by natural systems, but also spanned to a wider range of applications such as swarm of robots \cite{jadbabaie2003coordination,choi2019collisionless,king2023biologically}, as well as social sciences and economics, \cite{albi2016opinion,rainer2002opinion,toscani2006kinetic,dimarco2020wealth}, vehicular and pedestrian traffic \cite{cristiani2014multiscale,albi2019vehicular,bressan2014flows}.

These large ensemble of models incorporate rules governing the behaviour of individual entities within the system. By integrating such mechanisms, these models effectively capture the impact of each entity on others, taking into account their relative positions and velocities.  In this manuscript, our focus centres around the dynamics governing animal swarms, building upon the recent model proposed in \cite{cristiani2021all}. This model introduces spontaneous changes of direction within the swarm, which are assumed to be independent of external factors like for example predators. The concept has been thoroughly analyzed in \cite{attanasi2014information,attanasi2015emergence} , serving as the foundation of our research. In general, birds, such as European starling (\textit{Sturnus vulgaris}), move together in huge groups. Indeed, by acting as a unified entity, they can effectively reduce individual vulnerability to predators. However, within these cohesive flocks, some members may deviate from the established course in response to various stimuli, such as the presence of food sources.   This dynamic interplay between collective safety in numbers and individual responsiveness to environmental highlights the intricate balance between social cohesion and adaptive behavior within birds communities. As birds navigate their surroundings in unison, their synchronized movements reflect both the collaborative effort to evade threats and the opportunistic pursuit of sustenance. 
	In this context, we suppose to have two dynamic subpopulations labelled as leaders and followers. Birds who initiate turns are referred as leaders, while their nearest neighbours, who adjust their motion accordingly, are regarded as followers.   This change of labels is characterized as a stochastic process, where each occurrence represents a random event. Our primary interest lies in exploring the phenomenon of \textit{transient leadership}, wherein agents can alter their labels over time, as studied for example in \cite{albi2019leader,loy2020non}, and also at different scales in \cite{li2015cucker, morandotti2020mean, bernardi2021leadership,cristiani2023kinmacr}.
 Here, we examine an extended version of the second-order stochastic differential equations presented in \cite{cristiani2021all}. 
	Similarly, we consider that each agent (bird) can interact with a maximum of $M$ nearest neighbours, consistent with observations in \cite{ballerini2008empirical}. Although we assume to remove delay effect, our model retains the capability to consider significant positions of interest, such as the existence of food sources or nest locations.  
Our objective is to study such dynamics at the mesoscopic scale, formally introducing a kinetic model of the swarms with topological-type interaction dynamics and deriving the associated mean-field limit. For alternative mean-field and kinetic models with topological interactions, we refer to \cite{blanchet2017kinetic, blanchet2016topological, haskovec2013flocking}, and for rigorous derivations, we refer specifically to \cite{degond2019propagation, benedetto2022mean}.  
Another primary objective of this study is to efficiently perform numerical simulations of high-dimensional non-local dynamics. One of the main challenges arises from the presence of topological-type interactions, necessitating the implementation of ad-hoc methods to reduce the complexity of the nearest neighbour search process. To address this computational burden at the mesoscopic scale, we introduce a novel stochastic simulation algorithm for the simulation of kinetic models such as \cite{pareschi2013interacting, albi2013binary,carrillo2019particle}, adopting $k-$nearest neighbours search strategy ($k$-NN), following the approach proposed in \cite{friedman1977algorithm}.  
	The primary innovation compared to classical $k$-NN algorithms lies in generating an estimate of the nearest neighbors rather than computing them precisely. This is achieved by performing the $k$-NN search within a subset of the total sample size.  By implementing this method, we successfully reduce the computational complexity of the nearest neighbour search from quadratic to logarithmic scale, significantly enhancing the efficiency of numerical simulations.

The paper is organized as follows. In \ref{sec:micro_model} we introduce the microscopic model describing which are the forces that act on followers and on leaders. In \ref{sec:kinetic_model}, we extend the study to the kinetic level, describing the evolution of the densities and how the change of labels occurs. In \ref{sec:numerical_methods}, we introduce the algorithms that can be used to simulate the binary interaction rules and the change of labels. In \ref{sec:validation} we perform two different validations experiments, testing both the accuracy and the efficiency of the numerical methods introduced. In particular, we show how it is possible to reduce the computational costs in dealing with non-locality. In \ref{sec:2D3Dexperiments} we simulate the dynamics at the microscopic and kinetic level for both the two and three dimensional cases.

%%%%%%%%%%%%%%%%%%%%%%%%%%%%%%%%%%%%%%%%%%%%%
%%%%%%%%%%%%%%%%%%%%%%%%%%%%%%%%%%%%%%%%%%
%%%%%%%%%%%%%%%%%%%%%%%%%%%%%%%%%%%%%%%%%

\section{Swarming models with leaders-followers dynamics}\label{sec:micro_model}
We consider a large system of  $N$ interacting agents represented by points moving in a $d$-dimensional space with an evolving hierarchy of interactions ruled by follower-leaders dynamics. 
For every $i = 1,\ldots,N$, let $(x_i(t),v_i(t)) \in \RR^{2d}$ denote position and velocity of the $i$-th agent at time $t$, with $d=1,2,3$, and $\lambda_i(t)\in\Lambda \equiv {\left\{0,1\right\}}$ the space of labels indicating  at time $t$ the status of agent $i$ to be either {\em follower} ($F$) for  $\lambda_i(t)=0$, or {\em leader}  ($L$) for $\lambda_i(t)=1$.   Moreover we account for $N_{{src}}$ fixed target positions located at $x_k^{{src}}\in\RR^{d}$ for $k = 1,\ldots, N_{{src}} $, indicating positions of interested for the swarm such as nest, or foraging areas \cite{bernardi2018particle,bidari2019social}. 

We assume the system of agents  evolving according to ODEs system,
\begin{equation}\label{eq:dynamics} 
	\begin{aligned}
		\dot{x}_i &= v_i,\cr
		\dot{v}_i &= \frac{1}{M} \sum\limits_{\left\{j~:~x_j\in \mathcal{B}_M(x_i;\xx)\right\}}\left[ A^{rep}(x_i,x_j) + \left(  A^{ali}(v_i,v_j) + A^{att}(x_i,x_j)\right)(1-\lambda_i(t))\right] \cr
		&\qquad \qquad\qquad+ \left[ A^{src}(x_i) + A^{ctr}(x_i) +S(v_i) \right]\lambda_i(t),  \qquad i=1,\ldots,N, \cr
	\end{aligned}
\end{equation}
%where, the set $\mathcal{T}_M (x_i;\xx)$ is the topological index set of the $M$ closest agent to the agent $i$, i.e.
%\[
%\mathcal{T}_M (x_i;\xx) := \left\{j | x_j\in\mathcal{B}_M(x_i;\xx) \right\}
%\]
where we denoted by $\mathcal{B}_M(x_i;\xx)$ the ball centred at $x_i$, with $\xx=(x_1,\ldots,x_N)$, containing the $M$ nearest neighbours to $i$-agent,   assuming that in case of ambiguity, e.g. more than one agent is at the same distance from agent in position $x_i$, we select the first $M$ agents giving priority according to the indexing order. 
Hence, the dynamics encodes different behaviours according to the value of the label $\lambda_i(t)$. 
% with the lower index is included, and $\chi_{B_M(x;\xx,\yy)(y)}$ is the characteristic function of the topological ball $\mathcal{B}_M(x;\xx)$. Hence, 
\begin{itemize}
	\item 
	For $\lambda_i(t) = 0$, we have follower-type interactions characterized by
	\begin{itemize}
		\item  repulsion force
		\begin{equation}\label{eq:repulsion} 
			A^{rep}(x,x') =-C_{rep} \frac{ x' - x}{\Vert  x' - x\Vert^2 },
		\end{equation}
		\item alignment force
		\begin{equation}\label{eq:aligment} 
			A^{ali}(v,v') =C_{ali} \left( v' - v\right),
		\end{equation}
		\item and attraction force  
		\begin{equation}\label{eq:atraction} 
			A_i^{att}(x,x') =C_{att} \left( x' - x\right),
		\end{equation}
	\end{itemize} 
	where $C_{rep}\geq0$, $C_{ali}\geq0$ and $C_{att}\geq0$ are non-negative constants.
	\item For $\lambda_i(t)=1$, we have leaders-type dynamics  characterized by a repulsion force defined as in equation \eqref{eq:repulsion} and by a self-propulsion friction term $S(\cdot)$ defined as 
	\begin{equation}\label{eq:relaxation}
		S(v) = C_v (s-\Vert v \Vert^2)v,
	\end{equation}
	where $s$ is a given characteristic speed and $C_v\geq0$. In presence of sources terms, leaders are driven by 
	\begin{equation}\label{eq:food}
		A^{src}(x) = C_{src} \sum_{k=1}^{N_{src}}\varphi_\epsilon(\Vert x_k^{src}-x\Vert;\overline r) \frac{x_k^{src}-x}{\Vert x_k^{src}-x\Vert },
	\end{equation}
	\textcolor{black}{where $C_{src}\geq0$, $x_k^{src}$ denotes the position of the attraction source (nest, or food) and $ \varphi_\epsilon(\cdot)$ is a sigmoid function of the following type
		\begin{equation}
			\varphi_\epsilon(r;\overline r) := 	\frac{1}{1+\exp\{(r-\overline r)/\epsilon\}},
		\end{equation}
		with regularization parameter
		$\epsilon >0$, modelling a {perception area around the source} activating when the distance of the agent is below the threshold value $\overline r>0$. 
		Furthermore, leaders can  be forced to move toward the centre of mass $x_c$ according to the force
		\begin{equation}\label{eq:centre_of_mass}
			A^{ctr}(x) = C_{ctr} \Big(1-{\varphi_\epsilon} (\Vert x_c - x \Vert;\underline r)\Big) \frac{x_c - x}{\Vert x_c - x \Vert },
		\end{equation}
		where $C_{ctr}\geq0$, and when the distance with respect to the centre of mass is larger than  $\underline r$.
		%\begin{equation}
		%	{\varphi}(r) = (\frac{e^{[\bar{c}(r-\bar{r})]}}{1+\exp\{[\bar{c}(r-\bar{r})]\}}),
		%	\end{equation}
		%where $\bar{c}\geq0$ and $\bar{r}\geq0$. 
	}
\end{itemize}
\subsection{Stochastic process for leaders emergence}\label{sec:leaders_micro}
Agents have the ability to switch between being leaders and followers, and vice versa. Such a change in status is treated as a stochastic process, where each occurrence represents a random event governed by an assigned probability distribution. Each event is associated with a transition rate, which quantifies the probability of its occurrence per unit time. Therefore, for $\boldsymbol{\lambda} = (\lambda_1,\ldots,\lambda_N)$, each label $\lambda_i(t)$ will follow a jump process in this manner
\begin{itemize}
	\item if $\lambda_i(t) = 1$ then it switches to $0$ with rate $\pi_{L\to F}(t,x_i,v_i,\lambda_i;\xx,\boldsymbol{v},\boldsymbol{\lambda})$,
	\item if $\lambda_i(t) = 0$ then it switches to $1$ with rate $\pi_{F\to L}(t,x_i,v_i,\lambda_i;\xx,\boldsymbol{v},\boldsymbol{\lambda})$,
\end{itemize}
where the transition rates $\pi_{F\to L}(\cdot), \pi_{L\to F}(\cdot)$ in general are non-linear functions of the state of the system. In what follows we will consider different choices for the labels' switching rules, ranging from random, density dependent and aiming at organizing agents toward a common target. These choices will be detailed in \ref{sec:leaders_kinetic}.
 \begin{remark}
		Similar to the model presented in \cite{cristiani2021all}, we suppose that both leaders and followers experience repulsive forces to prevent collisions with other agents while followers are subjected also to alignment and attraction forces. Our model incorporates additional dynamics for the leaders, including a relaxation term towards a desired velocity to ensure that leaders remain in contact with followers. We also introduce two attraction terms: one towards a potential source, such as the nest or food position, and another towards the centre of mass to prevent excessive splitting of the flock. We also explore a distinct dynamics in the labels switching process compared to the seminal work.
\end{remark} 
%%%%%%%%%%%%%%%%%%%%%%%%%%%%%%%%%%%%%%%%%%%%%
%%%%%%%%%%%%%%%%%%%%%%%%%%%%%%%%%%%%%%%%%%
%%%%%%%%%%%%%%%%%%%%%%%%%%%%%%%%%%%%%%%%%
\section{Kinetic modelling of swarming dynamics}\label{sec:kinetic_model}
In this section, we will provide a kinetic description of the swarming model with leader emergence and topological interaction, we refer to \cite{albi2016opinion,albi2016invisible,loy2020non} for related studies in the context of kinetic models.

Thus, we associate to each agent a position and velocity $(x,v)\in\R^d\times\R^d$ and a leadership-level $\lambda$, as a discrete binary variable in the label space $\Lambda = \{0,1\}$. We are interested in the evolution of the probability density function 
\begin{equation}\label{eq:def_f}
	f=f(x,v,\lambda,t), \qquad f: \R^d\times\R^d\times\left\{0,1\right\} \times \R_+\rightarrow \R_+,
\end{equation}
where $t\in\RR^+$ denotes as usual the time variable. For each time $t\ge 0$, $\lambda\in\{0,1\}$, we have  the following marginal density
\begin{equation}\label{eq:integration_wc}
	\p(\lambda,t)=\int_{\R^{d}\times \R^d}f(x,v,\lambda,t)d(x,v), 
\end{equation}
which defines the quantity of agents with label $\lambda$ at time $t$. In the sequel, we will assume that the total number of agents is conserved, namely 
\be\label{eq:mass}
\p(1,t)+\p(0,t)=1.
\ee
Likewise, we define the marginal density  for agents in space and velocity
\begin{equation}\label{eq:total_density}
	g(x,v,t)=\sum_{\lambda} f(x,v,\lambda,t),\qquad \lambda\in \{0,1\}.
\end{equation}
Next, we assume the density $f(x,v,\lambda,t)$ to be solution of a kinetic equation accounting for pairwise interactions among agents, and for labels transition. 

\paragraph{Notational convention}  To ease the writing, we will use an equivalent notation for functions  depending on $\lambda$,  where we introduce the indexing given by the discrete label space $\Lambda$, as follows
\[
F_\lambda(\cdot) := F(\cdot,\lambda).
\]
Then, for example, the density $f(x,v,t,\lambda)$ will be denoted by $f_\lambda(x,v,t)$ or the mass $\rho(\lambda,t)$ by $\rho_\lambda(t)$.
\subsection{Povzner-Boltzmann-type model}\label{sec:kinetic_model_0}
We assume that each agent modifies its velocity through a binary interaction occurring with another agent within the topological ball $\mathit{B}_{r^*}(x,t)$, the ball centred in $x$ whose radius is defined, for a fixed $t \geq 0$, by the following variational problem 
\begin{equation}\label{eq:radius}
	r^*(x,t) = \arg \min_{\alpha>0} \left\{\sum_{\lambda} \int_{\mathit{B}_\alpha (x,t)\times\RR^{d}} f_\lambda(x,v,t) dx dv \geq \rho^* \ \right\},
\end{equation}
where $\rho^*\in(0,1]$ is the target topological mass, namely the ratio $\rho^* = M/N$  associated to the microscopic model  \eqref{eq:dynamics}.

Hence, we consider pairwise interactions among an agent with state $(x,v,\lambda)\in\RR^{2d}\times{\left\{0,1\right\}}$ and $(x_*,v_*,\lambda_*)\in\mathit{B}_{r^*}(x,t)\times\RR^d\times{\left\{0,1\right\}}$, where the post-interaction velocities are given by
\begin{equation}\begin{cases}\label{eq:binary}
		v'   &=v   + \alpha \mathcal{F}_\lambda(x,x_*,v,v_*), \\
		v_*' &=v_*,%+ \alpha \mathcal{F}(x_*,x,v_*,v,\lambda_*),
\end{cases}\end{equation}
where $v,v_*\in\mathbb{R}^d$ denote the pre-interaction velocities and $v',v_*'$ the velocities after the exchange of information between the two agents. In  \eqref{eq:binary} we assume 
\begin{equation}\label{eq:binary2}
	\begin{split}
		&\mathcal{F}_\lambda(x,x_*,v,v_*) =A^{rep}(x,x_*) +\bigl[ A^{ali}(v,v_*)+A^{att}(x,x_*)\bigl](1-\lambda) \\
		&\qquad \qquad\qquad \qquad \qquad \qquad + \bigl[ A^{src}(x) + A^{ctr}(x) + S(v)\bigr]\lambda. 
	\end{split}
\end{equation}
%% &\mathcal{A}^{food}(\cdot) = A^{food}(\cdot) M, \quad \mathcal{A}^{centre}(\cdot) = A^{centre}(\cdot) M,
%
%\end{equation}
For $\lambda\in \{0,1\}$, the evolution in time of the density function $f_\lambda(x,v,t)$ is described by a integro-differential equation of the Povzner-Boltzmann type  \cite{povzner1962boltzmann,fornasier2011fluid}  as follows
\begin{equation}\label{eq:boltz_lin}
	\partial_t f_\lambda(x,v,t) + v\cdot\nabla_x f_\lambda (x,v,t)-\T_\lambda[f](x,v,t)= Q_\lambda(f,f)(x,v,t),
\end{equation}
where  $\mathcal T_\lambda [f](\cdot)$ accounts for the evolution of the agents in the discrete label space and $Q_\lambda(\cdot,\cdot)$ is  the interaction operator defined as follows
\begin{small}
	\begin{equation}	\label{eq:Bo}
		Q_\lambda(f,f)(x,v,t) =\eta \sum_{\lambda_*}\int_{\Omega} \left(\dfrac{1}{J_\lambda}f_{\lambda}(x,'v,t)f_{\lambda_*}(x_*,'v_*,t)-f_\lambda(x,v,t)f_{\lambda_*}(x_*,v_*,t) \right)d(x_*,v_*),
	\end{equation}
\end{small}
where $\Omega = \mathit{B}_{r^*(x,t)}\times\mathbb{R}^d$, $('v,'v_*)$ are the pre-interaction velocities, and the term $J_\lambda$ denotes the Jacobian of the transformation $(v,v_*)\rightarrow (v',v_*')$ with  $ (v',v_*')$ the post-interaction velocities, and $\eta>0$ is a constant relaxation rate representing the interaction frequency.

% the kernels $'B,B$ define the binary interaction .
%We will consider constant interaction kernels of the following form 
%\begin{equation}
%	B_{(v,v_*)\rightarrow (v',v_*')}=\eta,
%\end{equation}

%where $\eta>0$ is a constant relaxation rate representing the interaction frequency. 

%Here and in the rest of the Section, for notation simplicity, the explicit dependence from the time variable is omitted. 
\subsection{Master equation for leaders transition}\label{sec:leaders_kinetic}
In the previous section, we have introduced the transition operator $\mathcal T_\lambda [f](x,v,t)=\mathcal T[f](x,v,\lambda,t)$ characterizing the evolution of the agents in the discrete space of labels $\Lambda = \{0,1\}$ (followers/leaders). Such operator is defined as follows
\begin{equation}\label{eq:master_0}
	\begin{split}
		\mathcal{T}_0[f](x,v,t) =&  \pi_{L\to F}f_1(x,v,t)-\pi_{F\to L} f_0(x,v,t), \\
		\mathcal{T}_1[f](x,v,t) =&   \pi_{F\to L}f_0(x,v,t)-\pi_{L\to F}f_1(x,v,t),
	\end{split}
\end{equation}
where $\pi_{F\to L}:= \pi_{F\to L}(x,v,t;f)$ and $\pi_{L\to F}:= \pi_{L\to F}(x,v,t;f)$ are certain transition rates.

Thus the evolution of the transition process of labels can be described by the evolution equation for $\rho_\lambda (t)= \rho(\lambda,t)$, 
\begin{equation}\label{eq:Ldef2}
	\dfrac{d}{dt}\p_\lambda(t) - \int_{\RR^{2d}} \T_\lambda[f](x,v,t)\,d(x,v)=0.
\end{equation}
From the definition of the transition operator $\T_\lambda[\cdot]$ and \eqref{eq:mass} it follows the  conservation of the mass,
\begin{equation}
	\dfrac{d}{dt}\sum_{\lambda}\p_\lambda(t) =  \sum_{\lambda}\int_{\RR^{2d}} \T_\lambda[f](x,v,t)\,d(x,v) = 0.
	\label{eq:tnc}
\end{equation}
In the sequel we list possible choices of transition rates in \eqref{eq:master_0}.

\paragraph{Constant rates}
Leaders emerge with rate $q_{FL}>0$ and return to the followers status with rate $q_{LF}>0$. Hence, the transition rates write as follows
\begin{equation}\label{eq:rates_test_0}
	%	\begin{split}
	\pi_{L\to F}  = q_{LF},\qquad \pi_{F\to L} = q_{FL}. 
	%	\end{split} 
\end{equation} 
%Generally, assuming the transition rates indepedent on $(x,v)$, 
%\begin{equation*}
%	\pi_{L\to F}(\cdot) = \alpha(\lambda,t),\qquad \pi_{F\to L}(\cdot) = \beta(\lambda,t),
%\end{equation*}
Thus, we can rewrite equation \eqref{eq:Ldef2} as 
\begin{equation}\label{eq:master_constant}
	\begin{aligned}
		\partial_{t} \rho_1(t) &= q_{FL} \rho_0(t) - q_{LF} \rho_1(t),\\
		\partial_{t} \rho_0(t) &= q_{LF}	\rho_1(t)-q_{FL} \rho_0(t),
	\end{aligned}
\end{equation}
and find the stationary solution of equation \eqref{eq:master_constant} that is 
\begin{equation}\label{eq:stationary}
	\rho_1^\infty = \frac{q_{FL}}{q_{LF} + q_{FL}},\qquad 	\rho_0^\infty = \frac{q_{LF}}{q_{LF} +q_{FL}}.
\end{equation}

\begin{figure}[tbhp]

	\tikzset{every picture/.style={line width=0.7pt}} %set default line width to 0.75pt        
	\centering
	\begin{tikzpicture}[x=0.7pt,y=0.7pt,yscale=-0.85,xscale=0.85]
		%uncomment if require: \path (0,255); %set diagram left start at 0, and has height of 255
		
		%Shape: Ellipse [id:dp7251636312579726] 
		\draw  [color={rgb, 255:red, 0; green, 0; blue, 255 }  ,draw opacity=1 ][fill={rgb, 255:red, 0; green, 0; blue, 255 }  ,fill opacity=1 ] (178.93,101.91) .. controls (178.93,99.8) and (180.65,98.09) .. (182.77,98.09) .. controls (184.88,98.09) and (186.6,99.8) .. (186.6,101.91) .. controls (186.6,104.02) and (184.88,105.73) .. (182.77,105.73) .. controls (180.65,105.73) and (178.93,104.02) .. (178.93,101.91) -- cycle ;
		%Shape: Ellipse [id:dp04553318231320125] 
		\draw  [color={rgb, 255:red, 0; green, 0; blue, 255 }  ,draw opacity=1 ][fill={rgb, 255:red, 0; green, 0; blue, 255 }  ,fill opacity=1 ] (274.93,171.65) .. controls (274.93,169.54) and (276.65,167.83) .. (278.77,167.83) .. controls (280.88,167.83) and (282.6,169.54) .. (282.6,171.65) .. controls (282.6,173.76) and (280.88,175.47) .. (278.77,175.47) .. controls (276.65,175.47) and (274.93,173.76) .. (274.93,171.65) -- cycle ;
		%Shape: Ellipse [id:dp7606316889747426] 
		\draw  [color={rgb, 255:red, 0; green, 0; blue, 255 }  ,draw opacity=1 ][fill={rgb, 255:red, 0; green, 0; blue, 255 }  ,fill opacity=1 ] (272.93,123.83) .. controls (272.93,121.72) and (274.65,120.01) .. (276.77,120.01) .. controls (278.88,120.01) and (280.6,121.72) .. (280.6,123.83) .. controls (280.6,125.94) and (278.88,127.65) .. (276.77,127.65) .. controls (274.65,127.65) and (272.93,125.94) .. (272.93,123.83) -- cycle ;
		%Shape: Ellipse [id:dp6505884439295662] 
		\draw  [color={rgb, 255:red, 0; green, 0; blue, 255 }  ,draw opacity=1 ][fill={rgb, 255:red, 0; green, 0; blue, 255 }  ,fill opacity=1 ] (300.93,147.74) .. controls (300.93,145.63) and (302.65,143.92) .. (304.77,143.92) .. controls (306.88,143.92) and (308.6,145.63) .. (308.6,147.74) .. controls (308.6,149.85) and (306.88,151.56) .. (304.77,151.56) .. controls (302.65,151.56) and (300.93,149.85) .. (300.93,147.74) -- cycle ;
		%Shape: Ellipse [id:dp8022301909851632] 
		\draw  [color={rgb, 255:red, 0; green, 0; blue, 255 }  ,draw opacity=1 ][fill={rgb, 255:red, 0; green, 0; blue, 255 }  ,fill opacity=1 ] (243.93,191.58) .. controls (243.93,189.47) and (245.65,187.76) .. (247.77,187.76) .. controls (249.88,187.76) and (251.6,189.47) .. (251.6,191.58) .. controls (251.6,193.69) and (249.88,195.4) .. (247.77,195.4) .. controls (245.65,195.4) and (243.93,193.69) .. (243.93,191.58) -- cycle ;
		%Shape: Ellipse [id:dp4554681458514984] 
		\draw  [color={rgb, 255:red, 0; green, 0; blue, 255 }  ,draw opacity=1 ][fill={rgb, 255:red, 0; green, 0; blue, 255 }  ,fill opacity=1 ] (237.93,107.89) .. controls (237.93,105.78) and (239.65,104.07) .. (241.77,104.07) .. controls (243.88,104.07) and (245.6,105.78) .. (245.6,107.89) .. controls (245.6,110) and (243.88,111.71) .. (241.77,111.71) .. controls (239.65,111.71) and (237.93,110) .. (237.93,107.89) -- cycle ;
		%Shape: Ellipse [id:dp8616395558126551] 
		\draw  [color={rgb, 255:red, 0; green, 0; blue, 255 }  ,draw opacity=1 ][fill={rgb, 255:red, 0; green, 0; blue, 255 }  ,fill opacity=1 ] (203.93,153.72) .. controls (203.93,151.61) and (205.65,149.9) .. (207.77,149.9) .. controls (209.88,149.9) and (211.6,151.61) .. (211.6,153.72) .. controls (211.6,155.83) and (209.88,157.54) .. (207.77,157.54) .. controls (205.65,157.54) and (203.93,155.83) .. (203.93,153.72) -- cycle ;
		%Shape: Ellipse [id:dp6742795532769135] 
		\draw  [color={rgb, 255:red, 0; green, 0; blue, 255 }  ,draw opacity=1 ][fill={rgb, 255:red, 0; green, 0; blue, 255 }  ,fill opacity=1 ] (232.93,133.79) .. controls (232.93,131.68) and (234.65,129.97) .. (236.77,129.97) .. controls (238.88,129.97) and (240.6,131.68) .. (240.6,133.79) .. controls (240.6,135.9) and (238.88,137.61) .. (236.77,137.61) .. controls (234.65,137.61) and (232.93,135.9) .. (232.93,133.79) -- cycle ;
		%Shape: Rectangle [id:dp5679418222296349] 
		\draw  [fill={rgb, 255:red, 155; green, 155; blue, 155 }  ,fill opacity=1 ][line width=1.5]  (78.6,155.55) -- (87.6,155.55) -- (87.6,164.51) -- (78.6,164.51) -- cycle ;
		%Straight Lines [id:da6551084893720522] 
		\draw [line width=1.5]    (241.77,107.89) -- (195.26,125.91) -- (162.73,138.52) ;
		\draw [shift={(159.93,139.6)}, rotate = 338.81] [color={rgb, 255:red, 0; green, 0; blue, 0 }  ][line width=1.5]    (14.21,-4.28) .. controls (9.04,-1.82) and (4.3,-0.39) .. (0,0) .. controls (4.3,0.39) and (9.04,1.82) .. (14.21,4.28)   ;
		%Straight Lines [id:da9613870635210426] 
		\draw [color={rgb, 255:red, 128; green, 128; blue, 128 }  ,draw opacity=1 ][line width=1.5]    (241.77,107.89) -- (253.24,59.82) ;
		\draw [shift={(253.93,56.9)}, rotate = 103.42] [color={rgb, 255:red, 128; green, 128; blue, 128 }  ,draw opacity=1 ][line width=1.5]    (14.21,-4.28) .. controls (9.04,-1.82) and (4.3,-0.39) .. (0,0) .. controls (4.3,0.39) and (9.04,1.82) .. (14.21,4.28)   ;
		%Shape: Arc [id:dp9164304233240115] 
		\draw  [draw opacity=0] (218.61,116.1) .. controls (218.91,100.68) and (230.87,88.23) .. (245.94,86.91) -- (248.6,116.69) -- cycle ; \draw  [color={rgb, 255:red, 255; green, 0; blue, 0 }  ,draw opacity=1 ] (218.61,116.1) .. controls (218.91,100.68) and (230.87,88.23) .. (245.94,86.91) ;  
		%Shape: Ellipse [id:dp7314376753521818] 
		\draw  [color={rgb, 255:red, 0; green, 0; blue, 255 }  ,draw opacity=1 ][fill={rgb, 255:red, 0; green, 0; blue, 255 }  ,fill opacity=1 ] (252.94,153.55) .. controls (253.03,151.44) and (254.82,149.81) .. (256.94,149.9) .. controls (259.05,150) and (260.69,151.78) .. (260.6,153.89) .. controls (260.5,156) and (258.71,157.63) .. (256.59,157.53) .. controls (254.48,157.44) and (252.84,155.65) .. (252.94,153.55) -- cycle ;
		%Shape: Ellipse [id:dp8244688649375185] 
		\draw  [color={rgb, 255:red, 0; green, 0; blue, 255 }  ,draw opacity=1 ][fill={rgb, 255:red, 0; green, 0; blue, 255 }  ,fill opacity=1 ] (182.53,187.6) .. controls (182.53,185.49) and (184.25,183.78) .. (186.37,183.78) .. controls (188.48,183.78) and (190.2,185.49) .. (190.2,187.6) .. controls (190.2,189.7) and (188.48,191.41) .. (186.37,191.41) .. controls (184.25,191.41) and (182.53,189.7) .. (182.53,187.6) -- cycle ;
		%Shape: Ellipse [id:dp20344636722308684] 
		\draw  [color={rgb, 255:red, 0; green, 0; blue, 255 }  ,draw opacity=1 ][fill={rgb, 255:red, 0; green, 0; blue, 255 }  ,fill opacity=1 ] (279.33,87.96) .. controls (279.33,85.85) and (281.05,84.14) .. (283.17,84.14) .. controls (285.28,84.14) and (287,85.85) .. (287,87.96) .. controls (287,90.07) and (285.28,91.78) .. (283.17,91.78) .. controls (281.05,91.78) and (279.33,90.07) .. (279.33,87.96) -- cycle ;
		%Shape: Ellipse [id:dp1428596208773334] 
		\draw  [color={rgb, 255:red, 0; green, 0; blue, 255 }  ,draw opacity=1 ][fill={rgb, 255:red, 0; green, 0; blue, 255 }  ,fill opacity=1 ] (220.93,59.26) .. controls (220.93,57.15) and (222.65,55.44) .. (224.77,55.44) .. controls (226.88,55.44) and (228.6,57.15) .. (228.6,59.26) .. controls (228.6,61.37) and (226.88,63.08) .. (224.77,63.08) .. controls (222.65,63.08) and (220.93,61.37) .. (220.93,59.26) -- cycle ;
		%Shape: Rectangle [id:dp49670518896625704] 
		%	\draw   (13.73,38.68) -- (320.93,38.68) -- (320.93,243.48) -- (13.73,243.48) -- cycle ;
		%Shape: Ellipse [id:dp3339554726485452] 
		\draw  [color={rgb, 255:red, 0; green, 0; blue, 255 }  ,draw opacity=1 ][fill={rgb, 255:red, 0; green, 0; blue, 255 }  ,fill opacity=1 ] (502.13,102.09) .. controls (502.13,99.98) and (503.85,98.27) .. (505.97,98.27) .. controls (508.08,98.27) and (509.8,99.98) .. (509.8,102.09) .. controls (509.8,104.2) and (508.08,105.91) .. (505.97,105.91) .. controls (503.85,105.91) and (502.13,104.2) .. (502.13,102.09) -- cycle ;
		%Shape: Ellipse [id:dp79446724261296] 
		\draw  [color={rgb, 255:red, 0; green, 0; blue, 255 }  ,draw opacity=1 ][fill={rgb, 255:red, 0; green, 0; blue, 255 }  ,fill opacity=1 ] (598.13,171.84) .. controls (598.13,169.73) and (599.85,168.02) .. (601.97,168.02) .. controls (604.08,168.02) and (605.8,169.73) .. (605.8,171.84) .. controls (605.8,173.95) and (604.08,175.66) .. (601.97,175.66) .. controls (599.85,175.66) and (598.13,173.95) .. (598.13,171.84) -- cycle ;
		%Shape: Ellipse [id:dp958687236272272] 
		\draw  [color={rgb, 255:red, 0; green, 0; blue, 255 }  ,draw opacity=1 ][fill={rgb, 255:red, 0; green, 0; blue, 255 }  ,fill opacity=1 ] (596.13,124.01) .. controls (596.13,121.9) and (597.85,120.19) .. (599.97,120.19) .. controls (602.08,120.19) and (603.8,121.9) .. (603.8,124.01) .. controls (603.8,126.12) and (602.08,127.83) .. (599.97,127.83) .. controls (597.85,127.83) and (596.13,126.12) .. (596.13,124.01) -- cycle ;
		%Shape: Ellipse [id:dp8524480928693141] 
		\draw  [color={rgb, 255:red, 0; green, 0; blue, 255 }  ,draw opacity=1 ][fill={rgb, 255:red, 0; green, 0; blue, 255 }  ,fill opacity=1 ] (624.13,147.92) .. controls (624.13,145.81) and (625.85,144.1) .. (627.97,144.1) .. controls (630.08,144.1) and (631.8,145.81) .. (631.8,147.92) .. controls (631.8,150.03) and (630.08,151.74) .. (627.97,151.74) .. controls (625.85,151.74) and (624.13,150.03) .. (624.13,147.92) -- cycle ;
		%Shape: Ellipse [id:dp9817247789923074] 
		\draw  [color={rgb, 255:red, 0; green, 0; blue, 255 }  ,draw opacity=1 ][fill={rgb, 255:red, 0; green, 0; blue, 255 }  ,fill opacity=1 ] (567.13,191.76) .. controls (567.13,189.65) and (568.85,187.94) .. (570.97,187.94) .. controls (573.08,187.94) and (574.8,189.65) .. (574.8,191.76) .. controls (574.8,193.87) and (573.08,195.58) .. (570.97,195.58) .. controls (568.85,195.58) and (567.13,193.87) .. (567.13,191.76) -- cycle ;
		%Shape: Ellipse [id:dp5213073787100686] 
		\draw  [color={rgb, 255:red, 0; green, 0; blue, 255 }  ,draw opacity=1 ][fill={rgb, 255:red, 0; green, 0; blue, 255 }  ,fill opacity=1 ] (561.13,108.07) .. controls (561.13,105.96) and (562.85,104.25) .. (564.97,104.25) .. controls (567.08,104.25) and (568.8,105.96) .. (568.8,108.07) .. controls (568.8,110.18) and (567.08,111.89) .. (564.97,111.89) .. controls (562.85,111.89) and (561.13,110.18) .. (561.13,108.07) -- cycle ;
		%Shape: Ellipse [id:dp11664204484675] 
		\draw  [color={rgb, 255:red, 0; green, 0; blue, 255 }  ,draw opacity=1 ][fill={rgb, 255:red, 0; green, 0; blue, 255 }  ,fill opacity=1 ] (527.13,153.9) .. controls (527.13,151.79) and (528.85,150.08) .. (530.97,150.08) .. controls (533.08,150.08) and (534.8,151.79) .. (534.8,153.9) .. controls (534.8,156.01) and (533.08,157.72) .. (530.97,157.72) .. controls (528.85,157.72) and (527.13,156.01) .. (527.13,153.9) -- cycle ;
		%Shape: Ellipse [id:dp3692934243634438] 
		\draw  [color={rgb, 255:red, 0; green, 0; blue, 255 }  ,draw opacity=1 ][fill={rgb, 255:red, 0; green, 0; blue, 255 }  ,fill opacity=1 ] (556.13,133.97) .. controls (556.13,131.87) and (557.85,130.16) .. (559.97,130.16) .. controls (562.08,130.16) and (563.8,131.87) .. (563.8,133.97) .. controls (563.8,136.08) and (562.08,137.79) .. (559.97,137.79) .. controls (557.85,137.79) and (556.13,136.08) .. (556.13,133.97) -- cycle ;
		%Shape: Rectangle [id:dp7788695712617792] 
		\draw  [fill={rgb, 255:red, 155; green, 155; blue, 155 }  ,fill opacity=1 ][line width=1.5]  (401.8,155.73) -- (410.8,155.73) -- (410.8,164.7) -- (401.8,164.7) -- cycle ;
		%Straight Lines [id:da1274363666127214] 
		\draw [line width=1.5]    (564.97,108.07) -- (518.46,126.1) -- (485.93,138.7) ;
		\draw [shift={(483.13,139.79)}, rotate = 338.81] [color={rgb, 255:red, 0; green, 0; blue, 0 }  ][line width=1.5]    (14.21,-4.28) .. controls (9.04,-1.82) and (4.3,-0.39) .. (0,0) .. controls (4.3,0.39) and (9.04,1.82) .. (14.21,4.28)   ;
		%Straight Lines [id:da40104240241374023] 
		\draw [color={rgb, 255:red, 128; green, 128; blue, 128 }  ,draw opacity=1 ][line width=1.5]    (564.97,108.07) -- (472.73,113.96) ;
		\draw [shift={(469.73,114.15)}, rotate = 356.35] [color={rgb, 255:red, 128; green, 128; blue, 128 }  ,draw opacity=1 ][line width=1.5]    (14.21,-4.28) .. controls (9.04,-1.82) and (4.3,-0.39) .. (0,0) .. controls (4.3,0.39) and (9.04,1.82) .. (14.21,4.28)   ;
		%Shape: Arc [id:dp2491950364747133] 
		\draw  [draw opacity=0] (541,117.33) .. controls (540.11,114.68) and (539.63,111.98) .. (539.51,109.31) -- (568.8,108.07) -- cycle ; \draw  [color={rgb, 255:red, 255; green, 0; blue, 0 }  ,draw opacity=1 ] (541,117.33) .. controls (540.11,114.68) and (539.63,111.98) .. (539.51,109.31) ;  
		%Shape: Ellipse [id:dp34381866520042537] 
		\draw  [color={rgb, 255:red, 0; green, 0; blue, 255 }  ,draw opacity=1 ][fill={rgb, 255:red, 0; green, 0; blue, 255 }  ,fill opacity=1 ] (576.14,153.73) .. controls (576.23,151.62) and (578.02,149.99) .. (580.14,150.09) .. controls (582.25,150.18) and (583.89,151.97) .. (583.8,154.07) .. controls (583.7,156.18) and (581.91,157.81) .. (579.79,157.72) .. controls (577.68,157.62) and (576.04,155.84) .. (576.14,153.73) -- cycle ;
		%Shape: Ellipse [id:dp624015446957777] 
		\draw  [color={rgb, 255:red, 0; green, 0; blue, 255 }  ,draw opacity=1 ][fill={rgb, 255:red, 0; green, 0; blue, 255 }  ,fill opacity=1 ] (505.73,187.78) .. controls (505.73,185.67) and (507.45,183.96) .. (509.57,183.96) .. controls (511.68,183.96) and (513.4,185.67) .. (513.4,187.78) .. controls (513.4,189.89) and (511.68,191.6) .. (509.57,191.6) .. controls (507.45,191.6) and (505.73,189.89) .. (505.73,187.78) -- cycle ;
		%Shape: Ellipse [id:dp5148262104198444] 
		\draw  [color={rgb, 255:red, 0; green, 0; blue, 255 }  ,draw opacity=1 ][fill={rgb, 255:red, 0; green, 0; blue, 255 }  ,fill opacity=1 ] (602.53,88.14) .. controls (602.53,86.03) and (604.25,84.32) .. (606.37,84.32) .. controls (608.48,84.32) and (610.2,86.03) .. (610.2,88.14) .. controls (610.2,90.25) and (608.48,91.96) .. (606.37,91.96) .. controls (604.25,91.96) and (602.53,90.25) .. (602.53,88.14) -- cycle ;
		%Shape: Ellipse [id:dp6824143765377642] 
		\draw  [color={rgb, 255:red, 0; green, 0; blue, 255 }  ,draw opacity=1 ][fill={rgb, 255:red, 0; green, 0; blue, 255 }  ,fill opacity=1 ] (544.13,59.45) .. controls (544.13,57.34) and (545.85,55.63) .. (547.97,55.63) .. controls (550.08,55.63) and (551.8,57.34) .. (551.8,59.45) .. controls (551.8,61.56) and (550.08,63.27) .. (547.97,63.27) .. controls (545.85,63.27) and (544.13,61.56) .. (544.13,59.45) -- cycle ;
		%Shape: Rectangle [id:dp8033994291082587] 
		%\draw   (336.93,38.87) -- (644.13,38.87) -- (644.13,243.67) -- (336.93,243.67) -- cycle ;
		
		% Text Node
		\draw (243.77,107.02) node [anchor=north west][inner sep=0.75pt]   [align=left] {$\displaystyle x$};
		% Text Node
		\draw (146,142.76) node [anchor=north west][inner sep=0.75pt]    {$\overline{x} -x$};
		% Text Node
		\draw (22,170.26) node [anchor=north west][inner sep=0.75pt]    {$Target\ \overline{x}$};
		% Text Node
		\draw (258,46.91) node [anchor=north west][inner sep=0.75pt]    {$\mathcal{G}[f]$};
		% Text Node
		\draw (136.4,15.2) node [anchor=north west][inner sep=0.75pt]    {$\lambda =0$};
		% Text Node
		\draw (566.97,107.21) node [anchor=north west][inner sep=0.75pt]   [align=left] {$\displaystyle x$};
		% Text Node
		\draw (469.2,142.95) node [anchor=north west][inner sep=0.75pt]    {$\overline{x} -x$};
		% Text Node
		\draw (345.2,170.45) node [anchor=north west][inner sep=0.75pt]    {$Target\ \overline{x}$};
		% Text Node
		\draw (426,87.1) node [anchor=north west][inner sep=0.75pt]    {$\mathcal{G}[f]$};
		% Text Node
		\draw (459.6,15.38) node [anchor=north west][inner sep=0.75pt]    {$\lambda =1$};
	\end{tikzpicture}
	\caption{On the left the case in which  agent $x$ remains or switches to follower status ($\lambda =0$), and on the right the case in which it remains or switches to leader status  ($\lambda = 1$).} 
	\label{fig:configurations} 
\end{figure}
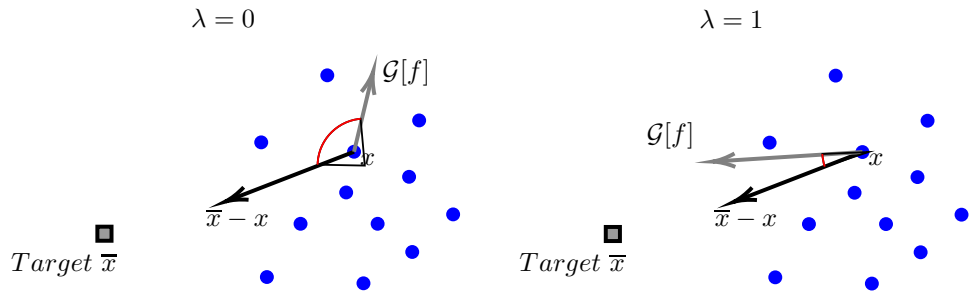

\paragraph{Density-dependent rates} Leaders emerge with higher probability where the followers density is higher and the leaders one is lower and they return to the followers status with higher probability if the followers concentration around them is lower, similarly to \cite{albi2022mean}. The transition rates read 
\begin{equation}\label{eq:rates_test}
	\begin{split}
		\pi_{L\to F} = q_F~ (1-\mathcal{D}_F[f](x,t)),\qquad
		\pi_{F\to L} = q_L~ (1-\mathcal{D}_L[f](x,t)),
	\end{split}
\end{equation} 
where $q_F$, $q_L$ are constant parameters and the functions $\mathcal{D}_F[f](x,t)$ and $\mathcal{D}_L[f](x,t)$ represent the concentration of leaders and followers in position $x$ and are defined as 
\begin{equation}\label{eq:concentration}
	\begin{split}
		&\mathcal{D}_F [f](x,t) = S_F (t) \int_{\RR^{2d}} e^{-\frac{\vert x - y\vert^2 }{\delta^2}} f_0(y,w) d(y,w),\\
		&\mathcal{D}_L [f](x,t) = S_L (t) \int_{\RR^{2d}} e^{-\frac{\vert x - y\vert^2}{\delta^2}} f_1(y,w) d(y,w),
	\end{split}
\end{equation}
with $S_F (t)$, $S_L (t)$ normalization constants to ensure that the above quantities are bounded by one and with $\delta>0$.
\paragraph{Target-oriented rates} Leaders emerge when their direction is oriented in the correct direction toward a target position, $\bar{x}$,  such as the nesting or foraging area. We consider the following rates
\begin{equation}\label{eq:rates_opt}
	\begin{aligned}
		\pi_{F\to L} &
		= \begin{cases}
			0,\qquad  \text{if } \alpha(x,v,t;f) < \overline{\alpha},\\
			1, \qquad \text{if } \alpha(x,v,t;f) \geq \overline{\alpha},
		\end{cases}\qquad
		\pi_{L\to F}&=
		\begin{cases}
			0,\qquad  \text{if } \alpha(x,v,t;f) \geq \underline{\alpha},\\
			1, \qquad \text{if } \alpha(x,v,t;f) < \underline{\alpha},
		\end{cases}
	\end{aligned}
\end{equation}
with $\underline{\alpha},\overline{\alpha}\in[-1,1]$ and
%\begin{equation}\label{eq:alpha}
%	\alpha(x,v,t;f) = \frac{\langle \bar{x}-x, G(x,v,t;f)\rangle}{\langle G(x,v,t;f),G(x,v,t;f)\rangle},
%\end{equation}
\begin{equation}\label{eq:alpha}
	\alpha(x,v,t;f) = \cos\left(\angle \left(\bar{x}-x, 	\mathcal{G}[f](x,v,t)\right)\right),
\end{equation}
with $\angle(\cdot,\cdot)$ denoting the angle between two vectors.
%\begin{equation}\label{eq:lambda_opt} 
%	\lambda = \sigma \left( \frac{\left\langle  \sum_{n=1}^d (\bar{x}-x), (G_i)_n\right\rangle }{\sum_{n=1}^d (G^2_i)_n}\right),
%\end{equation}
%where $\sigma(\alpha) = 1$ if $\alpha \geq0.5$ and $\sigma(\alpha) = 0$ if $\alpha <0.5$.
The functional $	\mathcal{G}[f](\cdot)$ accounts for the  directional information of agents according to
\begin{equation}\label{eq:G_test1}
	\mathcal{G}[f](x,v,t) =  S(v) - \mathcal{X}_c[f](x,t) - \mathcal{V}_c[f](x,v,t),
\end{equation} 
where $S(v)$ is the self-propulsion term, and the terms $\mathcal{X}_c[f](\cdot), \mathcal{V}_c[f](\cdot)$ account for the average influence induced by neighbours  as follows
\begin{equation*}
	\begin{split}	
		\mathcal{X}_c[f] (x,t) &= \int_{\mathit{B}_{r^*(x,t)}\times\mathbb{R}^d} A^{att}(x,x_*) f_\lambda(x_*,v_*,t) dx_* dv_*,\\
		\mathcal{V}_c[f] (x,v,t) &= \int_{\mathit{B}_{r^*(x,t)}\times\mathbb{R}^d} A^{ali}(v,v_*) f_\lambda(x_*,v_*,t) dx_* dv_*,
	\end{split}
\end{equation*}
with $A^{ali}(\cdot,\cdot)$, $A^{att}(\cdot,\cdot)$ defined as in equation \eqref{eq:aligment}-\eqref{eq:atraction}.
Note that in \eqref{eq:G_test1}, when the term $\mathcal{G}[f]$ is partially aligned with the target direction $\bar x - x$, i.e., $\alpha(x,v,t;f) \geq \overline{\alpha}$, agents switch to, or remain in, leader status, naturally steering their dynamics towards the target $\bar x$.
Conversely, if $\alpha(x,v,t;f) \leq \underline{\alpha}$, the agent with position and velocity $(x,v)$ remains in, or is switched to, follower status.  \ref{fig:configurations} illustrates two possible configurations.	
 
	\begin{remark}[Multiple-label case and continuous limit]\label{remark_continuos}
		We observe that the previous formulation can be extended to include multiple levels of leadership, up to a continuous space of labels \cite{cristiani2023kinmacr, albi2019leader}.
		Hence, we consider $\lambda \in \{\lambda_1,\ldots,\lambda_{N_\ell}\}$ such that $\lambda_k = k/N_\ell$ with $\lambda_1 = 0$ and $\lambda_{{N}_\ell} =1$, and we assume that the evolution of the label density $f_k(t) := f(x,v,\lambda_k,t)$ in \eqref{eq:boltz_lin} is determined according to the following master equation
		\begin{equation}\label{eq:master_multiple}
			\frac{d f_k(t)}{dt} = \sum_{j\neq k} \left( \pi_{jk}(t) f_j(t) - \pi_{kj}(t) f_k(t)\right),\quad k=1,\ldots,{N}_\ell.
		\end{equation}
		%where  we neglect the transport and interaction terms in \eqref{eq:boltz_lin}.
		%and for the boundary values $\lambda_1=0, \lambda_{{N}_\ell} = 1$ we have
		%\begin{equation}\label{eq:master_multiple_bc}
		%\begin{split}
		%&\frac{d f_1(t)}{dt} =  \Pi_{21}[f](t) f_2(t) - \Pi_{12}[f](t) f_1(t) ,\\
		%&\frac{d f_{N_\ell}(t)}{dt} =  \Pi_{N_\ell-1 N_\ell}[f](t) f_{\N_\ell-1}(t) - \Pi_{N_\ell N_\ell-1}[f](t) f_{N_\ell}(t) ,
		%\end{split}
		%	\end{equation}
		The non-negative transition rates  are defined as $\pi_{kj}(t) := {\pi_{\lambda_{k}\rightarrow\lambda_{j}}}(x,v,t;f)$ representing  the possibility to jump from label $\lambda_j$ to label $\lambda_k$. 
		We can introduce formally a continuous approximation of \eqref{eq:master_multiple} by scaling the density and the transition rates according to
		\begin{equation*}
			{f_k(t)} =	 f(\lambda_k,t)/N_\ell, \qquad
			{\pi}_{kj}(t) = {\omega}(\lambda_j|\lambda_k)/N_\ell,
		\end{equation*}
		where ${\omega}(\lambda_j|\lambda_k)$ are the continuous transition rate. Then,  for $N_\ell \to \infty$ we retrieve the master equation as
		\begin{equation}\label{eq:master_continuos} 
			{\partial_t{f(\lambda,t)}} = \int_{0}^1 \left( {\omega}(\lambda|\lambda_*) f(\lambda_*,t) - {\omega}(\lambda_*|\lambda) f(\lambda,t)\right) d\lambda_*, \quad \lambda \in [0,1],
		\end{equation}
		which encodes arbitrary jumps between different labels $\lambda\in[0,1]$. Equivalently, the master equation \eqref{eq:master_continuos}  can be written in terms of the Kramer-Moyal expansion as follows
		\begin{equation}\label{eq:KM_master} 
			{\partial_t{f(\lambda,t)}} = \sum_{m=1}^{+\infty}\frac{(-1)^m}{m!}\partial^{(m)}_\lambda \left( \mathcal{W}_m(\lambda) f(\lambda,t) \right), 
		\end{equation}  
		where the coefficient are defined as the moment functions
		\begin{equation}\label{eq:FP_master_coef}
			\mathcal{W}_m(\lambda,t) = \int_{0}^1 (\lambda_*-\lambda)^m {\omega}(\lambda_*|\lambda) d \lambda_*,   \qquad  m = 1,2,\ldots.
		\end{equation}
		From this description, assuming that small jumps occurs, i.e. the transition rates are localized over the label space, it is possible to show that the first term, and possibly the second order term of \eqref{eq:KM_master}  characterize exactly the microscopic label dynamics, see \cite{pawula1967approximation,van1992stochastic}. Thus we retrieve a non-linear Fokker-Planck opertor for the evolution of the continuous label dynamics.
		Finally we remark that, the type of transition rates $\pi_{kj}(t)$ in the Markov process \eqref{eq:master_multiple} is crucial to correctly treat the limit as $N_\ell$ becomes large and to derive further approximation such as Fokker-Planck in differential form, see \cite{ethier2009markov,van1992stochastic}. Furthermore, the Fokker-Planck approximation may fail to describe the correct  discrete dynamics, in particular when the transition rates have non-linear behaviors, see for example  \cite{doering2005extinction,lunz2021continuum}. 
	\end{remark}

\subsection{Grazing-collision limit and mean-field model} \label{sec:kinetic_model_1}
In order to retrieve asymptotic behaviour of the Boltzmann-type equation \eqref{eq:boltz_lin}, we  resort on a mean-field approximation of the interaction dynamics. Thus, we introduce a grazing collision limit for the interaction operator \eqref{eq:Bo}, following the approach in \cite{pareschi2013interacting,carrillo2010particle}. To this aim, we rescale the interaction frequency $\eta$ and the interaction propensity $\alpha$  to maintain asymptotically the memory of the microscopic interactions, as follows
\begin{equation}\label{eq:scaling}
	\alpha = \varepsilon, \qquad \eta = \frac{1}{\varepsilon},
\end{equation}
for $\varepsilon>0$,
which corresponds to the case where the interaction kernel concentrates on
binary interactions producing very small changes in the agents velocity but at the same
time the number of interactions becomes very large.
From now on, for simplicity we remove the dependence on time $t$. 
We introduce the test function $\psi(x,v)\in C^1_0(\RR^d\times\RR^d)$ and  we write the weak form of the scaled kinetic equation \eqref{eq:boltz_lin} %collision operator $Q(\cdot,\cdot)$ in weak form as follows
\begin{equation}
	\begin{split}\label{eq:boltz_weak_0}
		&\int_{\RR^{2d}} \left( \partial_t  f_\lambda(x,v) + v\cdot \nabla_x  f_\lambda(x,v)\right) \psi(x,v)d(x,v)- \int_{\RR^{2d}} \T_\lambda[ f](x,v)\psi(x,v)d(x,v) = \\
		&\quad\frac{1}{\varepsilon}\sum_{\lambda_*} \int_{\RR^{2d}}\int_{\Omega} \left(\psi(x,v')-\psi(x,v)\right)f_{\lambda_*}(x_*,v_*)  f_\lambda(x,v) d(x,v)d(x_*,v_*),
	\end{split}
\end{equation}
%\begin{equation}\begin{split}\label{eq:collisional_op}
%		\int_{\RR^{2d}} & Q(f,f)(x,v,\lambda)\psi(x,v)d(x,v) =\eta \sum_{\lambda_*\in\{0,1\}} \int_{\RR^{2d}}\int_{\mathit{B}_{r_*(x,t)}\times \RR^{d} }\ \left(\psi(x,v')-\psi(x,v)\right) df_* df ,
%\end{split}\end{equation}
%The scaled equation \cref{eq:boltz_weak_0} reads as
%\begin{equation}\label{eq:boltz_weak_scaled} 
%	\begin{split}
%		\int_{\RR^{2d}}& \left( \partial_t f(x,v,\lambda) + v\cdot \nabla_x f(x,v,\lambda)\right) \psi(x,v)d(x,v)- \int_{\RR^{2d}} \T[ f](x,v,\lambda)\psi(x,v)d(x,v) = \\
%		&\frac{1}{\varepsilon}\sum_{\lambda_*\in\{0,1\}} \int_{\RR^{2d}}\int_{\mathit{B}_{r_*(t,x)}\times \RR^{d} } \left(\psi(x,v')-\psi(x,v)\right)df_* df ,
%	\end{split}
%\end{equation}
with scaled  interactions \eqref{eq:binary} as follows
\begin{equation}\label{eq:binary_scaled}
	v'-v =  \varepsilon \mathcal{F}_\lambda(x,x_*,v,v_*).
\end{equation}
Since as $\varepsilon \to 0$, we have $v'\to v$ we can expand $\psi(x,v')$ in Taylor series centred in $(x,v)$ up to second order and rewrite the right hand side of equation \eqref{eq:boltz_weak_0} as
\begin{equation}\label{eq:boltz_weak_1}
	\begin{split}
		\frac{1}{\varepsilon} &\sum_{\lambda_*}\int_{\RR^{2d}}\int_{\Omega }  \left(\psi(x,v')-\psi(x,v)\right) df_{\lambda_*} df_{\lambda} =\\&\frac{1}{\varepsilon} \sum_{\lambda_*} \int_{\RR^{2d}} \int_{\Omega} \nabla_v \psi (x,v) \cdot (v'-v) df_* df+ R(\varepsilon),
	\end{split}
\end{equation} 
where we used the shorten notation $df_{\lambda_*}=f(x_*,v_*,\lambda_*)d(x_*,v_*)$,\\  $df_\lambda=f(x,v,\lambda)d(x,v)$, and where $R(\varepsilon)$ indicates the remainder which is given by 
\begin{equation}\label{eq:remainder}
	R(\varepsilon) = \frac{1}{2\varepsilon}   \sum_{\lambda_*} \int_{\RR^{2d}} \int_{\Omega } \left[ \sum_{i,j = 1}^d \partial_v^{(i,j)}\psi (x,\bar{v})(v'-v)_i(v'-v)_j \right]  df_{\lambda_*} df_{\lambda},
\end{equation}
with 
\begin{equation*}
	\bar{v}= \gamma v + (1-\gamma)v',
\end{equation*}
for some $\gamma \in [0,1]$.
Therefore, the scaled binary interaction term \eqref{eq:boltz_weak_1} reads 
\begin{equation}\label{eq:boltz_weak_1_scaled}
	\begin{split}
		\sum_{\lambda_*} \int_{\RR^{2d}} \int_{\Omega } \nabla_v \psi (x,v) \cdot \mathcal{F}_\lambda(x,x_*,v,v_*)  df_{\lambda_*} df_{\lambda} + R(\varepsilon).
	\end{split}
\end{equation} 
Integrating equation \eqref{eq:boltz_weak_1_scaled} by parts and taking the limit $\varepsilon \to 0$ we have
\begin{equation}\label{eq:boltz_weak_2}
	\begin{split}
		& \sum_{\lambda_*} \int_{\RR^{2d}}\int_{\Omega } \nabla_v 
		\psi (x,v) \cdot \mathcal{F}_\lambda(x,x_*,v,v_*)   df_{\lambda_*} df_{\lambda}=\\&-
		\sum_{\lambda_*}\left\langle  \nabla_v \cdot \Bigl[f_\lambda(x,v) \int_{\Omega }   \mathcal{F}_\lambda(x,x_*,v,v_*)    df_{\lambda_*}\Bigr], \psi(x,v)\right\rangle, 
	\end{split}
\end{equation}
where  we denoted the inner scalar product
\begin{equation}\label{eq:scalar_product}
	\left\langle h, \phi \right\rangle : = \int_{\RR^{2d}} h(x,v) \phi(x,v) d(x,v), 
\end{equation}
for any function $h(x,v)$, $\phi(x,v)$ for which the integral in \eqref{eq:scalar_product} is well defined. By similar arguments of \cite{albi2016invisible}, it can be shown rigorously that $R(\varepsilon) \to 0$, as $\varepsilon \to 0$. 
Thus, we can rewrite equation \eqref{eq:boltz_weak_0} as follows
\begin{equation}\begin{split}\label{eq:boltz_weak}
		&\left\langle \partial_{t} f_\lambda(x,v) + v \cdot \nabla_x f_\lambda(x,v) - \T_\lambda[ f](x,v),\psi(x,v)\right\rangle  = \\
		&	 -\sum_{\lambda_*}\left\langle  \nabla_v \cdot \Bigl[f_\lambda(x,v) \int_{\Omega }   \mathcal{F}_\lambda(x,x_*,v,v_*)   df_*\Bigr], \psi(x,v)\right\rangle.
\end{split}\end{equation}
Finally, we retrieve the mean-field equation as the strong form of \eqref{eq:boltz_weak} 
\begin{equation}
	\begin{aligned}\label{eq:boltz_strong}
		& \partial_{t} f_\lambda(x,v) + v \cdot \nabla_x f_\lambda(x,v) - \T_\lambda[f](x,v) =\\
		&\qquad\quad	 - \nabla_v \cdot \Bigl[f_\lambda(x,v) \int_{\Omega }   \mathcal{F}_\lambda(x,x_*,v,v_*)   \sum_{\lambda_*} f_{\lambda_*}(x_*,v_*)d(x_*,v_*)\Bigr].
	\end{aligned}
\end{equation}
Summing over the values of $\lambda$ in equation \eqref{eq:boltz_strong} the transition operator vanishes as in \eqref{eq:tnc} and we obtain the mean-field model for the total density $g(x,v)$ as 
\begin{equation}
	\begin{split}\label{eq:boltz_strong_g}
		\partial_{t} g(x,v) + &v \cdot \nabla_x g(x,v) = \\
		&- \nabla_v \cdot \Bigl[\sum_{\lambda} f_\lambda(x,v)\int_{\Omega }  \mathcal{F}_\lambda(x,x_*,v,v_*)  g(x_*,v_*) d(x_*,v_*)\Bigr].
	\end{split}
\end{equation}

\begin{remark}
	Note that the continuous mean-field model \eqref{eq:boltz_strong} and the microscopic one \eqref{eq:dynamics} are equivalent
	when we consider the empirical distribution of the $N$-particles
	\begin{equation}\label{eq:empirical_dist}
		f^N(x,v,\lambda,t) = \frac{1}{N} \sum_{i=1}^N \delta(x-x_i(t))\delta(v-v_i(t))\delta(\lambda-\lambda_i(t)),
	\end{equation}
	where $\delta(\cdot)$ indicates the Dirac-delta function.
\end{remark} 
\begin{remark}\label{remark_kincont}
		Here we model interactions and label transition in such a way that they possibly occur on different time scales. Thus, when the grazing collision scaling is considered, it is applied solely to the interaction term, while the operator $\mathcal{T}_\lambda[f](x,v)$ is treated separately for the evolution of discrete labels. Following the discussion of  \ref{remark_continuos}, one could consider continuous label evolution introducing different operator for transition dynamics, similarly to what proposed in \cite{albi2022mean, cristiani2023kinmacr,loy2020non}. In this framework, the equation \eqref{eq:boltz_lin} modifies as follows
		\begin{align}
			\label{eq:new}
			\partial_t f(x,v,\lambda,t) +v\cdot \nabla_v f(x,v,\lambda,t)  = \mathcal{Q}(f,f)(x,v,\lambda,t) + \mathcal{P}(f)(x,v,\lambda,t),
		\end{align}
		for $\lambda\in [0,1]$, and  $\mathcal{P}(\cdot)$ the operator accounting for the continuous label dynamics.
		% The type of interaction encoded in  $\mathcal{P}(\cdot)$ can be described in various ways, see for example \cite{loy2020non,cristiani2023kinmacr}. 
		To illustrate this further, we neglect in \eqref{eq:new} the transport and interaction terms, and we implicitly assume state dependence  on $(x,v,t)$, thus we introduce the following microscopic dynamics
		\begin{equation}\label{eq:micro_lambda}
			\lambda' = \lambda + \beta\Theta (\lambda_*-\lambda)(1+\chi),
		\end{equation}
		where $\beta\in(0,1]$ represents the size of the jump, $\lambda_*$ is drawn uniformly in the continuous label space $\Lambda =[0,1]$, $\Theta$ is a Bernoulli random variable such that $ \Theta\sim \text{Bernoulli} (p)$ with probability $p \equiv p(\lambda_*|\lambda)$, and $\chi$ an independent random variable of zero average, second moment $\mathbb{E}[\chi^2] = \varsigma^2$  and finite higher moments. The dynamics \eqref{eq:micro_lambda} describes the possibility of jumping from state $\lambda$ towards state $\lambda_*$, and allowing random perturbation due to $\chi$. In particular, for $\beta =1$ and $\chi\equiv 0$  if a jump occurs then $\lambda'=\lambda_*$, instead if $\chi\not\equiv 0$ and the support of $\chi$ is such that $\chi\in[-1,1/\beta -1]$ the interaction \eqref{eq:micro_lambda} preserves the bounds, i.e. $\lambda'\in[0,1]$.
		Thus, at mesoscopic level, the evolution of the density of labels in weak form  writes
		\begin{equation}\label{eq:weak_boltz_new}
			\frac d{dt}	\int_0^1 f(\lambda)\phi(\lambda) d \lambda = \nu\int_{0}^1 \left(\int_0^1 \mathbb{E}\left[ \phi(\lambda')- \phi(\lambda)\right]d\lambda_*\right) f(\lambda)d{\lambda}, 
		\end{equation}   
		where $\phi(\lambda)\in \mathcal{C}_0^\infty([0,1])$ is a test function, $\nu>0$ is the frequency rate, and $\mathbb{E}[\cdot]$ is the expected value associated to the random variables $\Theta$ and $\chi$.  
		Hence,  if the random perturbation is neglected ($\chi\equiv 0$), the series expansion of $\phi(\lambda')$ around $\lambda$ writes
		$$
		\mathbb{E}[\phi(\lambda')-\phi(\lambda)]=\sum_{m=1}^{+\infty}\frac{\beta^m p(\lambda_*|\lambda)}{m!} (\lambda_*-\lambda)^m\partial^{(m)}_\lambda \phi(\lambda),
		$$%+\frac{\beta^2}{2}\Theta^2(\lambda_*-\lambda)^2\partial^2_\lambda \phi(\lambda)+\frac{\beta^3}{6}\Theta^3(\lambda_*-\lambda)^3\partial^3_\lambda \phi(\bar\lambda),$$
		%where $\bar\lambda =  \theta\lambda  +(1-\theta)\lambda'$ with $\theta \in [0,1]$. 
		where we employ the property that  $\mathbb{E}[\Theta^m] = p$ for any $m>0$. If we return \eqref{eq:weak_boltz_new} to the strong form, we obtain
		\begin{equation}\label{eq:KM_master_2} 
			{\partial_t{f(\lambda,t)}} = \sum_{m=1}^{+\infty}\frac{(-1)^m}{m!}\partial^{(m)}_\lambda \left(\int_0^1\nu\beta^m p(\lambda_*|\lambda) (\lambda_*-\lambda)^m f(\lambda,t) \right), 
		\end{equation}  
		where consistency with  expansion in \eqref{eq:KM_master} is retrieved for $\beta = 1$, and transition rates such that $\omega(\lambda_*|\lambda)\equiv \nu p(\lambda_*|\lambda)$.
		Alternatively, if $\chi\not\equiv 0$, we can consider the expansion up to third order of $\phi(\lambda')$ around $\lambda$, and substitute it in \eqref{eq:weak_boltz_new} to obtain
		\begin{equation}\label{eq:weak_boltz_taylor}
			\begin{split}
				&\frac d{dt}	\int_0^1 f(\lambda)\phi(\lambda) d \lambda = \nu\beta \int_{0}^1 \left(\int_0^1 p(\lambda_*|\lambda)(\lambda_*-\lambda)\partial_\lambda \phi(\lambda)d\lambda_*\right) f(\lambda)d{\lambda}\cr
				&\qquad \frac{\nu\beta^2}{2}(1+\varsigma^2) \int_{0}^1 \left(\int_0^1 p(\lambda_*|\lambda)(\lambda_*-\lambda)^2\partial^2_\lambda \phi(\lambda)d\lambda_*\right) f(\lambda)d{\lambda}
				+\mathcal{R}(\nu\beta^2) ,
			\end{split}
		\end{equation}  
		where the reminder is defined as 
		\begin{equation}\label{eq:reminder}
			\begin{split}
				\mathcal{R}(\nu\beta^2) = \frac{\nu\beta^3}{3}(1+3\varsigma^2+\mu) \int_{0}^1 \left(\int_0^1 p(\lambda_*|\lambda)(\lambda_*-\lambda)^3\partial^2_\lambda \phi(\bar \lambda)d\lambda_*\right) f(\lambda)d{\lambda} ,
			\end{split}
		\end{equation}  
		%\begin{equation}
		%\begin{split}
		%\mathbb{E}[\phi(\lambda')-\phi(\lambda)] &= \beta p(\lambda_*|\lambda)(\lambda_*-\lambda)\partial_\lambda \phi(\lambda) + \frac{\beta^2}{2}(1+\varsigma^2)p(\lambda_*|\lambda)(\lambda_*-\lambda)^2\partial^2_\lambda \phi(\lambda)
		% \cr&\qquad \frac{\beta^3}{3}(1+3\varsigma^2+\mu)p(\lambda_*|\lambda)(\lambda_*-\lambda)^3\partial^3_\lambda \phi(\bar\lambda)
		%% + \frac{\beta}{2} p(\lambda_*|\lambda)(\lambda_*-\lambda)^2\partial^2_\lambda \phi(\lambda)
		%\end{split}
		%\end{equation}
		with $\bar\lambda = \tau \lambda'+ (1-\tau)\lambda$ for $\tau\in[0,1]$, and $\mu:= \mathbb{E}[\chi^3]$. Finally, introducing the grazing scaling in \eqref{eq:weak_boltz_taylor} as  
		$\nu=\bar \nu/\varepsilon$, $\beta = 1/\varepsilon$, and $\varsigma=1/\sqrt{\varepsilon}$,  in limit for $\varepsilon\to0$ we have that both first and second order terms of  \eqref{eq:weak_boltz_taylor} are preserved since $\nu\beta\to \bar \nu$ and $\nu\beta^2(1+\varsigma^2)\to \bar \nu$, whereas the reminder $\mathcal{R}(\varepsilon)\to 0$ vanishes. Passing to the strong form the limit of \eqref{eq:weak_boltz_taylor} restitutes the following Fokker-Planck equation 		
		\begin{equation}\label{eq:FP_boltz_1} 
			\partial_t f  = -{\partial_\lambda} \left(  f\int_0^1 (\lambda_*-\lambda) \bar\nu p(\lambda_*|\lambda) d \lambda_* \right) +\frac{1}{2}{\partial^2_\lambda} \left(  f \int_0^1 (\lambda_*-\lambda)^2 \bar \nu p(\lambda_*|\lambda) d \lambda_* \right).
			%= -{\partial_\lambda} \left(  f(\lambda,t)\mathcal{W}_1(\lambda,t)\right),
		\end{equation}  
		We note that if the perturbation $\chi$ is neglected in \eqref{eq:lambda_evolution}, i.e. $\varsigma$ and $\mu$ are set to zero in \eqref{eq:weak_boltz_taylor}, then in the grazing limit the diffusion term disappears, and the label dynamics in \eqref{eq:FP_boltz_1} are governed solely by the drift term. 
	\end{remark}

%%%%%%%%%%%%%%%%%%%%%%%%%%%%%%%%%%%%%%%%%%
%%%%%%%%%%%%%%%%%%%%%%%%%%%%%%%%%%%%%%%%%
%%%%%%%%%%%%%%%%%%%%%%%%%%%%%%%%%%%%%%%%
\section{Stochastic particle-based approximation}\label{sec:numerical_methods}
We aim at solving the large system of agent \eqref{eq:dynamics} for $N\gg1$, i.e. at solving  the mean-field model  \eqref{eq:boltz_strong}  by means of the scaled Boltzmann equation in the asymptotic regime \eqref{eq:scaling}.
In particular, we aim at developing  asymptotic stochastic algorithms for the simulation of the swarming dynamics, such as in \cite{albi2013binary, pareschi2013interacting}. These approaches, based on Monte-Carlo algorithms are based of direct simulation Monte-Carlo methods (DSMCs) for kinetic equations \cite{nanbu1986theoretical,pareschi2001time}. We mention also Random Batch Methods (RBMs) which, similarly, have been devised for simulating large systems of interacting agents \cite{jin2020random}.

\subsection{Asymptotic Nanbu-type algorithm}
In order to solve the mean-field dynamics we consider the Boltzmann-type equation \eqref{eq:boltz_lin} in the scaling limit \eqref{eq:scaling}, and we split the dynamics evaluating in three different steps, the free transport, the label evolution and the interaction process, as follows
\begin{align}
	\partial_t f_\lambda(x,v) & = -v\cdot \nabla_x f_\lambda(x,v) \label{eq:boltz_transport}\\
	\partial_t f_\lambda(x,v) & = \mathcal{T}_\lambda[f](x,v)\label{eq:boltz_label}\\
	\partial_t f_\lambda(x,v) & = Q^{\varepsilon}_\lambda(f_\lambda,f_\lambda)(x,v).\label{eq:boltz_collision}
\end{align}
%then solving the interaction step 
%\begin{equation}\label{eq:boltz_collision}
%	\partial_t f = Q(f,f),
%\end{equation}
%
%	\begin{equation} \label{eq:boltz_transport}
%		\partial_t f = -v\cdot \nabla_x f,
%	\end{equation}
%then solving the interaction step 
%\begin{equation}\label{eq:boltz_collision}
%	\partial_t f = Q(f,f),
%\end{equation}
In order to approximate the time evolution of the density $f_\lambda(x,v,t)$ we assume to sample $N_s$ particles $(x_i^0,v_i^0,\lambda_i^0)$ from the initial distribution. We consider a time interval $[0, T]$ discretized in $N_t$ intervals of size $\Delta t$.

\paragraph{Transport step}
First, we focus on the transport step in equation \eqref{eq:boltz_transport} and  we employ an explicit Euler scheme to compute the solution at time $t^{n+1}$ as follows 
\begin{equation}\label{eq:transport}
	x_i^{n+1} = x_i^{n} + \Delta t v_i^n, \qquad i = 1,\ldots,N_s
\end{equation}
\paragraph{Labels switching}
Secondly, we simulate how the labels change denoting by $f_\lambda^{n}$ the approximation of $f_\lambda(x, v,n\Delta t)$, and writing the discrete version of equation \eqref{eq:boltz_label}, for the transition operator \eqref{eq:master_0} as follows 
\begin{equation}\label{eq:lambda_evolution}
	\begin{split}
		f_0^{n+1} = (1-\Delta t ~\pi_{F\to L}^n)~ f_0^n+ \Delta t ~\pi_{L\to F}^n ~f_1^n,\\
		f_1^{n+1} = (1-\Delta t ~\pi_{L\to F}^n)~ f_1^n + \Delta t ~\pi_{F\to L}^n~ f_0^n,
	\end{split}
\end{equation}
%where $\pi_{F\to L}(\cdot)$ and $\pi_{L\to F}(\cdot)$ are the transition rates. then we 
The following  \eqref{alg_lambda} describes how to simulate equation \eqref{eq:lambda_evolution} in a time interval $[0,T]$ divided into $N_t$ time steps. 

\begin{algorithm}
	\caption{Labels switching}
	\label{alg_lambda}
	\begin{algorithmic}[1]
		\STATE  Given $N_s$ samples $(x_i^0,v_i^0,\lambda_i^0)$ from the initial distribution $f_\lambda^0$; 
		\FOR {$n=0$ \texttt{to} $N_t$} 
		\FOR{ $i=1$ \texttt{to}  $N_s$}
		\STATE compute $
		p_{L} =\Delta t ~\pi_{F\to L}$, and	$p_{F}= \Delta t ~\pi_{L\to F}$, 
		\IF {$\lambda_i^n = 0$} \STATE with probability $p_{L}$ agent $i$ becomes a leader: $\lambda_i^{n+1} = 1$,
		\ENDIF
		\IF {$\lambda_i^n = 1$} \STATE with probability $p_{F}$ agent $i$ becomes a follower: $\lambda_i^{n+1} = 0$,
		\ENDIF
		\ENDFOR
		\ENDFOR
	\end{algorithmic}
\end{algorithm}
\paragraph{Interaction step} Finally, we consider the interaction step \eqref{eq:boltz_collision} decomposing the interaction operator \eqref{eq:Bo} in its gain and loss part,
\[
Q^{\varepsilon}_\lambda(f_\lambda,f_\lambda) = \frac{1}{\varepsilon} \left[Q^{\varepsilon,+}_\lambda(f_\lambda,f_\lambda)- \rho^* f_\lambda \right],
\]
where $\rho^*=M/N$ is the topological mass. Considering a forward discretization we obtain
\begin{equation}\label{eq:collision}
	f_\lambda^{n+1} =\left( 1-\frac{\rho^*\Delta t}{\varepsilon} \right)  f_\lambda^n + \frac{\rho^*\Delta t}{\varepsilon}  \frac{Q^{\varepsilon,+}(f_\lambda^n,f_\lambda^n )}{\rho^*}.
\end{equation} 
Equation \eqref{eq:collision} can be interpreted as follows. With probability $1-\rho^*\Delta t/\varepsilon$ an individual in position $x$, velocity $v$ and label $\lambda$ will not interact with other individuals and, with probability $\rho^*\Delta t /\varepsilon$, it will interact with another individual according to
\begin{equation}\label{eq:binary_interaction}
	v^{n+1}_i = v^{n}_i + \varepsilon \mathcal{F}_{\lambda_i^n}(x^n_i,x^n_j,v_i^n,v_j^n),
\end{equation}
for any $i=1,\ldots,N_s$, and where $(x^n_j,v_j^n)$ is selected randomly among the nearest neighbours belonging to the topological ball $\mathit{B}_{r^*}(x_i,t)$. We will assume $\rho^*\Delta t = \varepsilon$ to maximize the total number of interactions and ensure that at each time step all agents interact with another individual with probability one.

Note that the sampling procedure of agents from the topological ball  $\mathit{B}_{r^*}(x_i,t)$ can have extremely high computational costs, especially when the sample size is large, since it requires the explicit computation of the distances between each agent $i$ and all the others agents.  In order to improve the computational efficiency of this step different algorithms can be applied. Among them we recall the $k$-Nearest Neighbours search. With this algorithm, the computational costs are reduced from quadratic to logarithmic. We refer to \cite{friedman1977algorithm} for further details. The main novelty introduced in this paper is to consider an improved version of the $k$-NN search algorithm. In particular, in order to approximate the topological ball, we consider a subsample of size $N_c\ll N_s$, where $N_s$ is the total number of particles. We define the approximation to radius of the topological ball as follows 
	\begin{equation}\label{eq:topological_ball}
		\tilde r_*(x_i,t) = \arg\min_{r>0}\left\{\frac{1}{N_c} \sum_{k=1}^{N_c}\chi_{{B}_{r}(x_i)}(x_k)\geq \mathcal \rho^* \right\},
	\end{equation}    
	where $\rho^*$ is the target topological mass. As in the classical algorithm, we perform a $k$-NN search over a $k-d$ binary tree, but, instead of constructing the binary tree over the whole sample, we built it over a subsample. Again we are able to partition the space and organize the points optimally dividing them according to their medians.  Then, we use a $k$-NN algorithm over the binary tree to find the $\rho^*N_c$ nearest neighbours to a given agent $i$, using the tree structure. This strategy provides an estimate of which are the nearest agents but leads to an improvement in terms of efficiency, as we will show in the next subsection.  
%In order to improve the computational efficiency of this step we propose a procedure based on two steps: $a)$ an approximation of the topological ball, $b)$ $k$--Nearest Neighbours ($k$--NN) search.
%
%{\em $a)$ Topological ball approximation.} 
%\textcolor{blue}{The main novelty of our approach consists in approximating the topological ball by consider a subsample of size $N_c$ of the $N_s$ selected particles such that $N_c<<N_s$}. To this aim, we define the approximation to radius of the topological ball as follows
%\begin{equation}\label{eq:topological_ball}
%	\tilde r_*(x_i,t) = \arg\min_{r>0}\left\{\frac{1}{N_c} \sum_{k=1}^{N_c}\chi_{{B}_{r}(x_i)}(x_k)\geq \mathcal \rho^* \right\},
%\end{equation}    
%where $\rho^*$ is the target topological mass. 
%
%{\em $b)$ \textit{$k$--NN} search.}  We perform a \textit{ $k$--NN search} over a $k$-$d$ binary tree.  First, we construct the binary tree on the subsample in such a way to partition the space and organize the points optimally dividing them according to their medians. We assume that every leaf-node contains at most $N_l$ points. Then, we use a $k$-NN algorithm to find the $\rho^*N_c$ nearest neighbours to a given agent $i$, using the tree structure.  
%We will show in the numerical experiments that this algorithm reduces the computational costs from the original quadratic to logarithmic. 
%We refer to \cite{friedman1977algorithm} for further details about this procedure. 

Algorithm \ref{alg_binary} describes how to solve equation \eqref{eq:collision} in a time interval $[0,T]$ divided into $N_t$ time steps. 
\begin{algorithm}
	\caption{Asymptotic Nanbu algorithm}
	\label{alg_binary}
	\begin{algorithmic}[1]
		\STATE Give  $N_s$ samples $(x_i^0,v_i^0,\lambda_i^0)$ from the initial distribution $f_\lambda^0$; 
		\STATE Set the value of the topological mass $\rho^*$ and of the subsample size $N_c$;
		\FOR {$n=0$ \texttt{to} $N_t$}
		\STATE Select a subsample of size $N_c$;
		\STATE Construct a binary tree over the subsample;\FOR {$i=1$ \texttt{to}  $N_s$}
		\STATE Find the $\rho^*N_c$ nearest agents using a $k$--NN search algorithm on the tree;
		\STATE Select randomly an index $j$ among the   nearest neighbours;
		\STATE Compute the velocity change $v_i^{n+1}$ as in equation \eqref{eq:binary_interaction};
		\STATE Update the position $x_i$ according to \eqref{eq:transport}, with $\rho^*\Delta t = \varepsilon$.
		\ENDFOR
		\ENDFOR
	\end{algorithmic}
\end{algorithm}

%\begin{remark}\label{rmk:cost}
%	Note that the main advantage of Algorithm \cref{alg_binary} is that it allows us to construct the topological ball $\mathcal{B}_{r_M^*(t,x_i)} $ reducing the computational cost from the original $\mathcal{O}(N_s^2)$ to $\mathcal{O}(M N_s \log{N_s})$. Indeed, in general to find the nearest neighbors to a given agent we should first compute the distances between the selected agent and all the others, then sort the results and finally identify the index of the $M$ nearest agents. We will call this procedure \textit{exhaustive search}. In Algorithm \cref{alg_binary} instead we perform a $k$-nn search on a binary tree. The computational cost required to construct the binary tree over the subsample of size $N_c$ is $\mathcal{O}(N_s \log{N_s})$ and the one needed to search for each agent its $M$ nearest neighbor is $\mathcal{O}(N_s M \log{N_s})$. Hence the total cost is  $\mathcal{O}(N_s M \log{N_s})$. For further details we refer to \cite{friedman1977algorithm}.
%\end{remark}
%%%%%%%%%%%%%%%%%%%%%%%%%%%%%%%%%%%%%%%
%%%%%%%%%%%%%%%%%%%%%%%%%%%%%%%%%%%%%%
%%%%%%%%%%%%%%%%%%%%%%%%%%%%%%%%%%%
%%%%%%%%%%%%%%%%%%%%%%%%%%%%%%%%%%%%%%%%%%%%%%%%%%%%%%%%%%%%
%%%%%%%%%%%%%%%%%%%%%%%%%%%%%%%%%%%%%%%%%%%%%%%%%%%%%%%%%%%%%%
%%%%%%%%%%%%%%%%%%%%%%%%%%%%%%%%%%%%%%%%%%%%%%%%%%%%%%%%%%%%%%
%\section{Numerical experiments} \label{sec:numerical_experiments} 

\subsection{Numerical validation}\label{sec:validation} 
\begin{figure}[tbhp]
	\centering
	\includegraphics[width=0.327\linewidth]{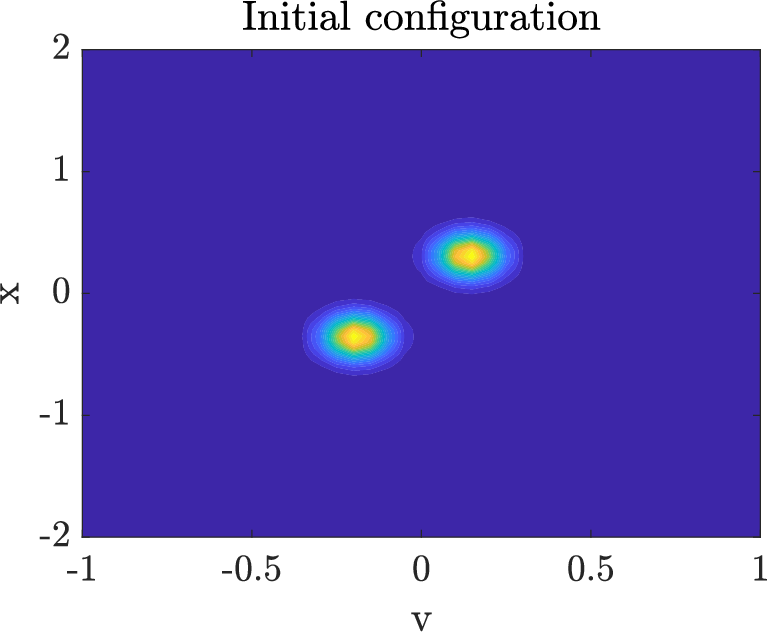}
	\includegraphics[width=0.327\linewidth]{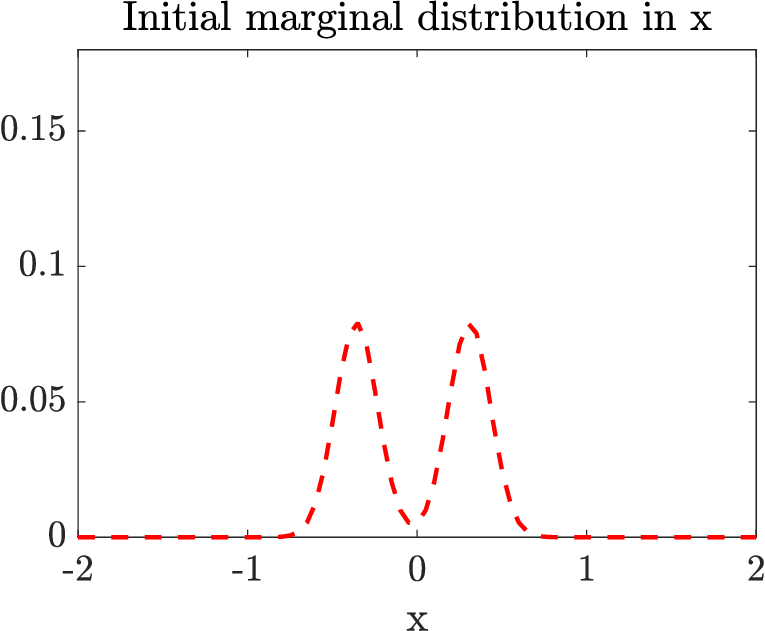}
	\includegraphics[width=0.327\linewidth]{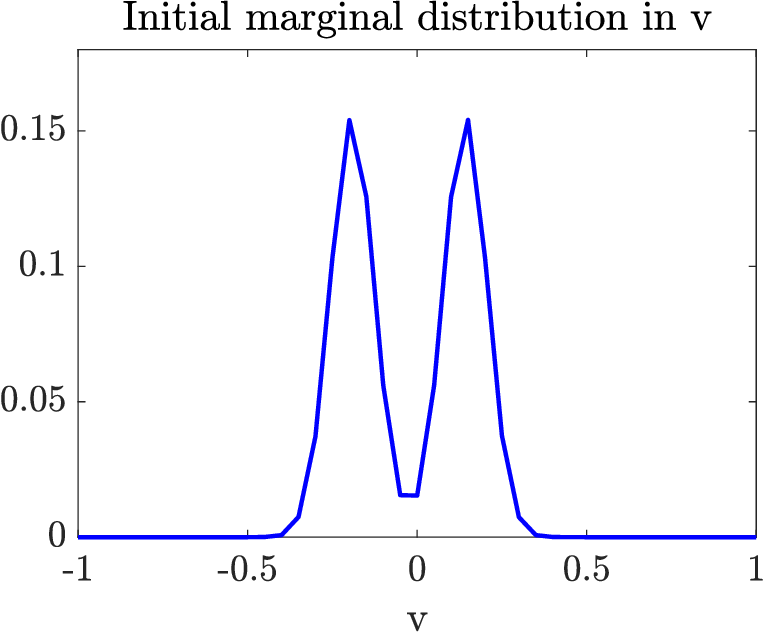} 
	\caption{Validation test: initial configuration and its marginals in $x$ and $v$.}
	\label{fig:validation_ord2_ali_IC} 
\end{figure}

\begin{figure}[tbhp]
	\centering
	\includegraphics[width=0.327\linewidth]{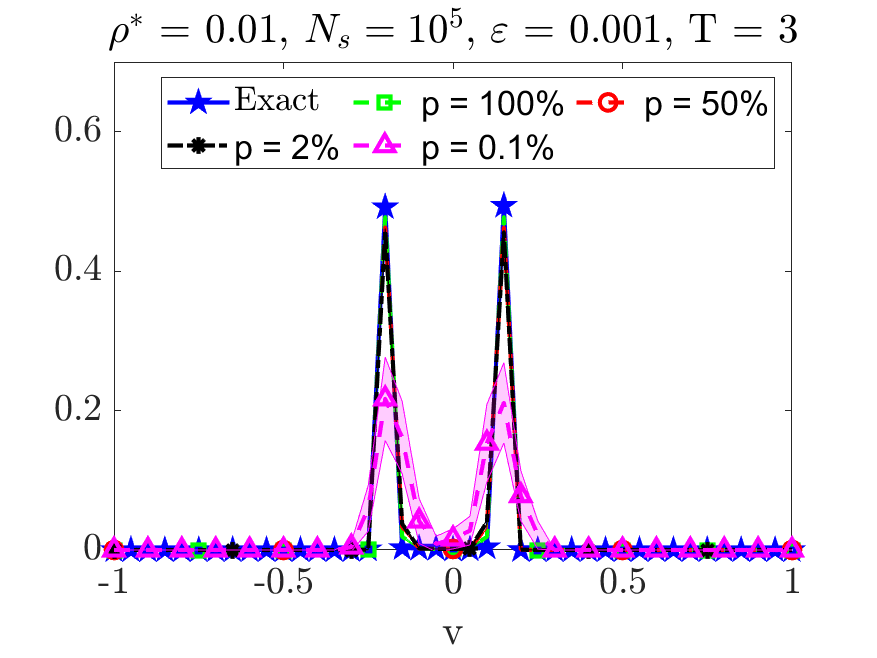}
	\includegraphics[width=0.327\linewidth]{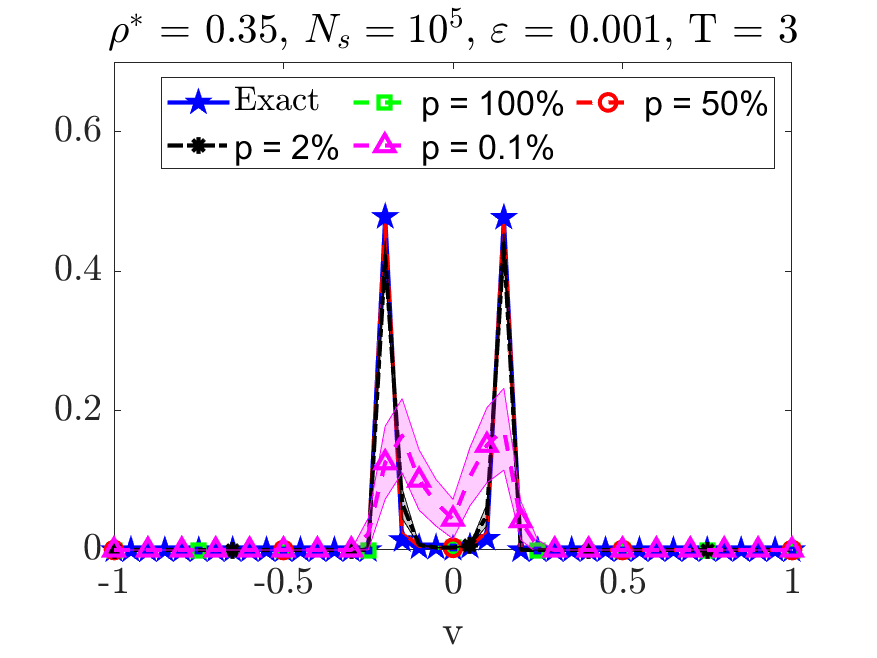}
	\includegraphics[width=0.327\linewidth]{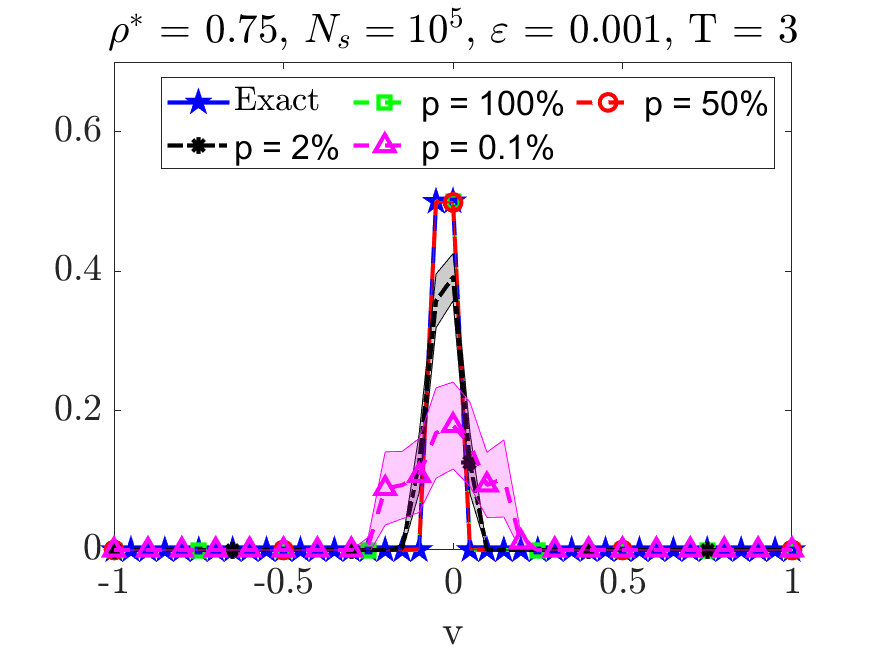}
	\caption{Validation test: comparison between the solution to the kinetic equation in \eqref{eq:boltzmann_ali} computed by means of Asymptotic Nanbu algorithm \ref{alg_binary} and the exact solution in \eqref{eq:exact_sol}. Mean (dashed line) and standard deviation (shaded area) of the velocity distribution computed over $S=100$ simulations  for different values of $p$. From the left to the right $\rho^* = 0.01, 0.35, 0.75$. Markers have been added just to indicate different lines.}
	\label{fig:validation_ord2_ali}
\end{figure}

\begin{figure}[tbhp]
	\centering
	\includegraphics[width=0.327\linewidth]{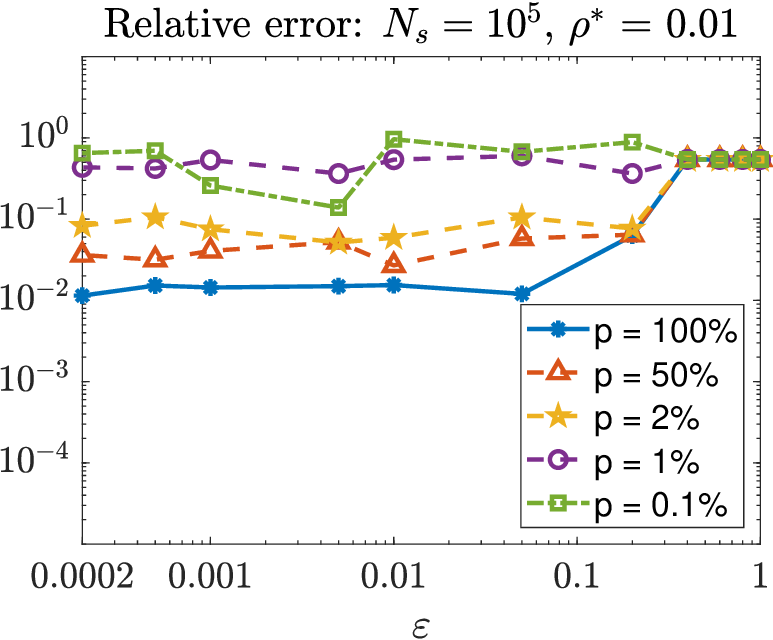}
	\includegraphics[width=0.327\linewidth]{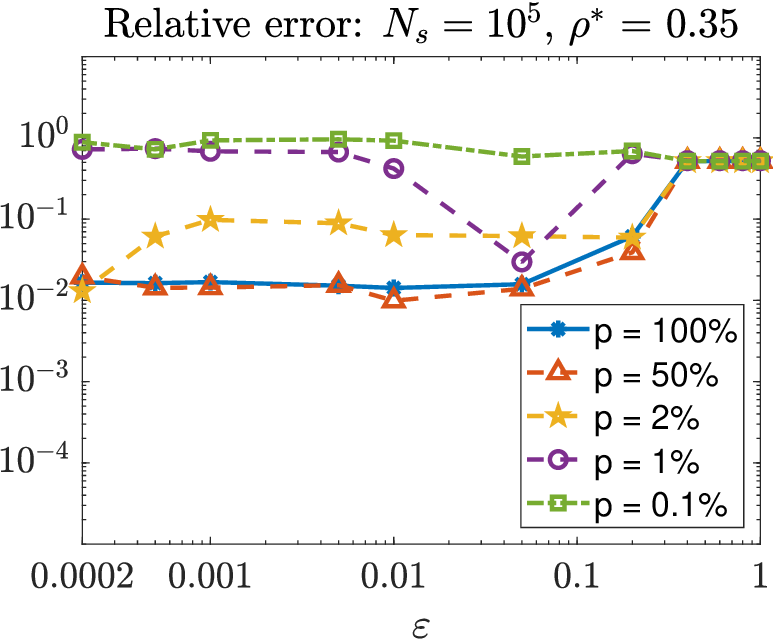}
	\includegraphics[width=0.327\linewidth]{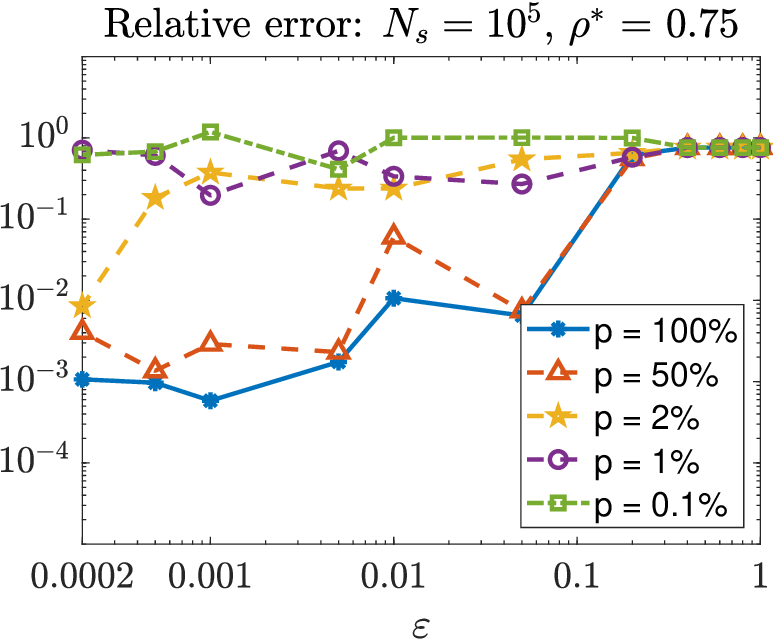}
	\caption{Validation test: $L_2$- norm of the error between the solution to the kinetic equation in \eqref{eq:boltzmann_ali} simulated by means of the Asymptotic Nanbu Algorithm \ref{alg_binary} (one simulation) and the exact solution in \eqref{eq:exact_sol}. From the left to the right $\rho^* =0.01, 0.35, 1$. Markers represent the error for the different values of $\varepsilon$. }
	\label{fig:error_ord2_ali}
\end{figure}

In this section we perform different numerical experiments to test both the accuracy and the efficiency of the Asymptotic Nanbu Algorithm  \ref{alg_binary} with $k$--NN search. 
\paragraph{Accuracy}\label{sec:accuracy} 
Consider a model in which $N_s$ agents with position $x_i$ and velocity $v_i$ interact with their nearest $M$ neighbours without changing their labels and their position. Assume agents are subjected just to alignment forces. Hence, their dynamics at the microscopic level is governed by the following ODE  for $i=1,\ldots,N_s$,
\begin{equation}\label{eq:validation_model_ord1}
	\dot{v}_i =\frac{1}{\rho^*}\frac{1}{N_s} \sum_{j=1}^{N_s}  (v_j-v_i)  \chi_{\mathcal{B}_M(x_i;\xx)}(x_j),
\end{equation}
where $\rho^* = M/N_s$ is the target topological mass.
At the kinetic level, suppose that agents modify their velocity according to binary interactions. Assume that at any time step an agent with position and velocity $(x,v)$ meets another agent with position and velocity $(x_*,v_*)\in \mathit{B}_{r^*_M}(t,x)$ where $r^*_M$ is defined as in \eqref{eq:topological_ball}. Its post-interaction velocity is given by
\begin{equation}\label{eq:binary_1}
	v' = v + \varepsilon (v_*-v),
\end{equation}
where $\varepsilon >0$ is a small parameter. 
Recall that the ball $B_{r^*_M}(t,x)$ by definition contains a certain percentage of mass, that we suppose to be $\rho^*$. If we denote by $f(x,v,t)$ the density of agents at time $t$ with position $x$ and velocity $v$, then the kinetic equation describing its evolution reads 
\begin{equation}\label{eq:boltzmann_ali}
	\partial_t f(x,v,t) = -\nabla_v \cdot  \Big[ f(x,v,t) \int_{\mathit{B}_{r^*(t,x)}\times\mathbb{R}^d} (v_*-v) f(x_*,v_*,t) dx_* dv_* \Big].
\end{equation} 
 The microscopic model in \eqref{eq:validation_model_ord1} can be solved exactly and the evolution of the velocity is given by 
	\begin{equation}\label{eq:exact_sol}
		v_i(t) = v_i(0)e^{-t} + e^{-t} \int_{0}^t \bar{v}_i(t) e^t dt,\quad \text {with } \quad\bar{v}_i(t) = \frac{1}{\rho^*}\frac{1}{N_s} \sum_{j=1}^{N_s}  v_j(t)  \chi_{\mathcal{B}_M(x_i;\xx)}(x_j).
	\end{equation}
	We choose as initial distribution the sum of two 2d-Gaussian in the plane $(x,v)$ one with mean $(-0.33,-0.16)$ and the other with mean $(0.33,0.16)$, and both with standard deviation $(0.12,0.06)$.
	The dynamics at the kinetic level is simulated with Algorithm \ref{alg_binary}, where we compute the velocity change as in equation \eqref{eq:binary_1}.  We suppose $N_s=10^5$, and $\varepsilon = 10^{-5},\ldots,10^0$. We perform the computations assuming the subsample is made with the $p = 100N_c/N_s\%$ of the total mass, for a certain $p$.  In \eqref{fig:validation_ord2_ali_IC} we plot the initial distribution in the $v$-$x$ plane and the marginals in $x$ and $v$. 
 Then, we run $S=100$ simulations and we plot in Figure \ref{fig:validation_ord2_ali} the mean and the standard deviation as a shaded area of the velocity distribution at time $T=3$ for $N_s=10^5$, $\rho^*=0.01, 0.35,0.75$ and $\varepsilon =  10^{-3}$ for different values of $p$.
	In Figure \ref{fig:error_ord2_ali} for different values of $\rho^*$, the $L_2$-norm of the error between the solution to the kinetic equation in \eqref{eq:boltzmann_ali}, simulated by means of the Asymptotic Nanbu  Algorithm \ref{alg_binary} (one simulation) for different values of $p$, and the exact solution in \eqref{eq:exact_sol}. Note that we observe a saturation effect for $\varepsilon \approx 10^{-2}$.  

\begin{figure}[tbhp]
	\centering
	\includegraphics[width=0.327\linewidth]{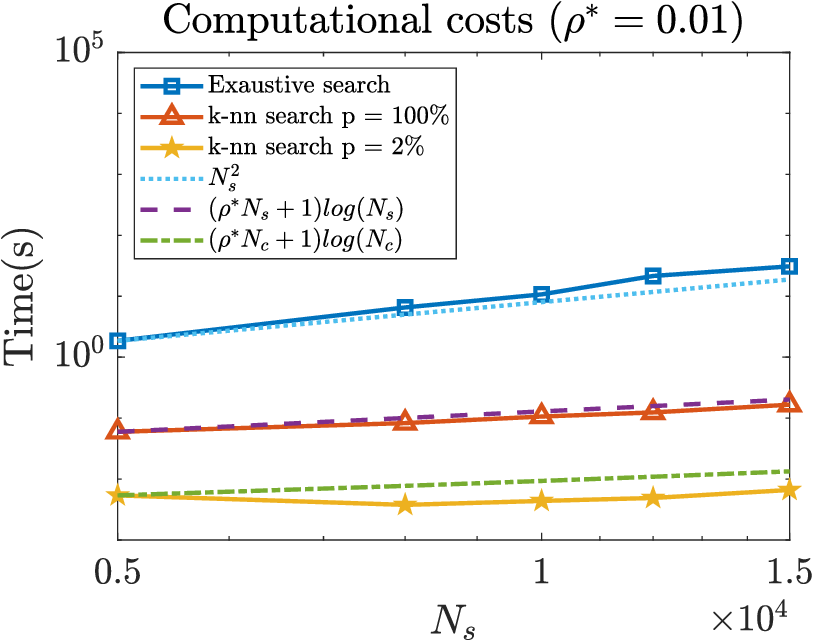}
	\includegraphics[width=0.327\linewidth]{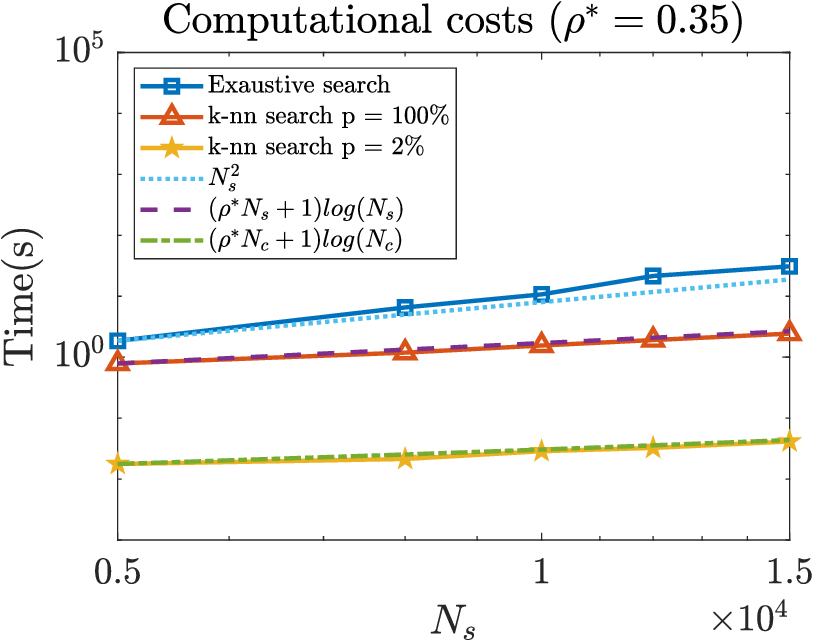}
	\includegraphics[width=0.327\linewidth]{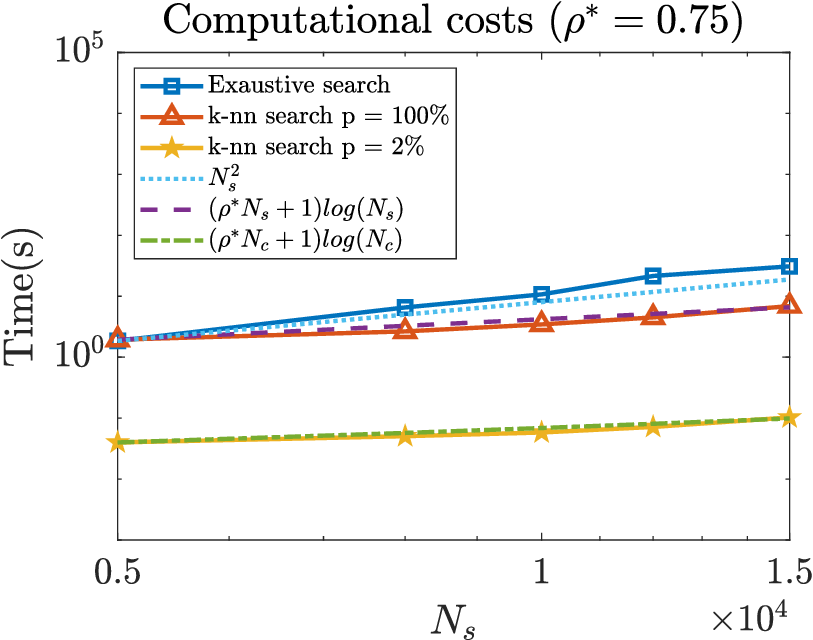}
	\caption{Comparison between the computational costs of the exhaustive search and the $k$--NN search for different values of $N_s$ and as the subsample percentage size varies from $p=100\%$ to $p=2\%$. From the left to the right $\rho^* = 0.01,0.35,0.75$.  Markers represent the computational costs relative to the different values of $N_s$.  }
	\label{fig:cost}
\end{figure}

\begin{figure}[tbhp]
	\centering
	\includegraphics[width=0.327\linewidth]{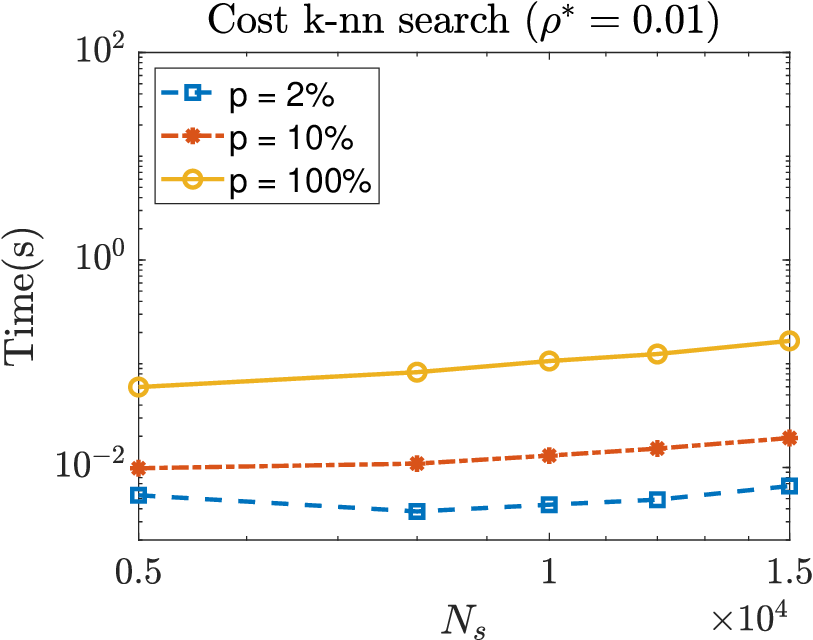}
	\includegraphics[width=0.327\linewidth]{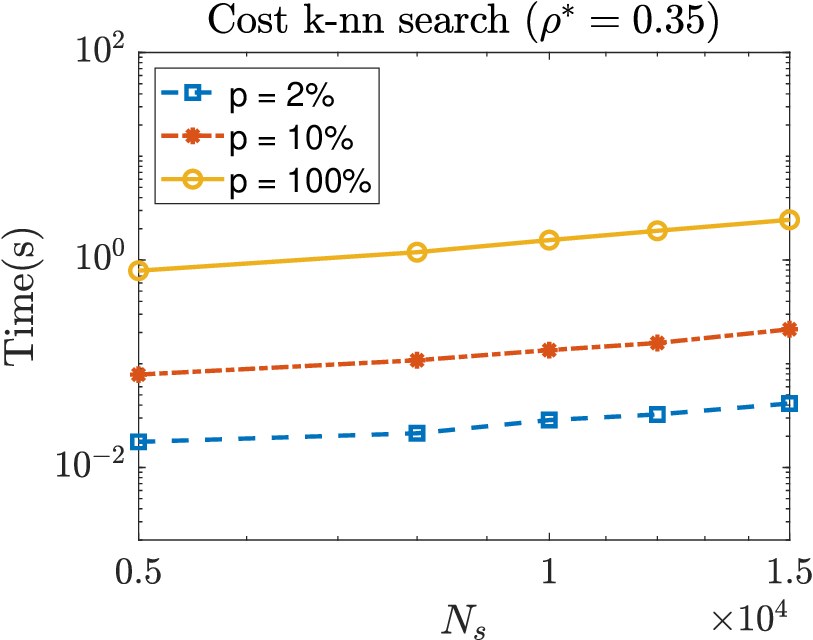}
	\includegraphics[width=0.327\linewidth]{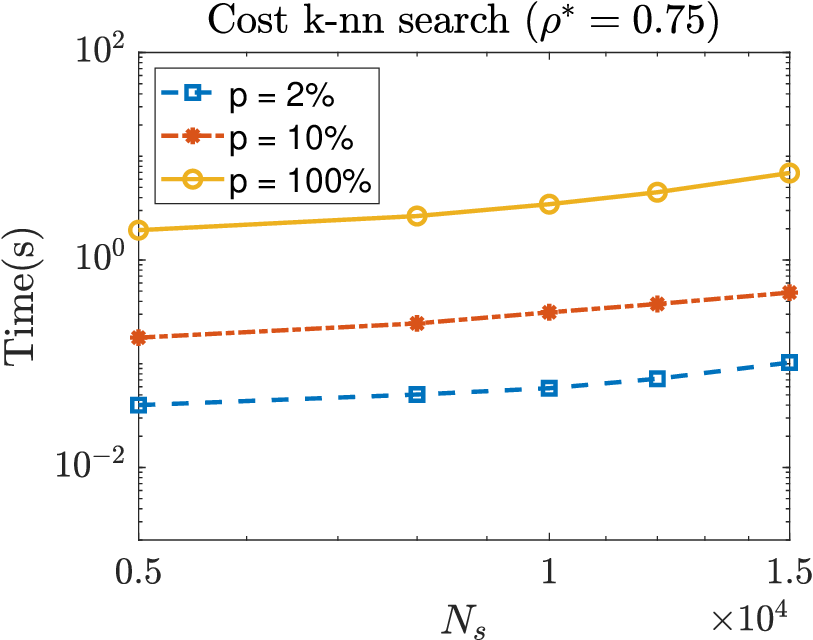}
	\caption{Comparison between the computational costs of the $k$--NN search for different values of $N_s$ and of the percentage $p$ of the subsample size. From the left to the right $\rho^* = 0.01,0.35,0.75$.  Markers represent the computational costs relative to the different values of $N_s$.  }
	\label{fig:cost_1}
\end{figure}

%%%%%%%%%%%%%%%%%%%%%%%%%%%%%%%%%%%%%%%%%%%%%%%
%%%%%%%%%%%%%%%%%%%%%%%%%%%%%%%%%%%%%%%%%%%%%%%%
\paragraph{Computational costs}\label{sec:costs} We now compare the computational costs of the exhaustive search and the $k$--NN search. 
The computational cost of an exhaustive search is $\mathcal{O}(d N_s^2)$, where $d$ is the space dimension and $N_s$ the number of particles. Indeed, first one needs to compute the distances between each point and all the others, with a cost of $\mathcal{O}(d N_s^2)$, and then to sort them, with a cost of $\mathcal{O}(N_s^2  log(N_s))$. The cost of a $k$--NN search is logarithmic in time. First one needs to organize agents optimally with a $k$-$d$ tree. The cost of this operation is proportional to $N_s log(N_s)$. Then the idea is to perform a search over the tree to select which are the nearest agents. It can be shown (see \cite{friedman1977algorithm}) that the $k$--NN search algorithm examines the nodes in optimal order, that is in order of increasing dissimilarities, and that the number of nodes that should be examined is proportional to $((\rho^*N_s)^{1/d}+1)^d$. Hence, the total cost to construct a $k$-$d$ tree and to perform the search over it is
\begin{equation}\label{eq:cost_knn}
	\mathcal{O}(max(((\rho^*N_s)^{1/d}+1)^d log(N_s), N_s log(N_s)).
\end{equation}
If instead we consider a subsample of size $N_c<<N_s$, the total cost is further reduced to 
	\begin{equation}\label{eq:cost_knn_p}
		\mathcal{O}(max(((\rho^*N_c)^{1/d}+1)^d log(N_c), N_c log(N_c)).
	\end{equation}
	In Figure \ref{fig:cost} we see the comparison between the computational cost to perform one exhaustive and one $k$--NN search as $N_s$ varies for different values of $\rho^*$. In particular, we compare the classical $k$--NN search, performed on the whole sample $(p = 100\%)$, and the one performed on a subsample of size $N_c$, with $N_c = 10^{-2} N_s p$, being $p$ the subsample percentage.   We set $N_s = 5\times 10^{3},\ldots, 1.5\times 10^4$ and $p=2\%$.  The $k$--NN computational cost increases as $\rho^*$ increases, and it is proportional to \eqref{eq:cost_knn}-\eqref{eq:cost_knn_p}, decreasing as the subsample percentage size $p$ decreases. More in details, in Figure \ref{fig:cost_1}  we see a comparison between the computational costs of a $k$--NN search as $N_s$ varies for different subsample percentage sizes $p$. Again, we see that the computational cost decreases proportionally to the subsample size, as expected. 
%%%%%%%%%%%%%%%%%%%%%%%%%%%%%%%%%%%%%%%%%%%%%%%%%%%%%%%%%%
%%%%%%%%%%%%%%%%%%%%%%%%%%%%%%%%%%%%%%%%%%%%%%%%%%%%%%%%%%%%%
%%%%%%%%%%%%%%%%%%%%%%%%%%%%%%%%%%%%%%%%%%%%%%%%%%%%%%%%%%%%%%%
\section{Numerical experiments} \label{sec:2D3Dexperiments} 
We present different numerical experiments simulating the two and three dimensional dynamics both at the microscopic and mesoscopic levels.  The dynamics at microscopic level is discretized by a forward Euler scheme with a time step $\Delta t = 0.01$, whereas the evolution of the
kinetic dynamics is approximated by the Asymptotic Nanbu Algorithm \ref{alg_binary} with $\varepsilon = 0.01$. The time evolution of the labels is computed with Algorithm \ref{alg_lambda} at both the microscopic and mesoscopic levels. In the microscopic case we set $N=400$. In the mesoscopic case we choose a sample
of $O(N_s)$ particles, with $N_s =5\times 10^5$, and a subsample of $O(N_c)$ particles, with $N_c=10^4$ that corresponds to a percentage $p=2\%$ of the total mass, for the approximation of the density. We assume the topological target mass to be $\rho^* = 0.01$. Table \ref{tab:all_parameters} reports the parameters of
the model that remain unchanged in the various scenarios.
\begin{table}[tbhp]
	\footnotesize
	\caption{Model parameters for the different scenarios.}\label{tab:all_parameters}
	\begin{center}
		\begin{tabular}{|c|c|c|c|c|c|c|c|c|c|}\hline
			&$C_{rep}$ & $C_{ali}$ & $C_{att}$ & $C_{v}$ &s&$\overline r$ & $\underline{r}$&$\epsilon$ \\
			\hline
			2D model & 100 & 12 &  0.7 &5&10&200&1&200\\
			3D model & 100 & 12 &  0.7 &5&10&350&20&150\\
			\hline
		\end{tabular}			
	\end{center}
\end{table}
The other parameters will be specified later.
%%%%%%%%%%%%%%%%%%%%%%%%%%%%%%%%%%%%%%%%%%%%%%%%%%
%%%%%%%%%%%%%%%%%%%%%%%%%%%%%%%%%%%%%%%%%%%%%%%%%%%
%%%%%%%%%%%%%%%%%%%%%%%%%%%%%%%%%%%%%%%%%%%%%%%%%%%
%%%%%%%%%%%%%%%%%%%%%%%%%%%%%%%%%%%%%%%%%%%%%
%%%%%%%%%%%%%%%%%%%%%%%%%%%%%%%%%%%%%%%%%%%%

\subsection{Numerical test in two spatial dimensions}\label{sec:test2D} 
We consider the swarming dynamics  evolving on the spatial space $(x,y)\in \R^2$ and velocity space $(v_x,v_y)\in \R^2$.
\subsubsection{Test 2D with no food sources}\label{sec:2Dtest_nofood} 
\begin{figure}[tbhp]
	\centering
	\includegraphics[width=0.48\linewidth]{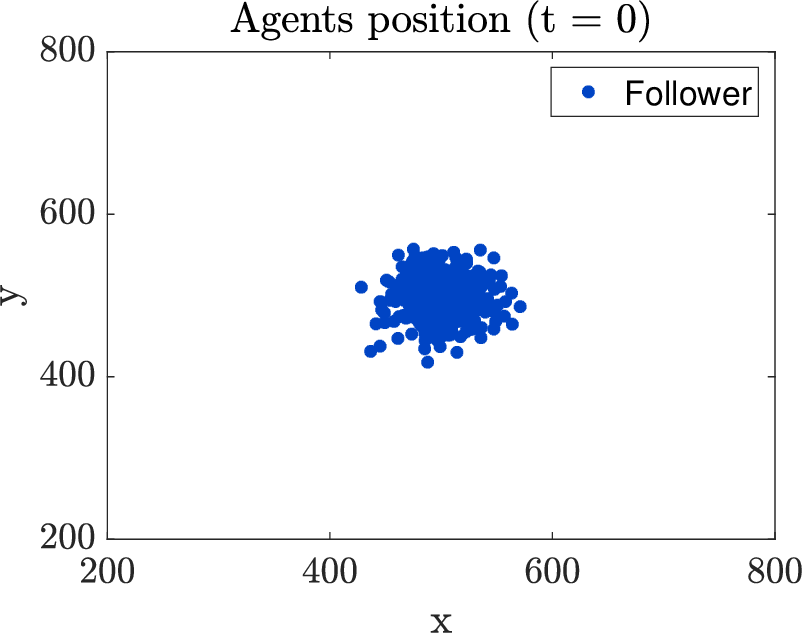}
	\includegraphics[width=0.505\linewidth]{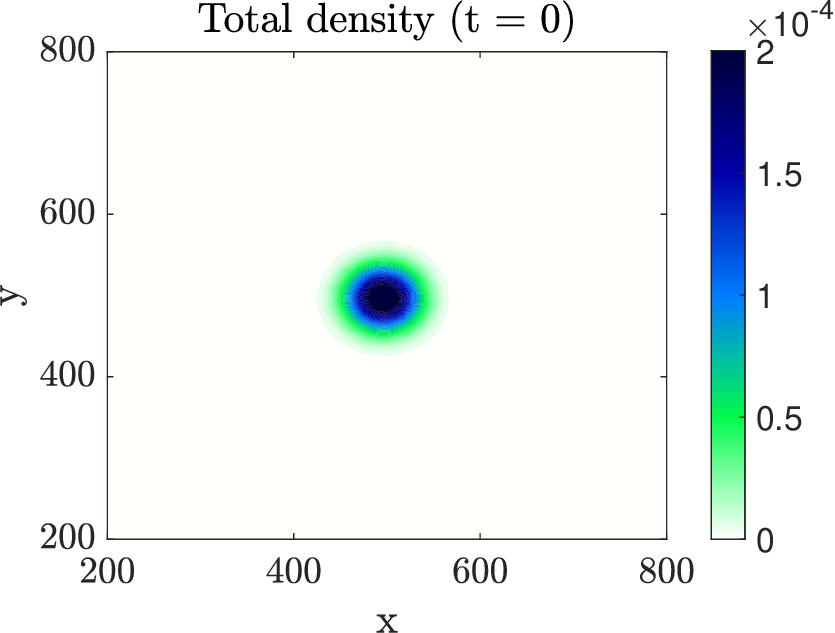}
	\caption{Initial configuration in the spatial 2D case with no food sources.}
	\label{fig:2D_initial_configuration_nofood}
\end{figure}

Suppose the model includes no food sources, i.e. $C_{src} = 0$, and no attraction to the centre of mass, i.e. $C_{ctr} = 0$. We simulate the dynamics up to time $T = 500$, and we report in Figure \ref{fig:2D_initial_configuration_nofood} the initial configuration for both the microscopic and mesoscopic dynamics.

\begin{figure}[tbhp]
	\centering
	\includegraphics[width=0.327\linewidth]{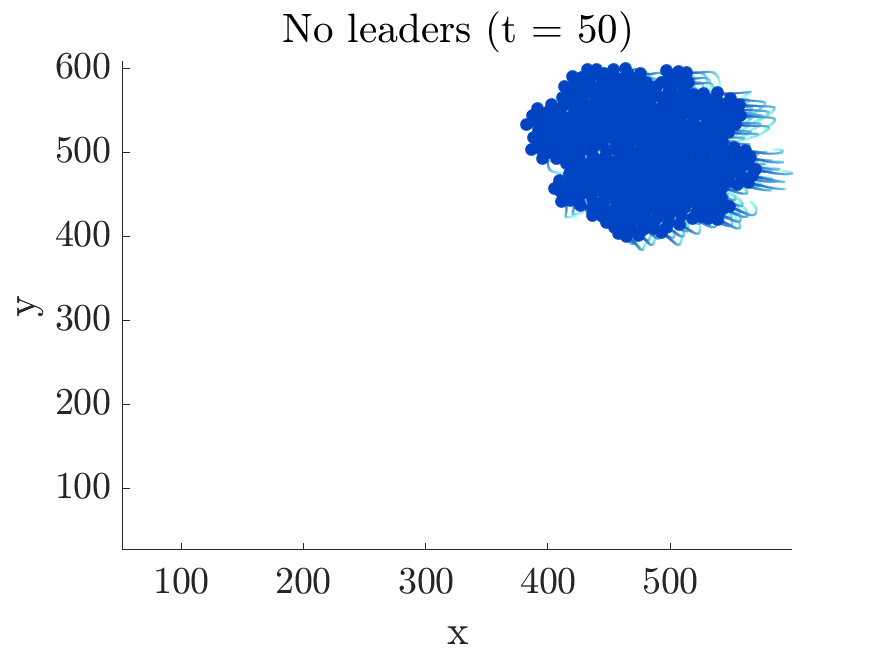}
	\includegraphics[width=0.327\linewidth]{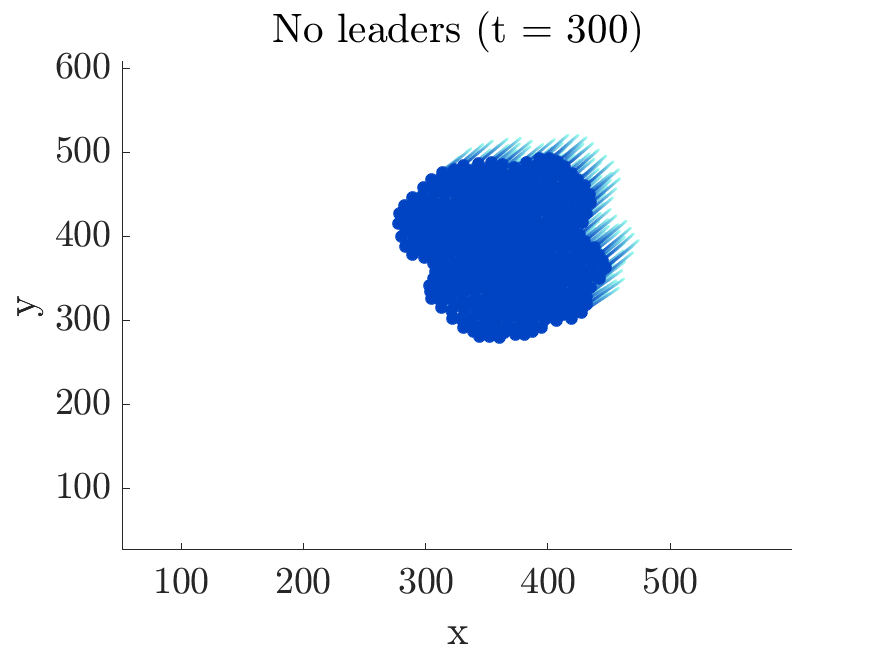}
	\includegraphics[width=0.327\linewidth]{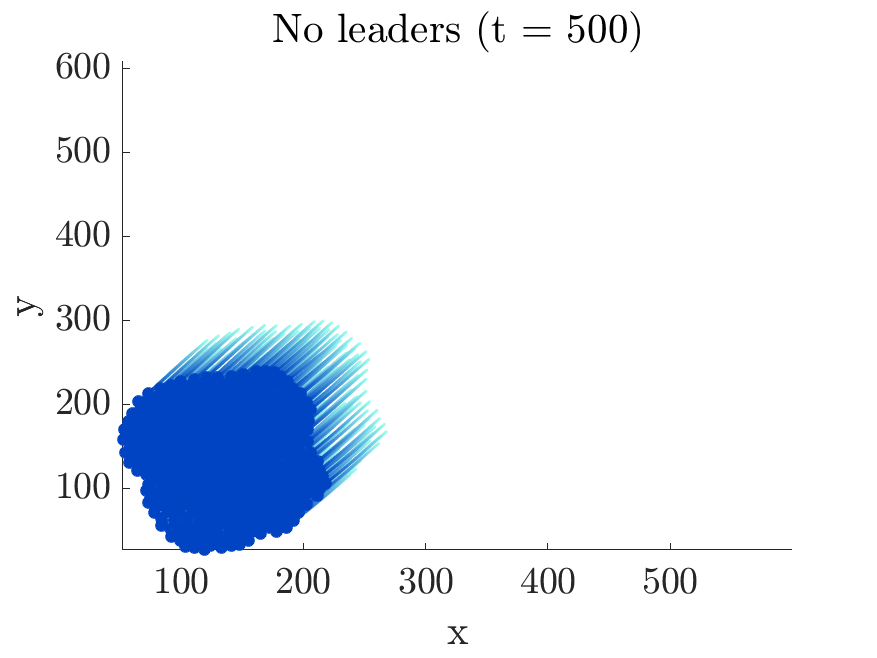}\\
	\includegraphics[width=0.327\linewidth]{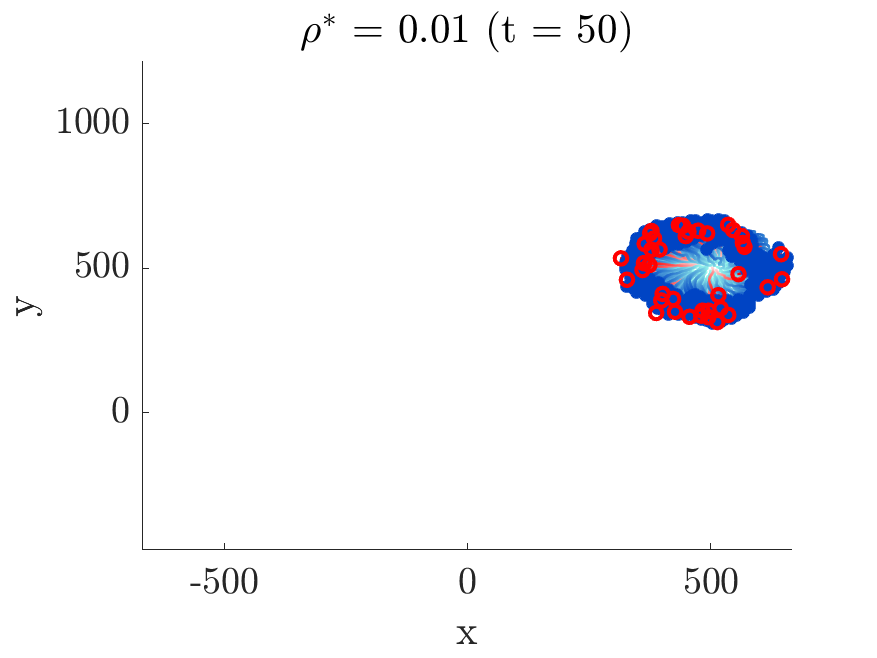}
	\includegraphics[width=0.327\linewidth]{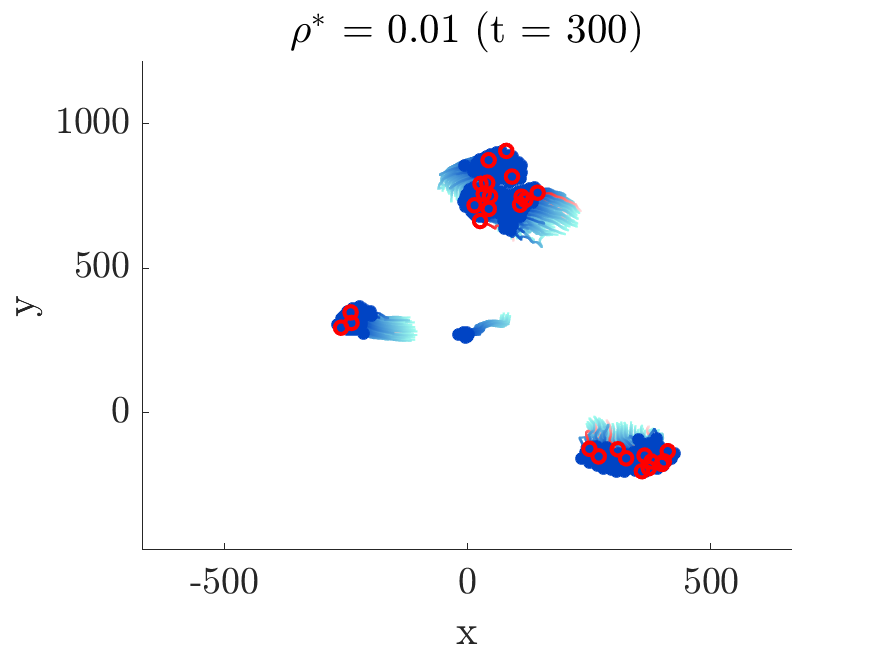}
	\includegraphics[width=0.327\linewidth]{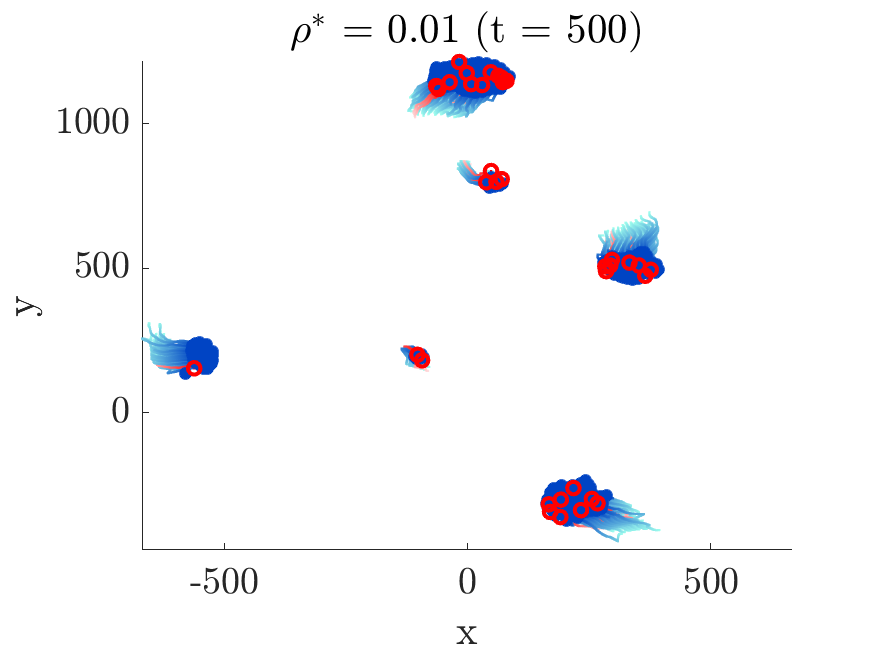}
	\caption{Three snapshots of the 2D dynamics at the microscopic level described in \eqref{eq:dynamics} with $\lambda$ evolving with rates \eqref{eq:rates_test_0} taken at time $t=50$, $t=300$ and $t=500$ and assuming $C_{ctr} = C_{src}=0$. In the first row, the dynamics without leaders. In the second row, the dynamics with leaders assuming $\rho^*=0.01$. We represent in blue the agents in the followers status and in red the ones in the leaders status.
	}
	\label{fig:micro2D_nofood_dynamics}
\end{figure}

At time $t=0$, agents are normally distributed with mean $\mu = 500$ and variance $\sigma^2 = 25^2$ and are in the followers status. Labels change according to the transition rates defined in \eqref{eq:rates_test_0} with $q_{FL}=2\times 10^{-4}$ and $q_{LF}= 4\times 10^{-3}$.
\paragraph{Microscopic case}
In Figure \ref{fig:micro2D_nofood_dynamics} we report three snapshots of the dynamics at time $t=50$, $t=300$ and $t=500$, for the dynamics without leaders' emergence (top row) and with leaders' emergence (bottom row). We observe that, without leaders, agents align and form a compact swarm, whereas  with leaders' emergence we observe the formation of different groups. The splitting is not symmetric since leaders' emergence occurs randomly and this is reflected in the cluster formation.

\begin{figure}[tbhp]
	\centering
	\includegraphics[width=0.327\linewidth]{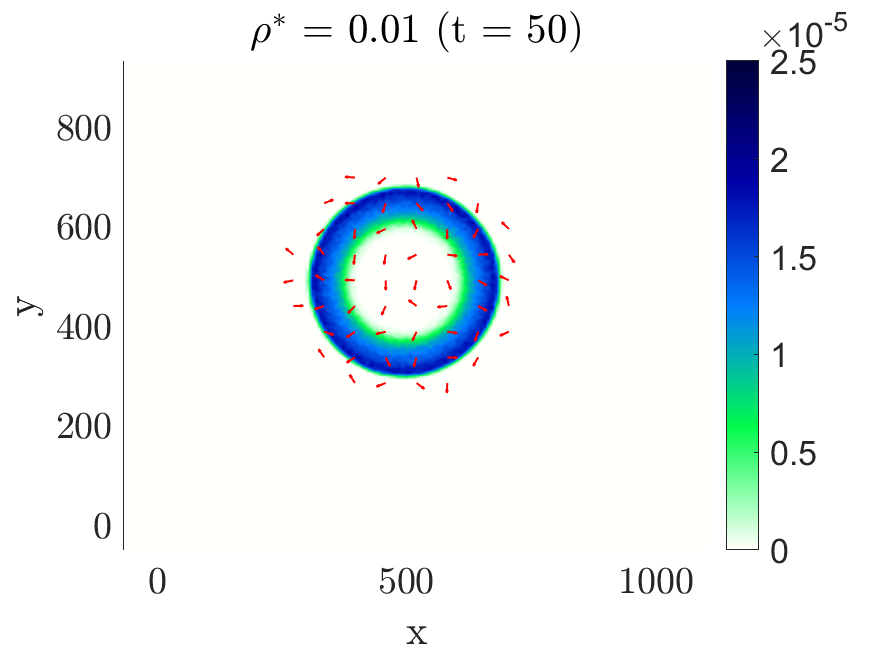}
	\includegraphics[width=0.327\linewidth]{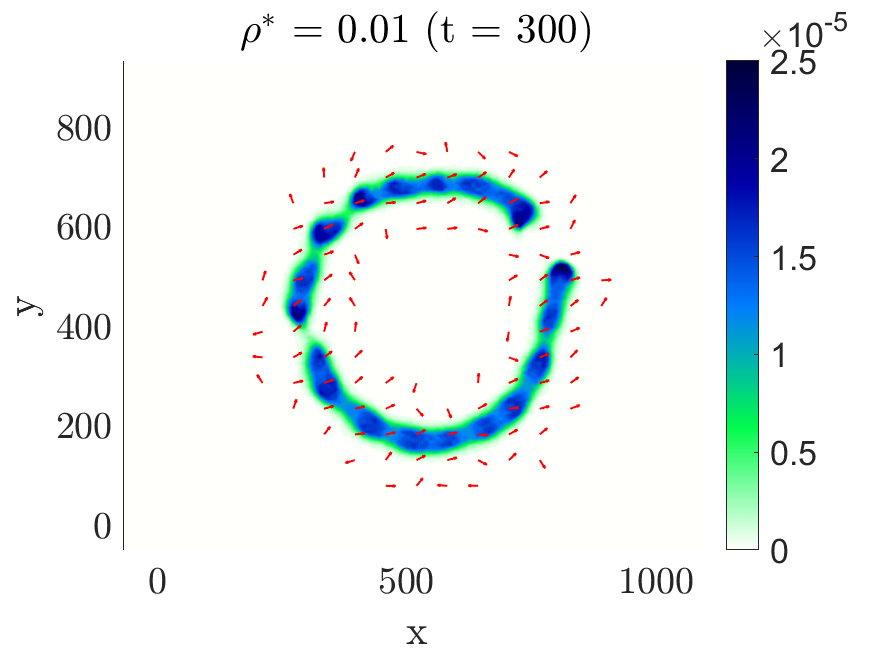}
	\includegraphics[width=0.327\linewidth]{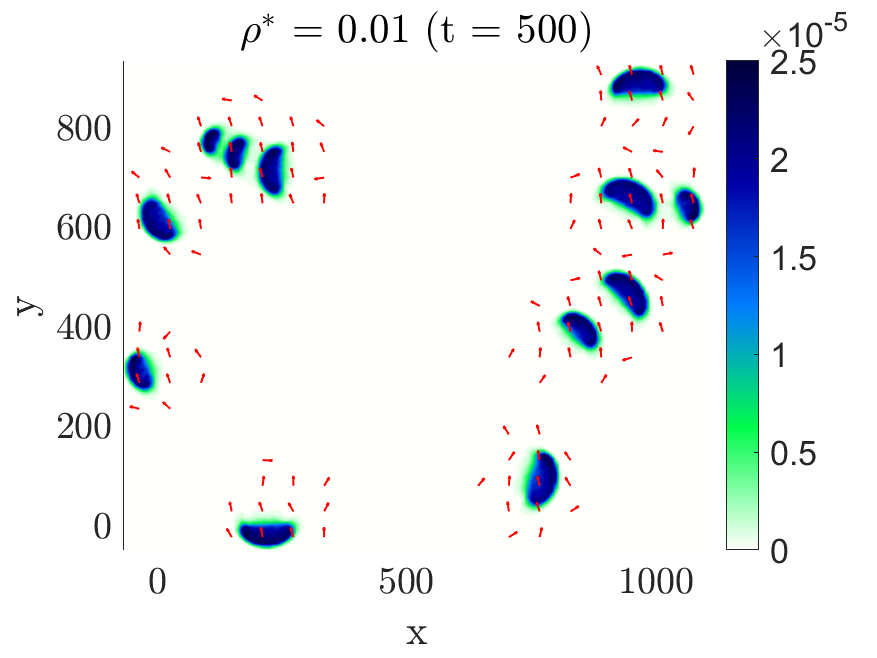}\\
	\includegraphics[width=0.327\linewidth]{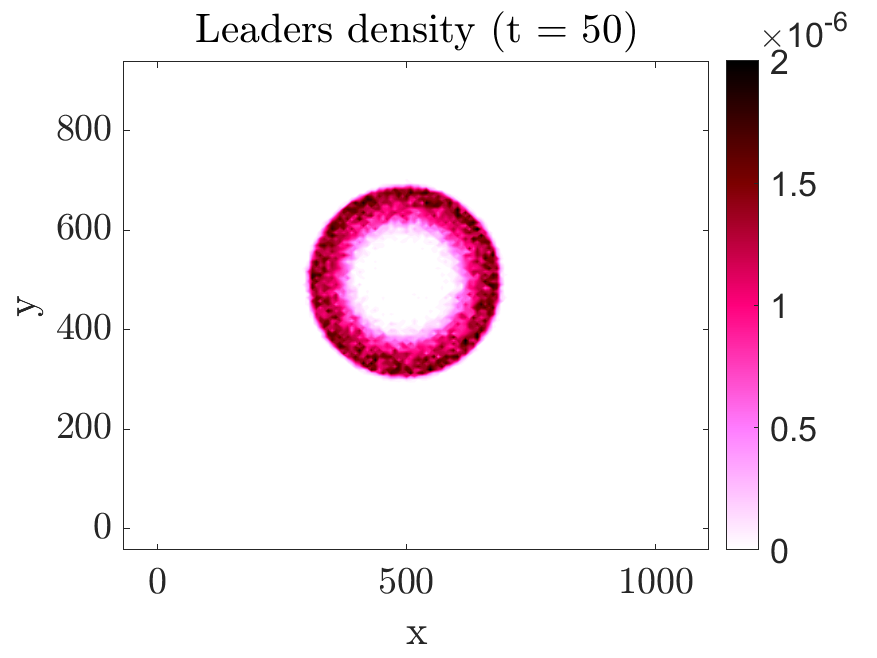}
	\includegraphics[width=0.327\linewidth]{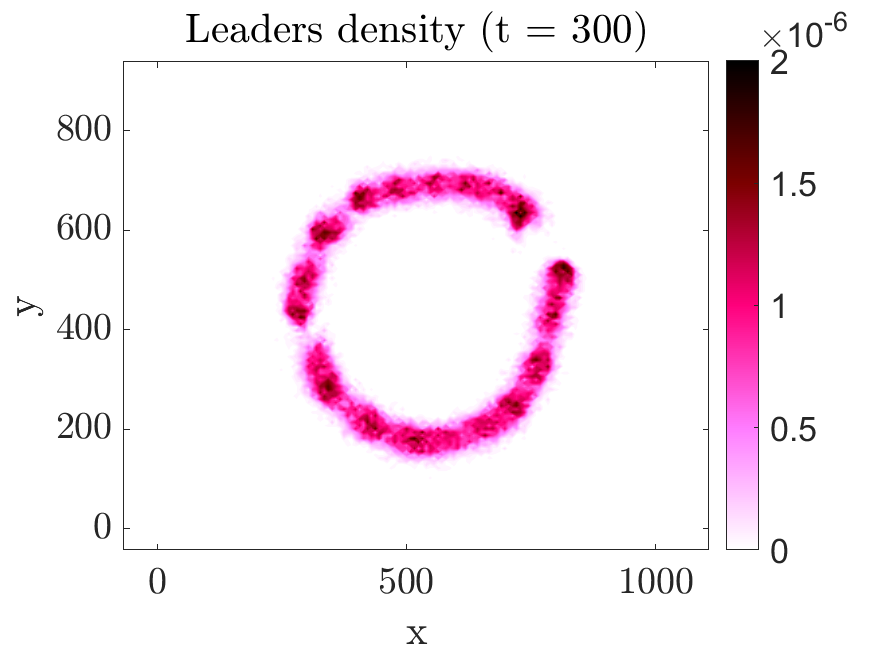}
	\includegraphics[width=0.327\linewidth]{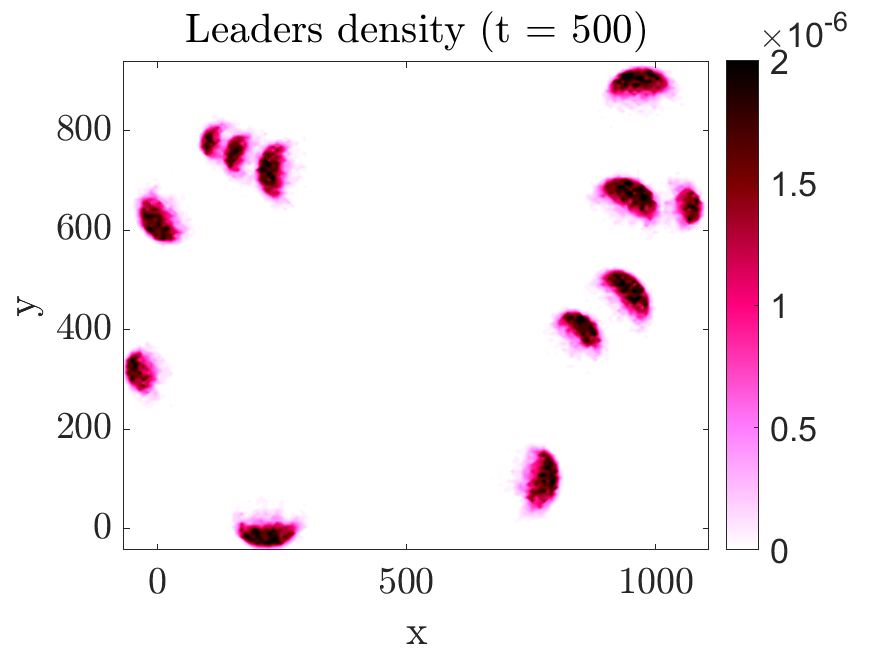}\\
	\caption{Three snapshots of the 2D dynamics at the mesoscopic level described in \ref{eq:boltz_strong} and simulated by means of the Asymptotic Nanbu  Algorithm \ref{alg_binary} with $\lambda$ evolving with rates \eqref{eq:rates_test_0} as in  Algorithm \ref{alg_lambda}, taken at time $t=50$, $t=300$ and $t=500$ and assuming $C_{ctr} = C_{src}=0$. In the first row, the dynamics of the total mass and in red the velocity vector field of the leaders. In the second row, the leaders' dynamics. 
	}
	\label{fig:meso2D_nofood_dynamics}
\end{figure}

\paragraph{Mesoscopic case}
In Figure \ref{fig:meso2D_nofood_dynamics} we report three snapshots of the dynamics at time $t=50$, $t=300$, $t=500$. In the first row, the time evolution of the total density and in red the velocity vector field of the leaders. In the second row, the evolution of the leaders' density. The behavior is similar to the one of the microscopic case, where we observe the formation of various clusters, and the emergence of leaders uniformly over the swarm density.

\begin{figure}[tbhp]
	\centering
	\includegraphics[width=0.48\linewidth]{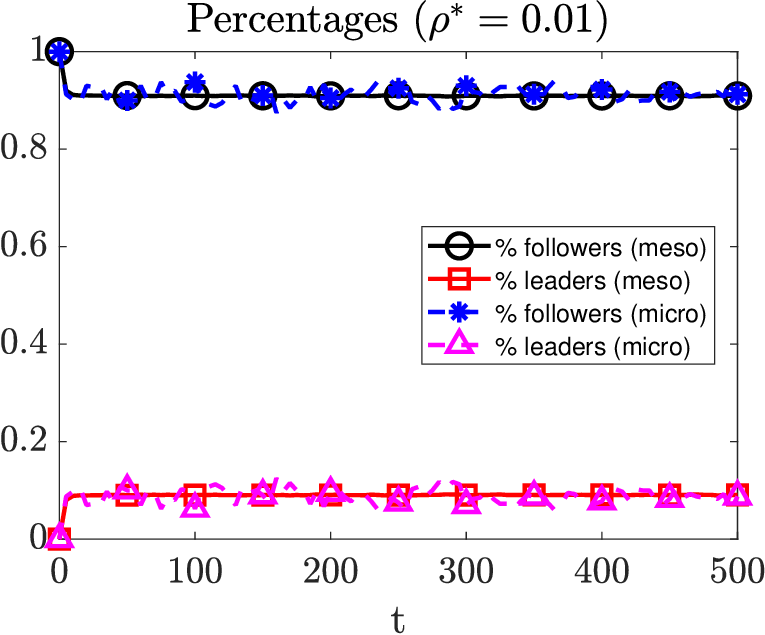}
	\caption{Agents percentages. Markers have been added just to distinguish the different lines. 
	}
	\label{fig:meso2D_nofood_percentages}
\end{figure}

In Figure \ref{fig:meso2D_nofood_percentages} the agents percentages for the dynamics with leaders. 
The videos of the simulations of this subsection are available at \href{https://drive.google.com/drive/folders/1VsO4ffzQvbMb5mvG3pHLoEpA-QyKlkbL?usp=share_link}{[VIDEO]}.
%%%%%%%%%%%%%%%%%%%%%%%%%%%%%%%%%%%%%%%%%%%%%%%%%%%%%%%%%%%%%%%
%%%%%%%%%%%%%%%%%%%%%%%%%%%%%%%%%%%%%%%%%%%%%%%%%%%%%%%%%%%%%%%

\subsubsection{Test 2D: two food sources}\label{sec:2Dtwofood} 
\begin{figure}[tbhp]
	\centering
	\includegraphics[width=0.327\linewidth]{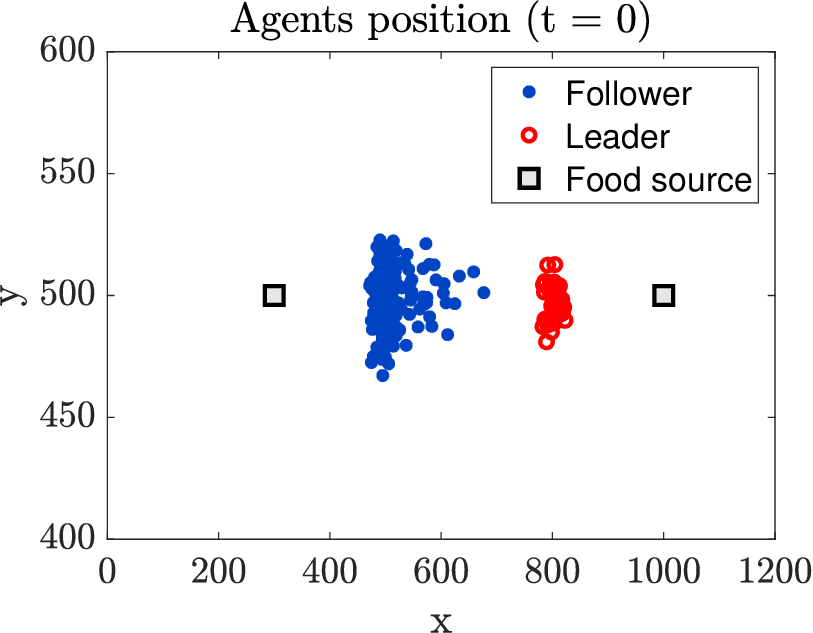}
	\includegraphics[width=0.327\linewidth]{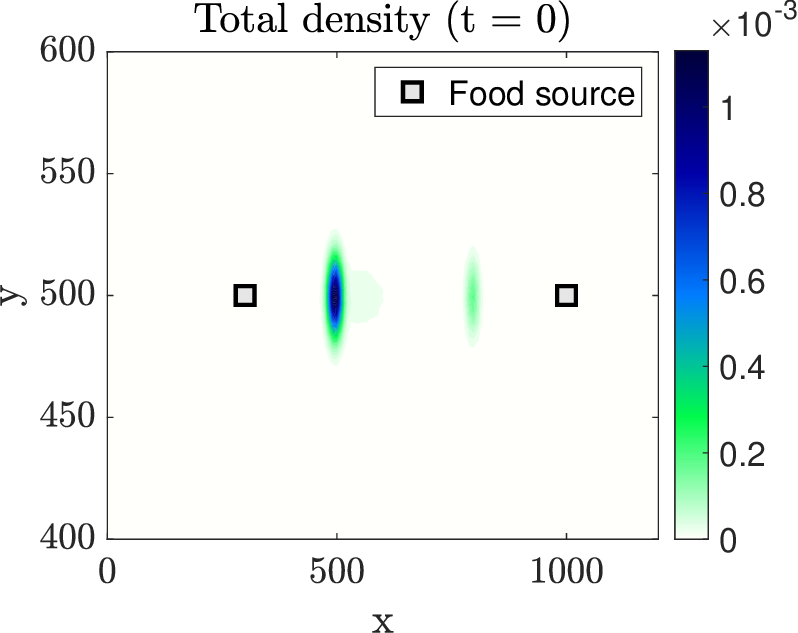}
	\includegraphics[width=0.327\linewidth]{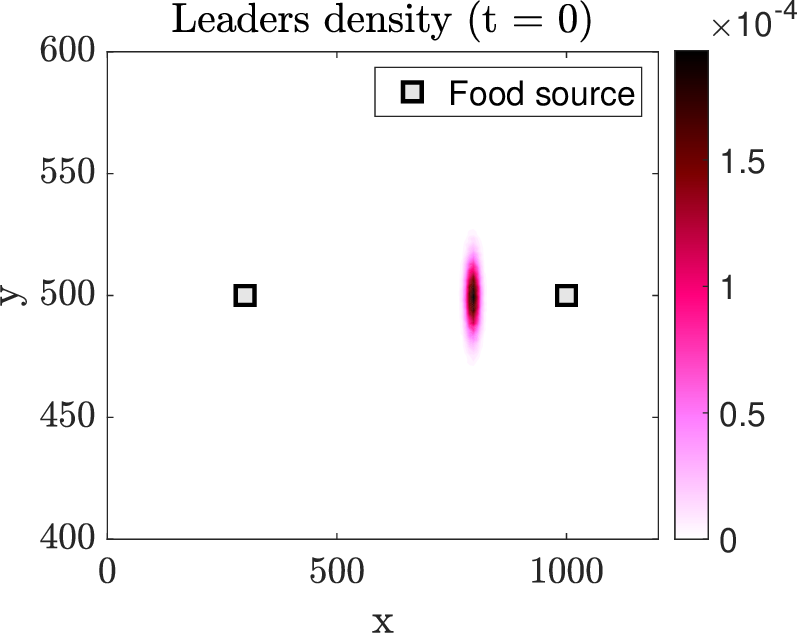}
	\caption{Initial configuration in the 2D case with two food sources.}
	\label{fig:meso2D_initial_configuration}
\end{figure}

Assume the model includes two food sources located in $x_1^{src} = (300,500)$ and $x_2^{src} = (1000,500)$. Assume $C_{ctr} = C_{src}=0.75$. Run the dynamics until time $T=200$.  In Figure \ref{fig:meso2D_initial_configuration} the initial configuration for the microscopic and mesoscopic case. 
We assume that initially the $87.5\%$ of agents is in the followers status. Among them the $75\%$ is normally distributed with mean $\mu = 550$ and variance $\sigma^2 = 10^2$ while the $12.5\%$ is  normally distributed with mean $\mu = 650$ and variance $\sigma^2 = 50^2$. The remaining $12.5\%$ is in the leaders status and it is normally distributed with mean $\mu = 800$ and variance $\sigma^2 = 10^2$. New leaders emerge with higher probability where the followers concentration is higher. Leaders return in the follower status with higher probability if the followers concentration around their position is lower.   Hence we consider density dependent transition rates defined in equation \eqref{eq:rates_test} with $q_L = 4\times 10^{-3}$ and $q_F = 3\times 10^{-3}$ and $\delta = 10^3$.

\begin{figure}[tbhp]
	\centering
	\includegraphics[width=0.327\linewidth]{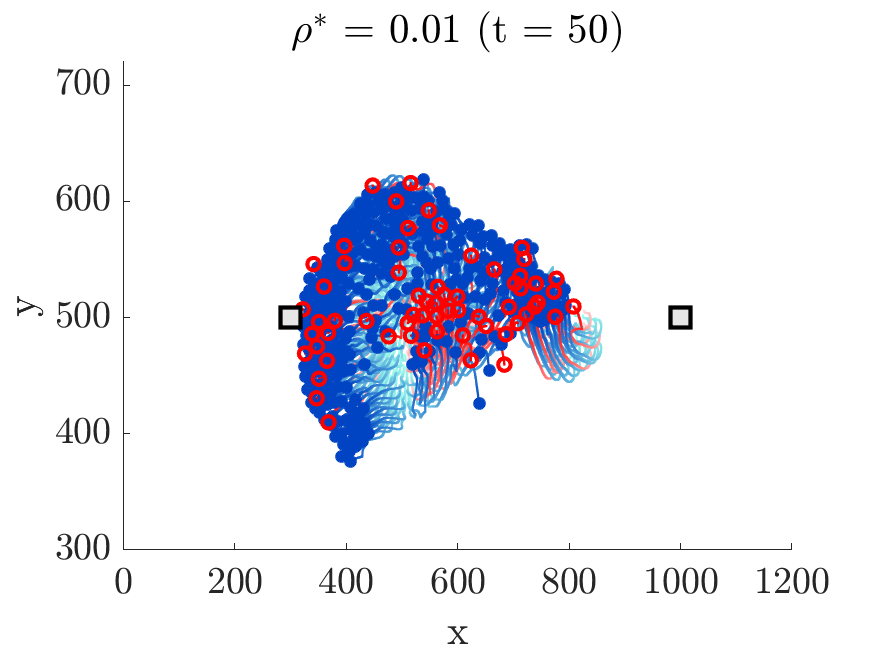}	\includegraphics[width=0.327\linewidth]{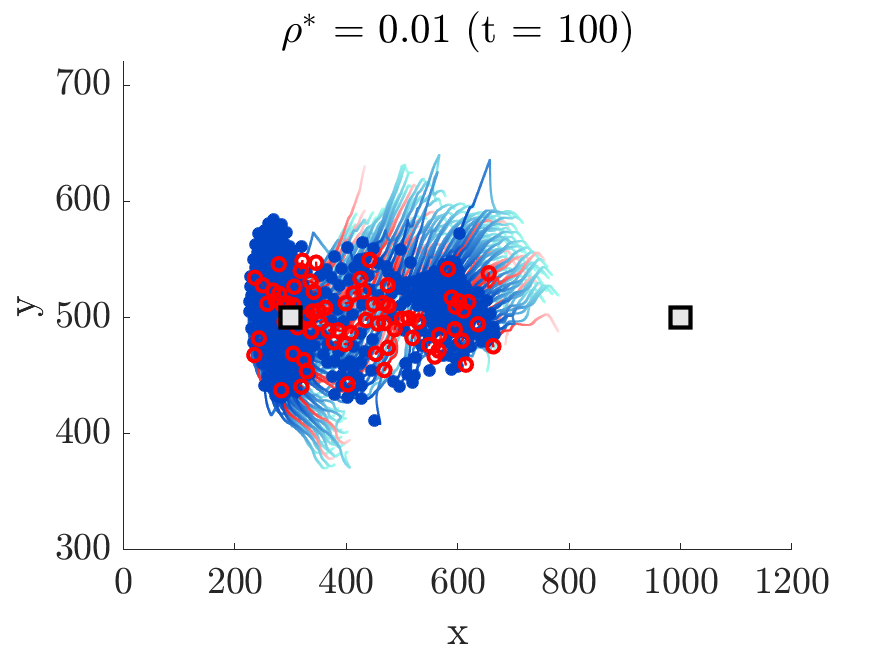}
	\includegraphics[width=0.327\linewidth]{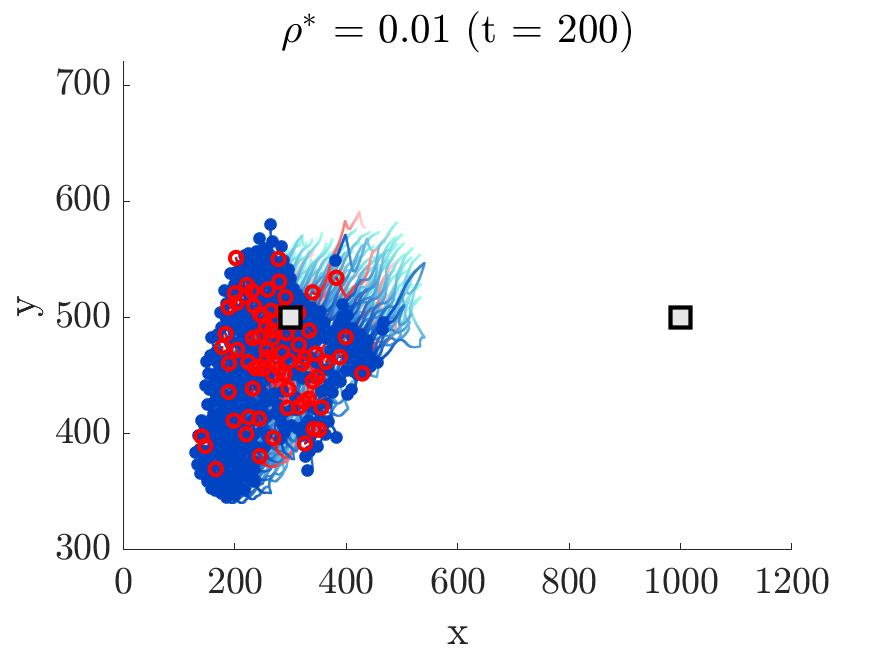}
	\caption{Three snapshots of the 2D dynamics at the microscopic level described in \eqref{eq:dynamics} with $\lambda$ evolving with rates \eqref{eq:rates_test} taken at time $t=50$, $t=100$ and $t=200$ and assuming $C_{ctr} = C_{src}=0.75$. We represent in blue the agents in the followers status and in red the ones in the leaders status.
	}
	\label{fig:micro2D_dynamics}
\end{figure}

\begin{figure}[tbhp]
	\centering
	\includegraphics[width=0.33\linewidth]{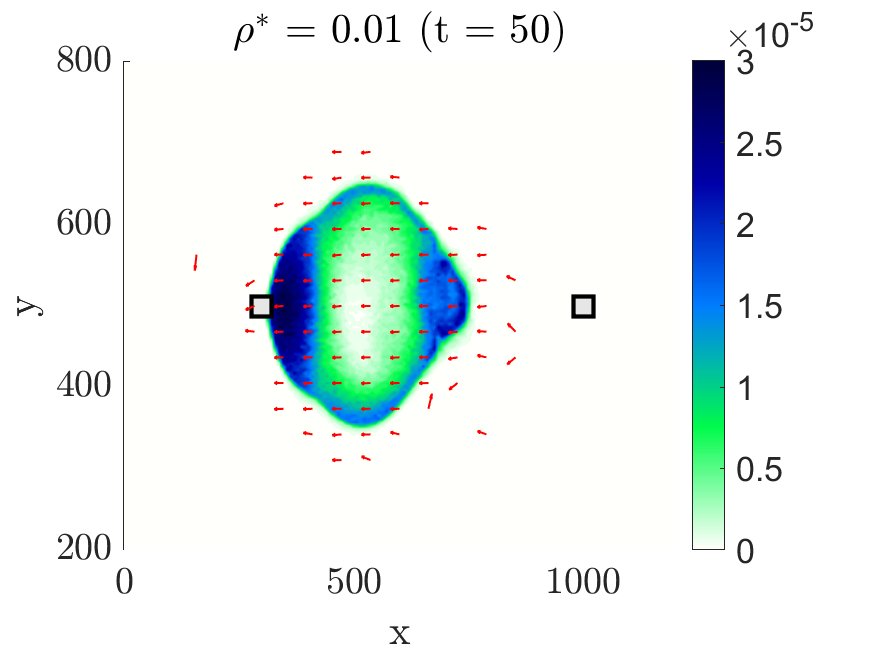}
	\includegraphics[width=0.327\linewidth]{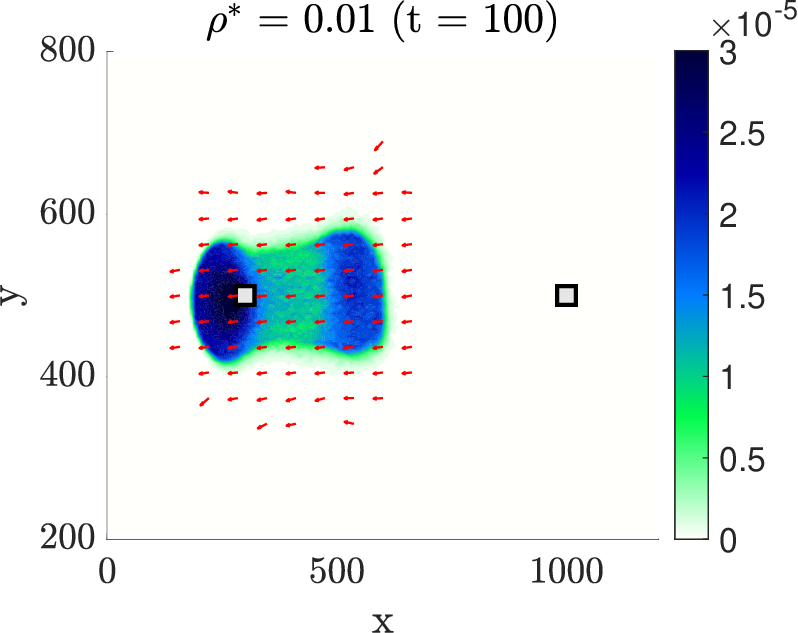}
	\includegraphics[width=0.327\linewidth]{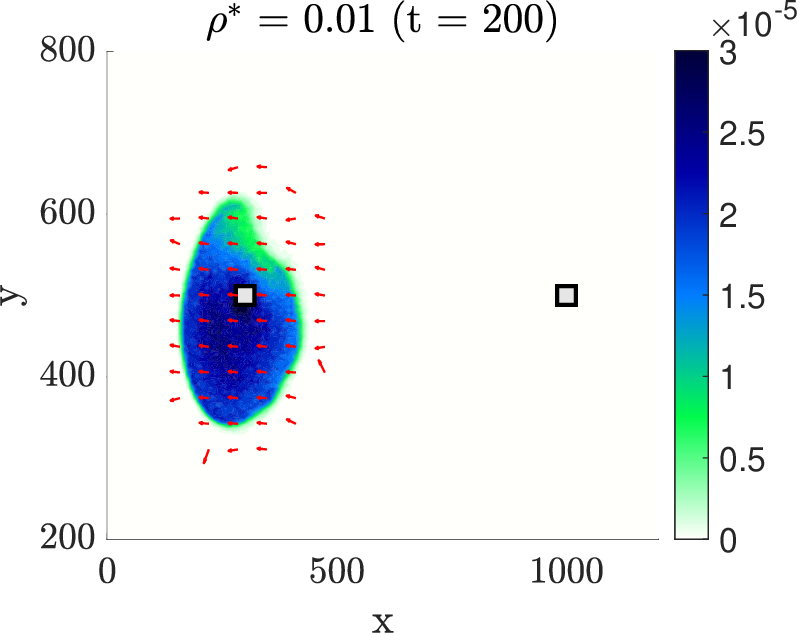}\\
	\includegraphics[width=0.327\linewidth]{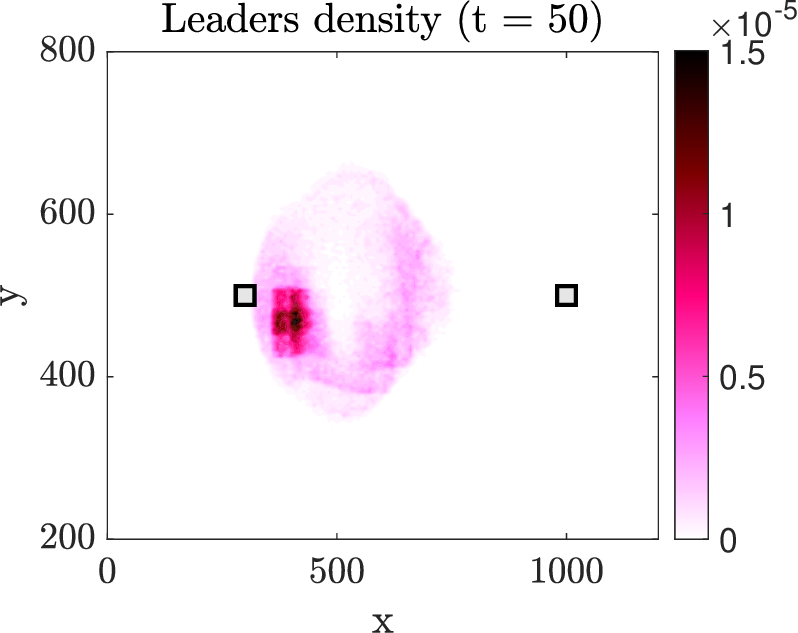}	\includegraphics[width=0.327\linewidth]{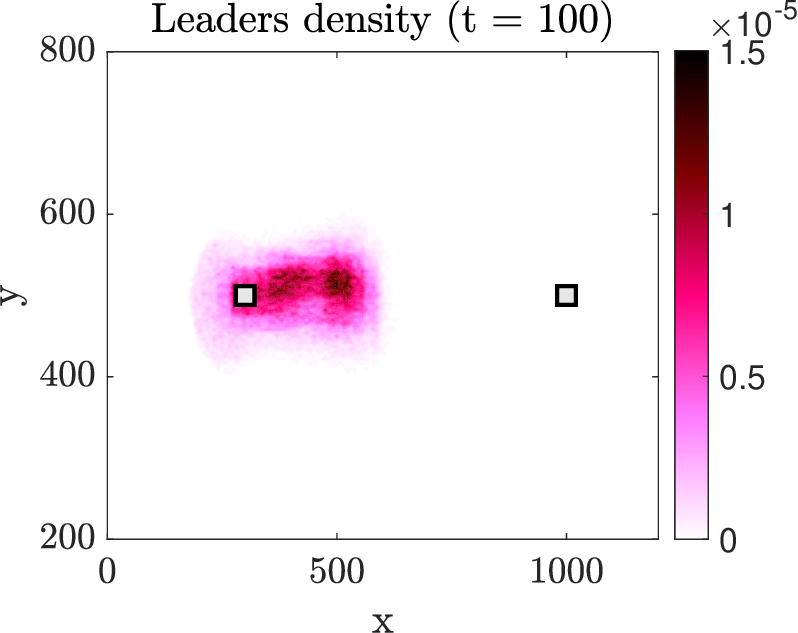}
	\includegraphics[width=0.327\linewidth]{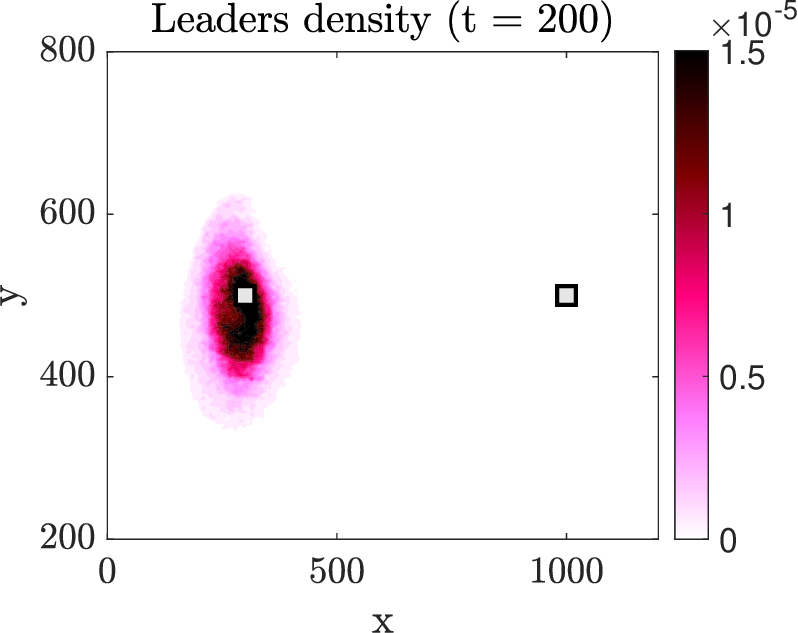}
	\caption{Three snapshots of the 2D dynamics at the mesoscopic level described in \eqref{eq:boltz_strong} and simulated by means of the Asymptotic Nanbu Algorithm \ref{alg_binary} with $\lambda$ evolving with rates \eqref{eq:rates_test} as in Algorithm \ref{alg_lambda}, taken at time $t=50$, $t=100$ and $t=200$ and assuming $C_{ctr} = C_{src}=0.75$.  First row: evolution of the total density and in red the velocity vector field of the leaders' density. Second row: evolution of the leaders' density.  }
	\label{fig:meso2D_dynamics}
\end{figure}

\paragraph{Microscopic case}
In Figure \ref{fig:micro2D_dynamics} three snapshots of the dynamics at time $t=50$, $t=100$, $t=200$. Agents that at time $t=0$ were in the leaders status change immediately  their labels since no followers are positioned around them. A large group is attracted by the food source on the left while the remaining part moves subjected just to attraction, repulsion and alignment forces without being attracted by the other food source. Once this smaller group  moves far away from the main group, leaders start to be attracted to the centre of mass. In late time, all agents join and move toward the food source on the left.  The initial separation can be observed more clearly at the mesoscopic level, as you can see in Figure \ref{fig:meso2D_dynamics}, or in the video simulation (link at the end of the subsection). 
%In Figure \cref{fig:micro2D_twofood_percentages} the agents percentages for the dynamics with leaders. 
%\begin{figure}[H]
%	\centering
%	\includegraphics[width=0.495\linewidth]{figure/micro2D/percentages_micro2D_two_food}
%	
%	\caption{Agents percentages: microscopic two dimensional model with two food sources. In black and in red the agents percentages computed by counting the effective number of followers and leaders per time steps.
%	}
%	\label{fig:micro2D_twofood_percentages}
%\end{figure}
%%%%%%%%%%%%%%%%%%%%%%%%%%%%%%%%%%%%%%%%%%%%%%%%%%%%%%%%%%%%%%%%%%%%%%%%%%%%%%%%%%%%%%%%%%%%%%%%%%%%%%%%%%%%%%%%
\paragraph{Mesoscopic case}
In Figure \ref{fig:meso2D_dynamics} three snapshots of the dynamics at time $t=50$, $t=100$, $t=200$. In the first row the time evolution of the total density and in red the velocity vector field of the leaders. In the second row the time evolution of the leaders' density. The behaviour is similar to the one observed in the microscopic case.
%% In Figure \cref{fig:meso2D_percentage} the mass percentages in time.  On the left $\phi = 0$ and on the right $\phi = 0.5$.

\begin{figure}[tbhp]
	\centering
	\includegraphics[width=0.48\linewidth]{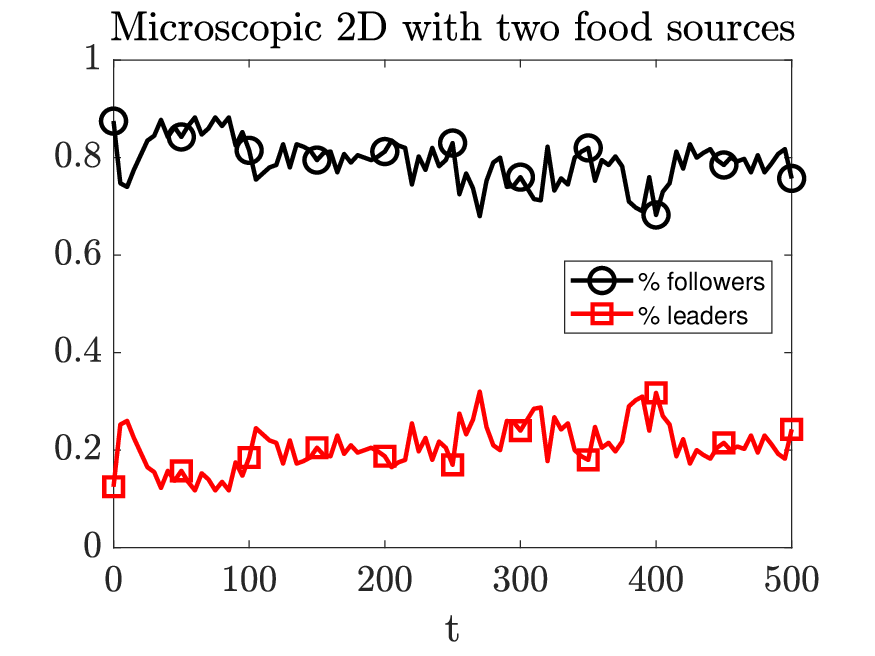}
	\includegraphics[width=0.48\linewidth]{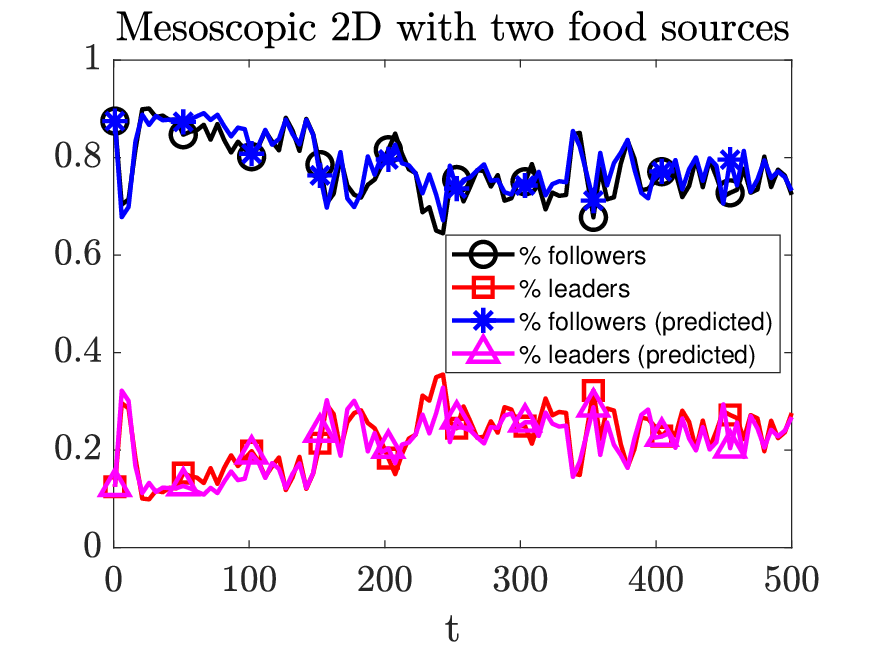} 
	\caption{ Agents percentages: microscopic (on the left) and mesoscopic (on the right) two dimensional model with two food sources. In black and in red the agents percentages computed by counting the effective number of followers and leaders per time steps. In blue and in magenta the stationary solution to the master equation \eqref{eq:master_constant}.   }
	\label{fig:percentages}
\end{figure}

In Figure \ref{fig:percentages} the time evolution of the percentages of leaders and followers for the two dimensional spatial test with two food sources. The agents percentages have been computed both by counting the effective number of followers and leaders per time steps, and as stationary solution to the master equation \eqref{eq:master_constant}. 
The densities reach the positive equilibrium defined in equation \eqref{eq:stationary}. Indeed, the transition rates defined in \eqref{eq:rates_test} 
can be approximated by constant values. In particular, for any fixed time $t>0$, for any $\lambda \in \{0,1\}$ and for any $x\in \RR^{d}$ we have 
\begin{equation}
	\pi_{L\to F}(x,\lambda;f,t) =  \mathbb{E}_x\left( 	\pi_{L\to F}(\cdot)\right), \qquad 	\pi_{F\to L}(x,\lambda;f,t) = \bar{\beta}(t) =  \mathbb{E}_x\left( 	\pi_{F\to L}(\cdot)\right),
\end{equation}
%where 
%\begin{equation}
%	\bar{\alpha}(t) = \mathbb{E}_x\left( 	\pi_{L\to F}(\cdot)\right) , \qquad \bar{\beta}(t) =  \mathbb{E}_x\left( 	\pi_{F\to L}(\cdot)\right) ,
%\end{equation}
with $\mathbb{E}_x(\cdot)$ denoting the mean value with respect to $x$. \\
Similar results can be obtained for the 2D model without food sources since the transition rates are constants values,  by definition. 
The videos of the simulations of this subsection are available at \href{https://drive.google.com/drive/folders/1KZ2vOeMNzmx8MBMX-_0kulGGcjQr_shi?usp=share_link}{[VIDEO]}.
%%%%%%%%%%%%%%%%%%%%%%%%%%%%
%%%%%%%%%%%%%%%%%%%%%%%%%%%%%%%%%%%%%%%%%%%%%

\subsubsection{Test 2D: one food source}\label{sec:test2D_optimal_control} 
\begin{figure}[tbhp]
	\centering
	\includegraphics[width=0.46\linewidth]{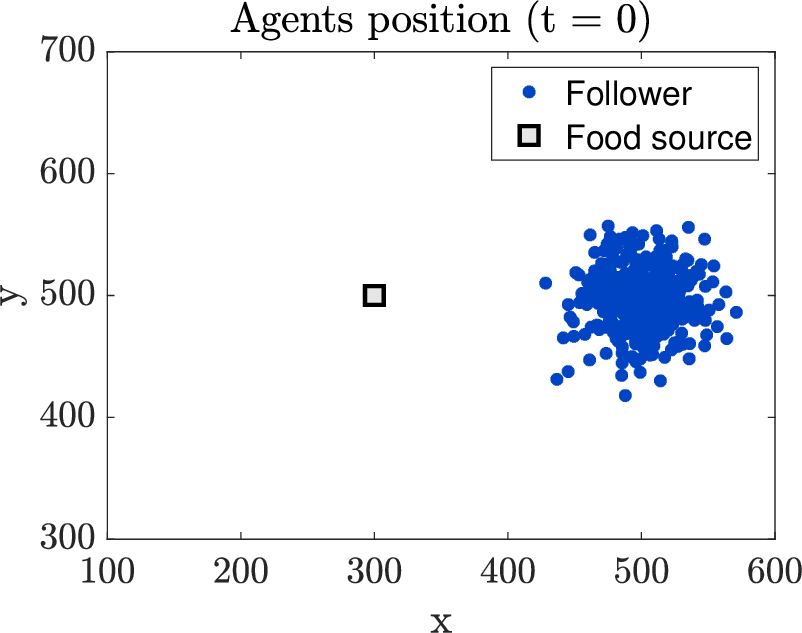}
	\includegraphics[width=0.48\linewidth]{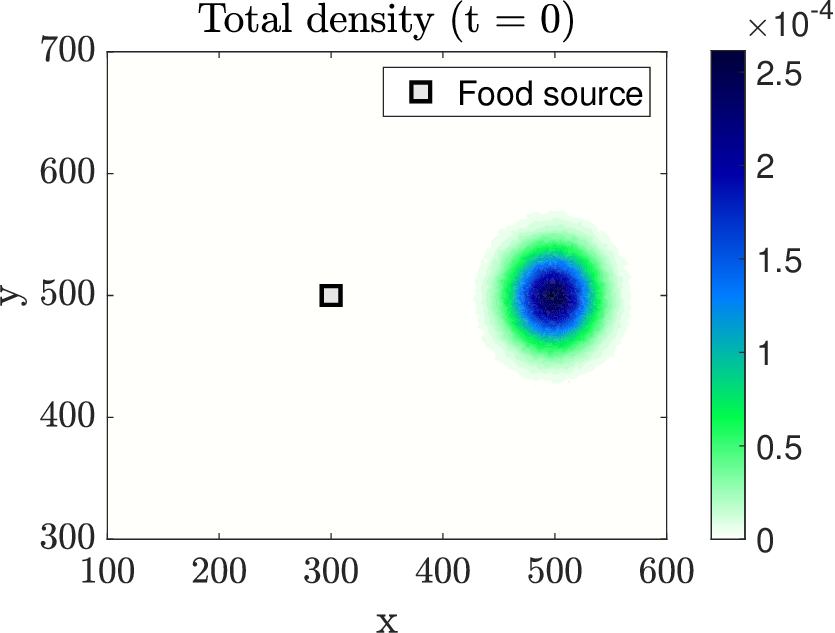}
	\caption{Initial configuration in the 2D case with one food source.}
	\label{fig:initial_configuration_opt}
\end{figure}

Assume the model includes one food source located in $x_1^{src} = (300,500)$. Run the simulation until time $T = 120$. 
Suppose labels change aiming at organizing agents toward a common target, that in this case is supposed to be the food source $x_1^{src}$. In particular, assume $\lambda$ varies with rates \eqref{eq:rates_opt} with $\bar{\alpha} = 0.7$ and $\underline{\alpha} = 0.3$.
In Figure \ref{fig:initial_configuration_opt} the initial configuration for both the microscopic and mesoscopic case.

\begin{figure}[tbhp]
	\centering
	\includegraphics[width=0.327\linewidth]{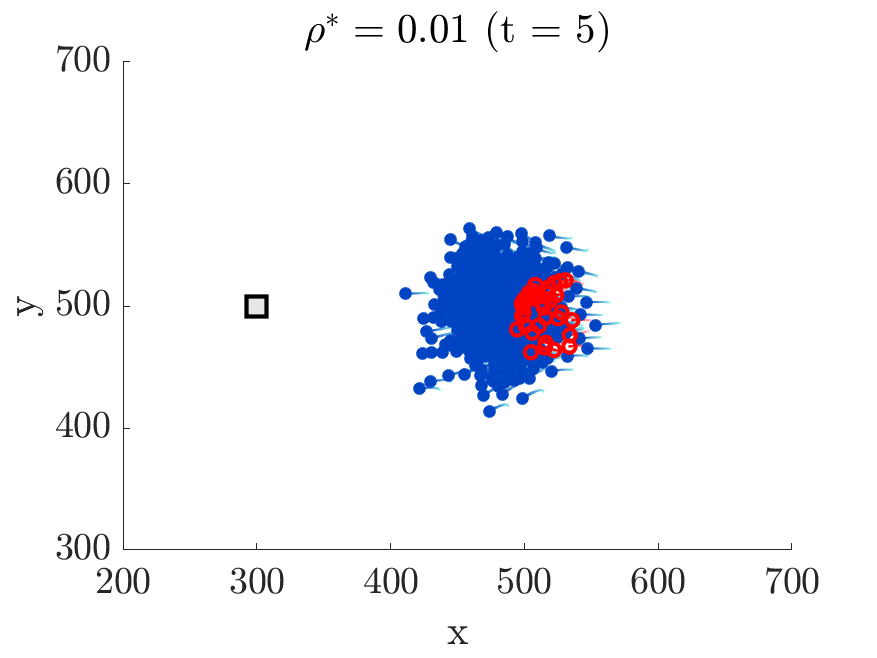}
	\includegraphics[width=0.327\linewidth]{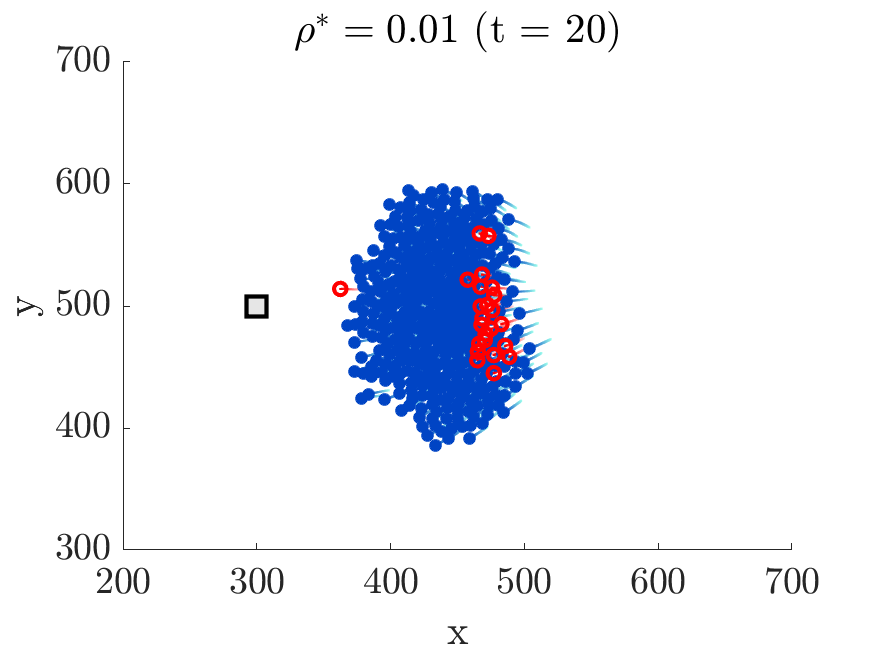}
	\includegraphics[width=0.327\linewidth]{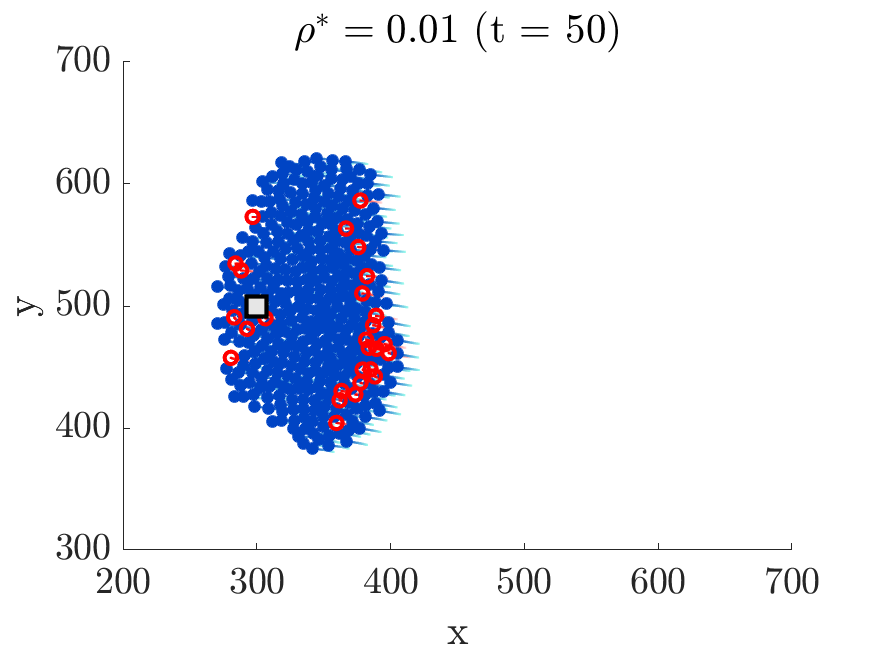}\\
	\caption{Three snapshots of the 2D dynamics at the microscopic level described in \eqref{eq:dynamics} with $\lambda$ evolving with rates \eqref{eq:rates_opt} taken at time $t=5$, $t=20$ and $t=50$ and assuming $C_{ctr} = C_{src} = 0$.   We represent in blue the agents in the followers status and in red the ones in the leaders status.
	}
	\label{fig:micro2D_dynamics_opt}
\end{figure}

\paragraph{Microscopic case}
In Figure \ref{fig:micro2D_dynamics_opt} three snapshots of the dynamics at time $t=5$, $t=20$ and $t=50$ with leaders and with $\mathcal{G}[f](\cdot)$ chosen as in \eqref{eq:G_test1}. Followers are driven by leaders and reach the target position. 

\begin{figure}[tbhp]
	\centering
	\includegraphics[width=0.327\linewidth]{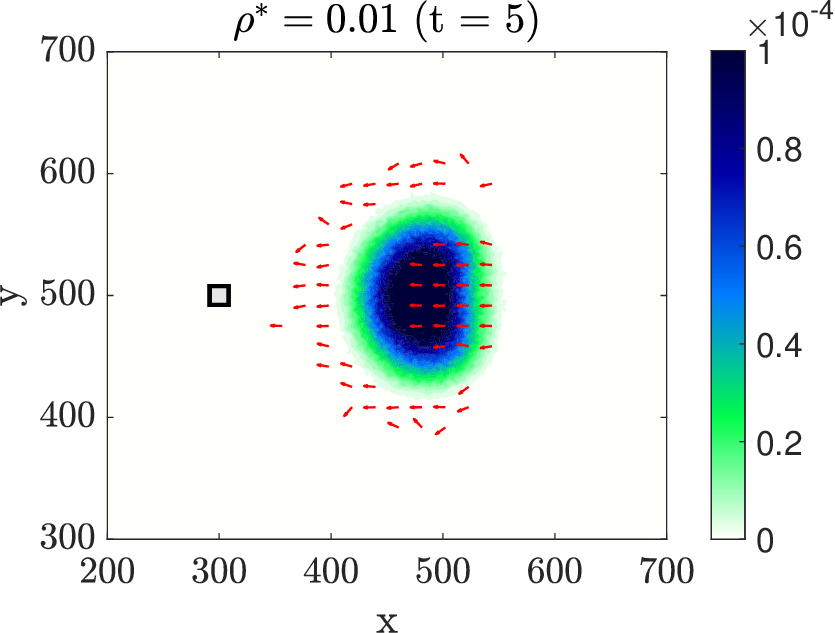}
	\includegraphics[width=0.327\linewidth]{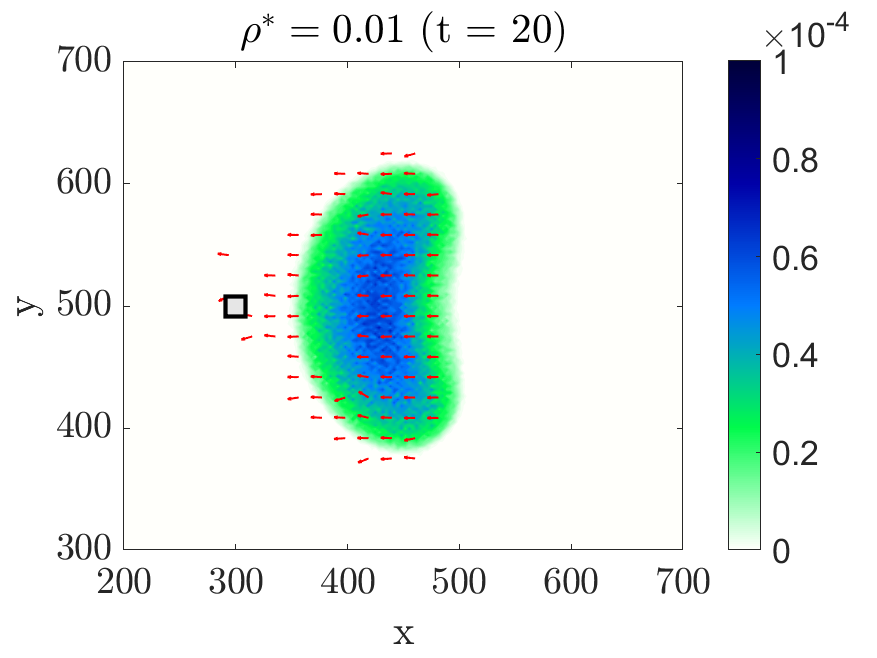}
	\includegraphics[width=0.327\linewidth]{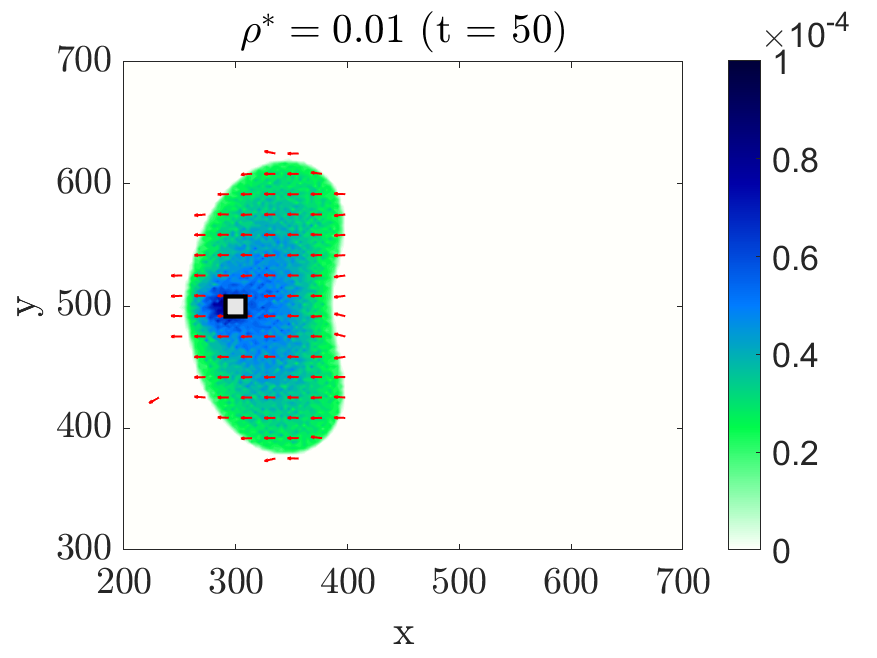}\\
	\includegraphics[width=0.327\linewidth]{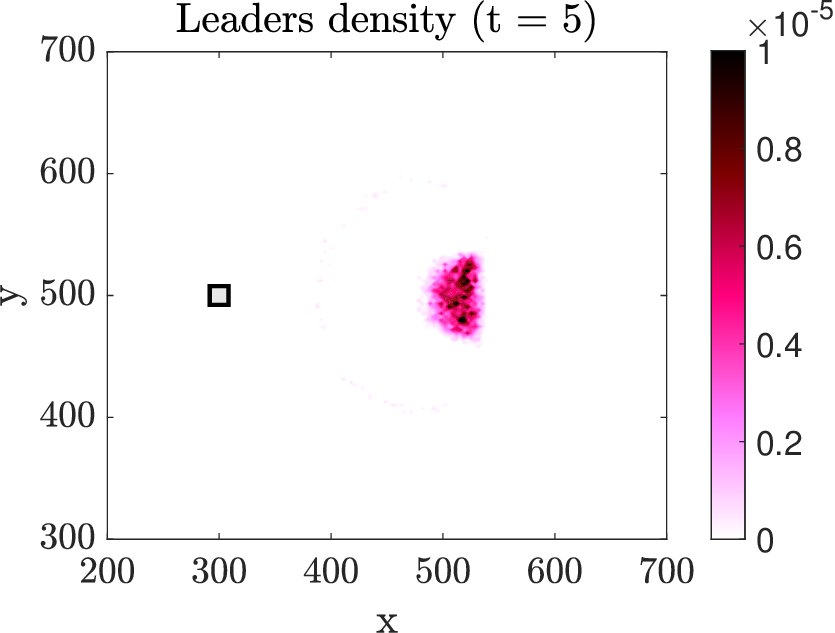}
	\includegraphics[width=0.327\linewidth]{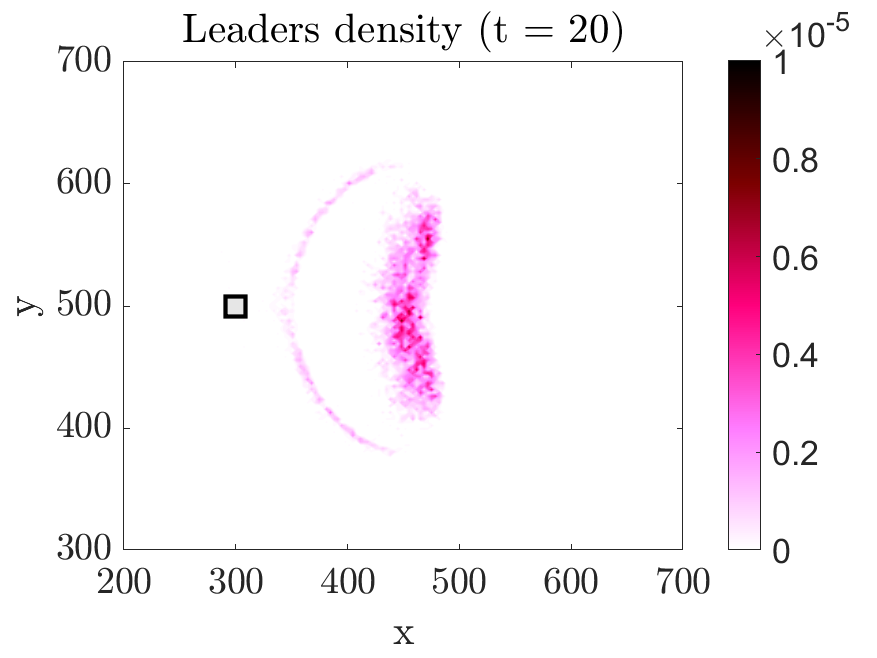}
	\includegraphics[width=0.327\linewidth]{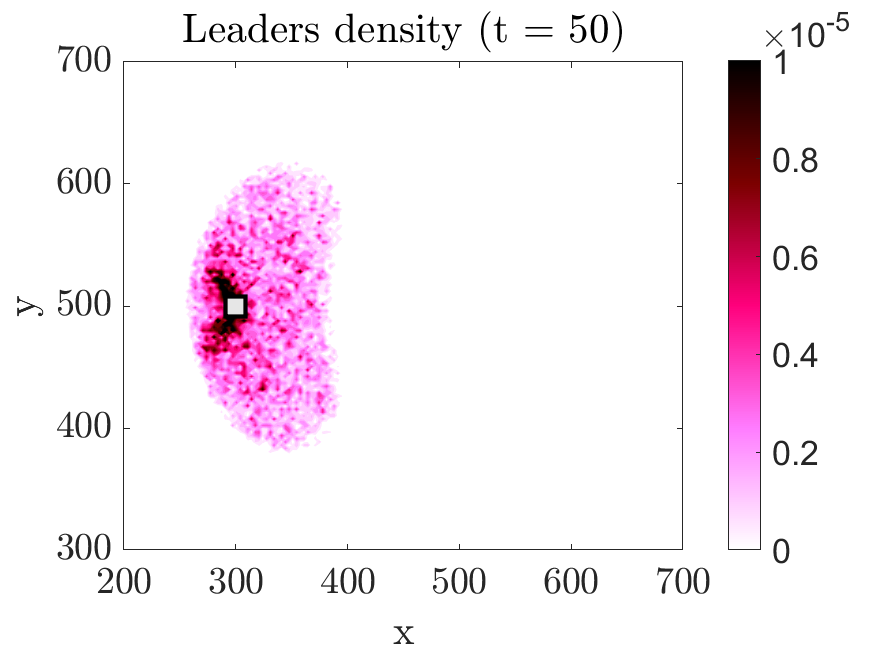}
	\caption{Three snapshots of the 2D dynamics  at the mesoscopic level described in \eqref{eq:boltz_strong} and simulated by means of the Asymptotic Nanbu Algorithm \eqref{alg_binary} with $\lambda$ evolving with rates \eqref{eq:rates_opt} as in \eqref{alg_lambda}, taken at time $t=5$, $t=20$ and $t=50$ and assuming $C_{ctr} = C_{src} = 0$. First row: evolution of the total density and in red the velocity vector field of the leaders. Second row: evolution of the leaders density. 
	}
	\label{fig:meso2D_dynamics_opt}
\end{figure}

\paragraph{Mesoscopic case}
In Figure \ref{fig:meso2D_dynamics_opt} three snapshots of the dynamics at time $t=5$, $t=20$, $t=50$ with transition rates depending on the orientation according to \eqref{eq:rates_opt}. In the first row, the evolution of the whole mass, and in red the velocity vector field of the leaders. In the second row, the evolution of the leaders mass. 

\begin{figure}[tbhp]
	\centering
	\includegraphics[width=0.327\linewidth]{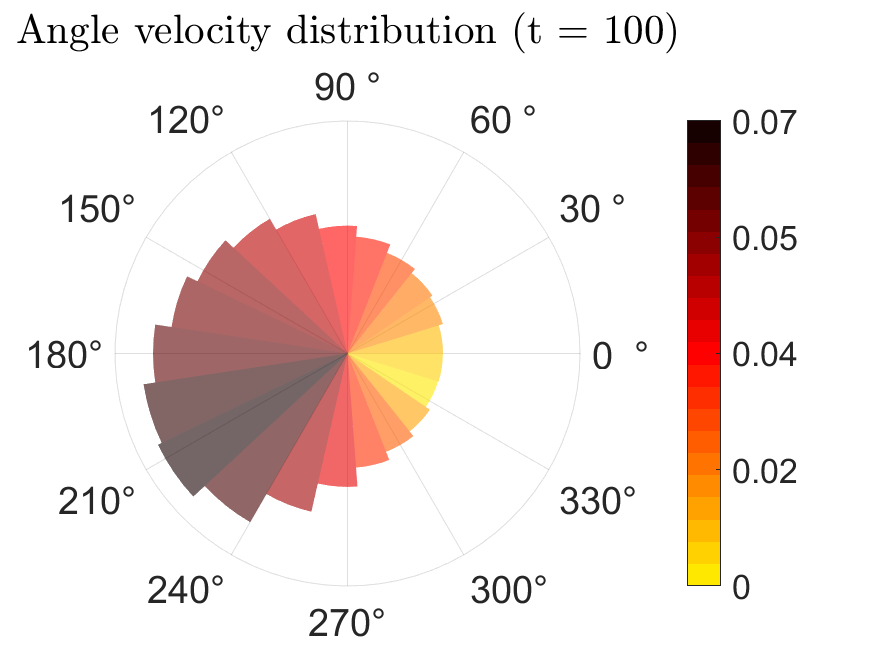}
	\includegraphics[width=0.327\linewidth]{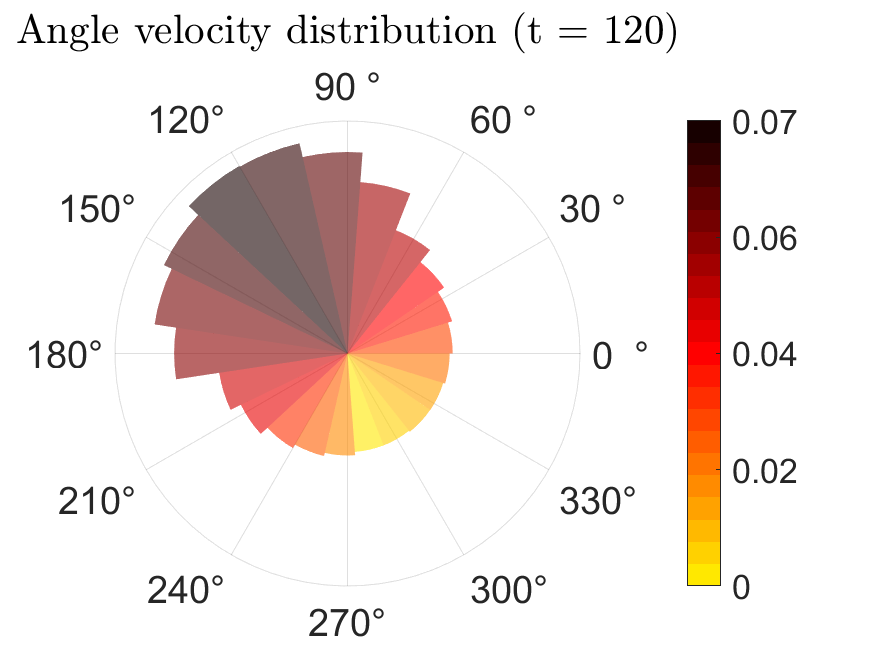}
	\includegraphics[width=0.327\linewidth]{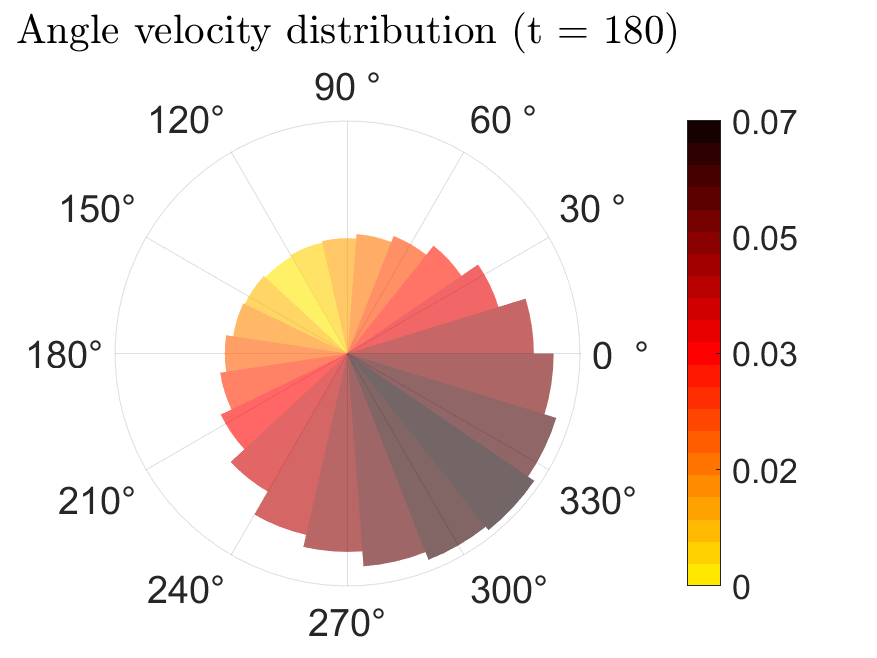}\\
	\includegraphics[width=0.327\linewidth]{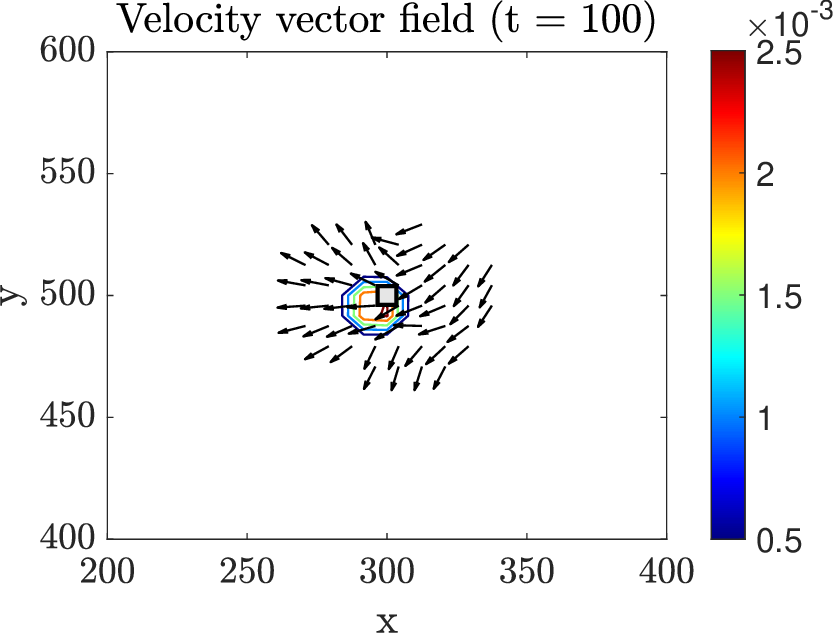}
	\includegraphics[width=0.327\linewidth]{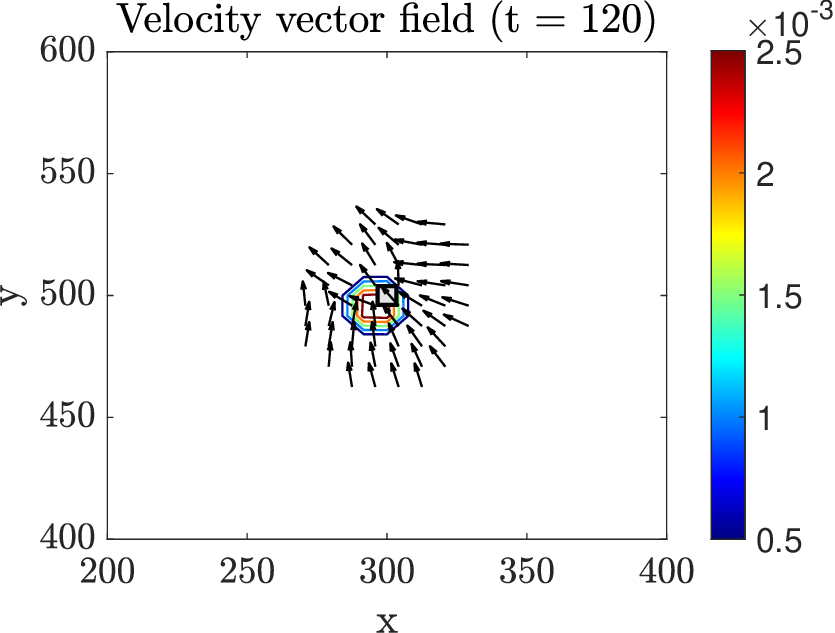}
	\includegraphics[width=0.327\linewidth]{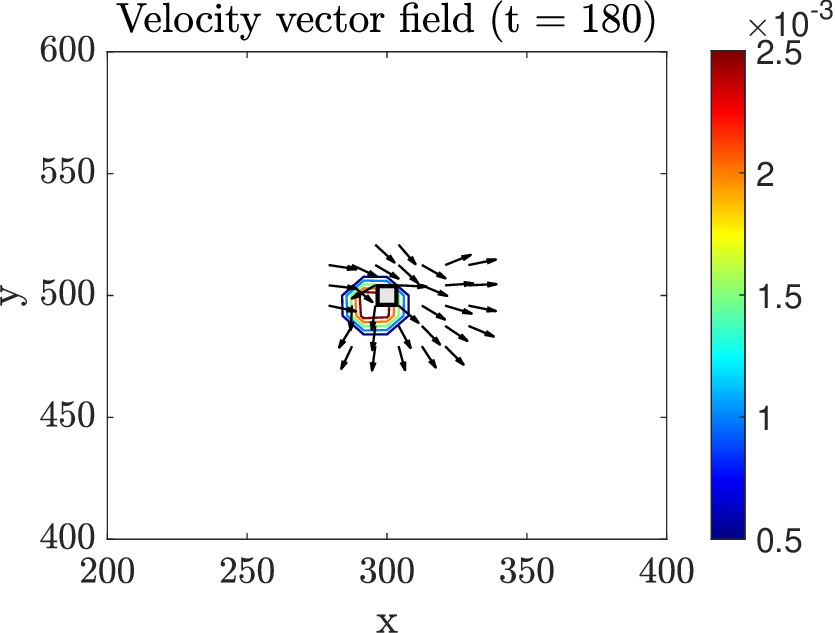}
	\caption{Angle velocity distribution (first row) and velocity vector field (second row) at time $t = 100$, $t=120$, $t=180$ obtained by simulating the dynamics in \eqref{eq:boltz_strong} by means of the Algorithm \ref{alg_binary} with $\lambda$ evolving with rates \eqref{eq:rates_opt} as in Algorithm \ref{alg_lambda}. 
		\label{fig:velocity}
	}
\end{figure}

In Figure \ref{fig:velocity} we report in the first row  the angle velocity distribution at time $t=100$, $t=120$, $t=180$ and in the second row the correspondent velocity vector field, outlining the milling behaviour around the target positions $\overline x=x_1^{src}$.
The videos of the simulations of this subsection are available at \href{https://drive.google.com/drive/folders/1A-opPh85UoA8POSsZb6gr73eF4vW81b8?usp=drive_link}{[VIDEO]}.
%%%%%%%%%%%%%%%%%%%%%%%%%%%%%%%%%%%%%%%%%%%%%
%%%%%%%%%%%%%%%%%%%%%%%%%%%%%%%%%%%%%%%%%%%%%

%%%%%%%%%%%%%%%%%%%%%%%%%%%%%%%%%%%%%%%%%%%%%%%%%%%%%%%%%%
%%%%%%%%%%%%%%%%%%%%%%%%%%%%%%%%%%%%%%%%%%%%%%%%%%%%%%%%%%%%%%%%%%%%%%%%%%%%%%%%%%%%%%%%%%%%%%%%%%%%%%%%%%%%%%%%%%
\subsection{Numerical test in 3D with two food sources} \label{sec:test3D} 
\begin{figure}[tbhp]
	\centering
	\includegraphics[width=0.48\linewidth]{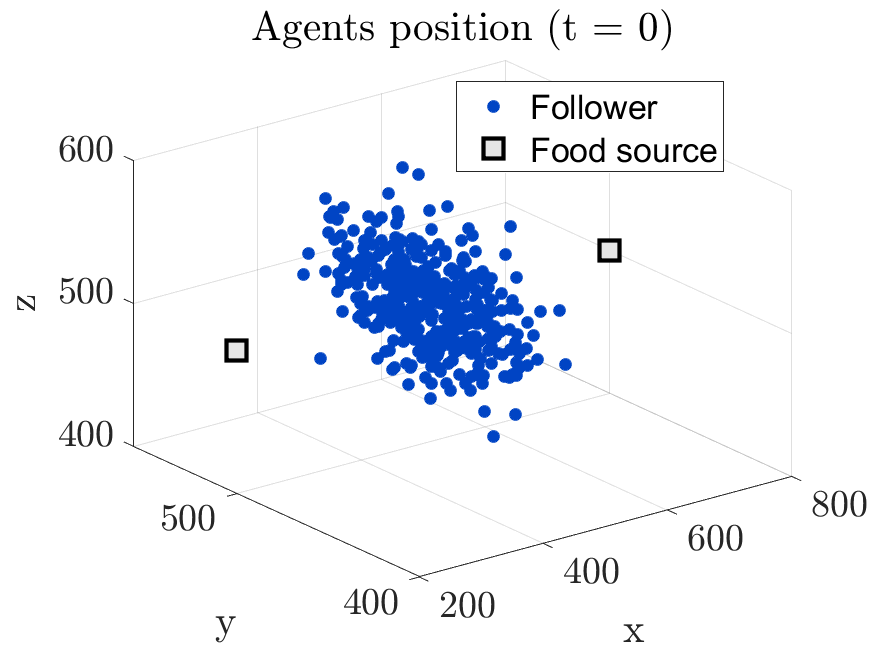}
	\includegraphics[width=0.48\linewidth]{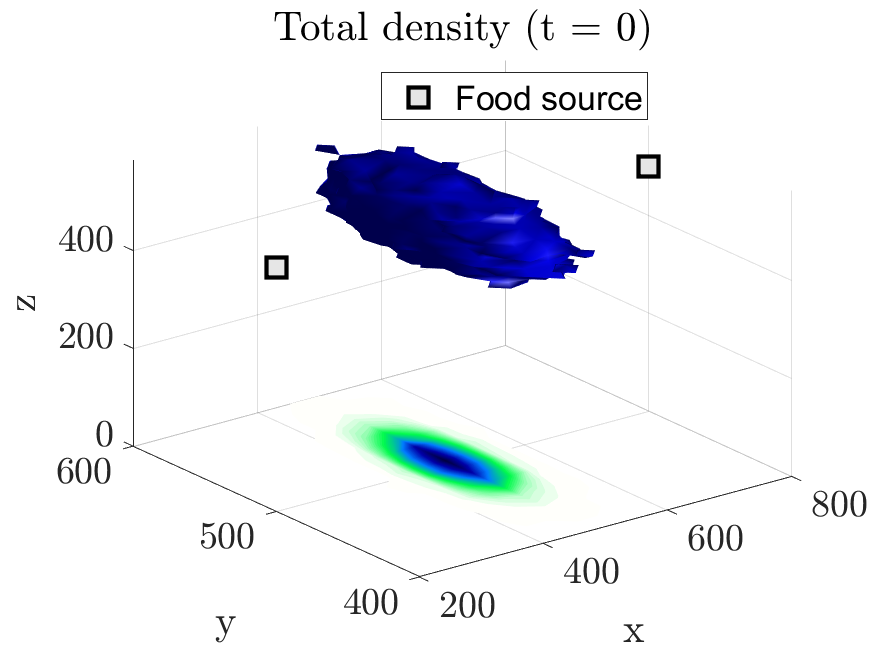}
	\caption{Initial configuration in the 3D case with food sources.}
	\label{fig:micro3D_initial_configuration}
\end{figure}

We consider the three dimensional model in space and velocity, simulating the swarming dynamics up time $T=200$.
Initially agents are normally distributed with mean $\mu = 500$ and
variance $\sigma^2 = 25^2$ in both spatial and velocity dimension, and are all in the followers status.
We report in Figure \ref{fig:micro3D_initial_configuration} the initial configuration for both the microscopic and mesoscopic case. For the mesoscopic case, we also depict on the $(x,y)$ plane the projection of the spatial density. 

\begin{figure}[tbhp]
	\centering
	\includegraphics[width=0.327\linewidth]{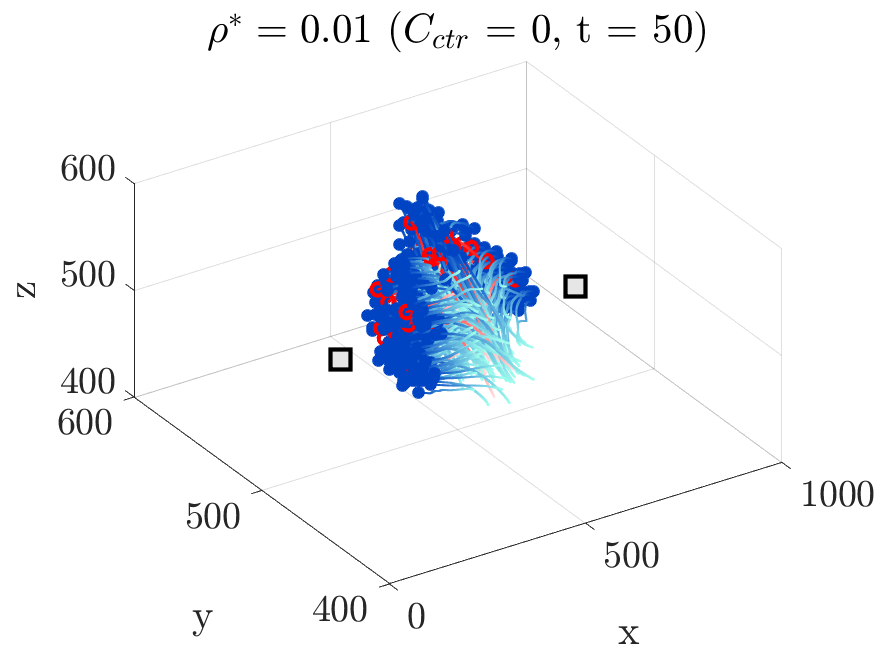}	\includegraphics[width=0.327\linewidth]{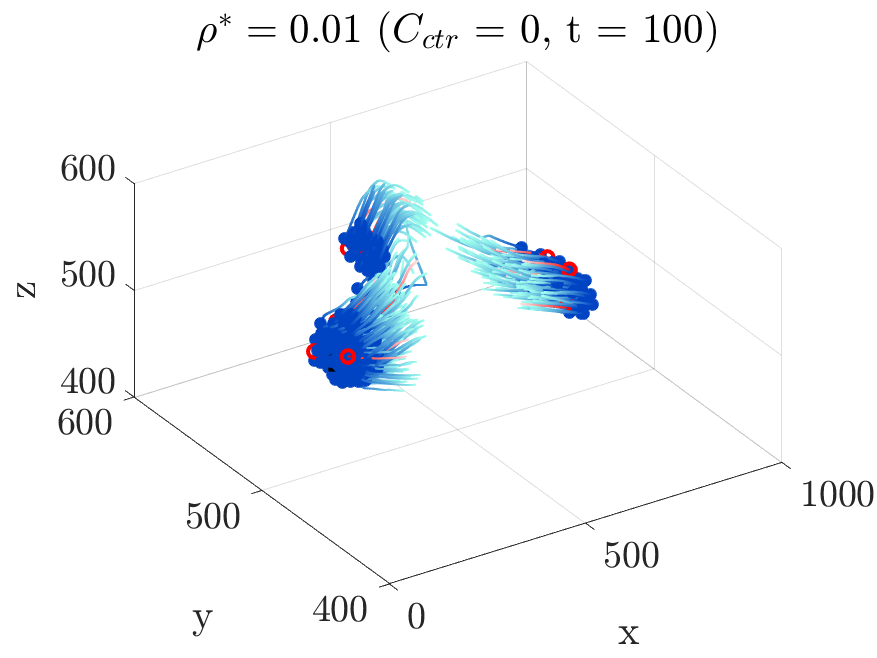}
	\includegraphics[width=0.327\linewidth]{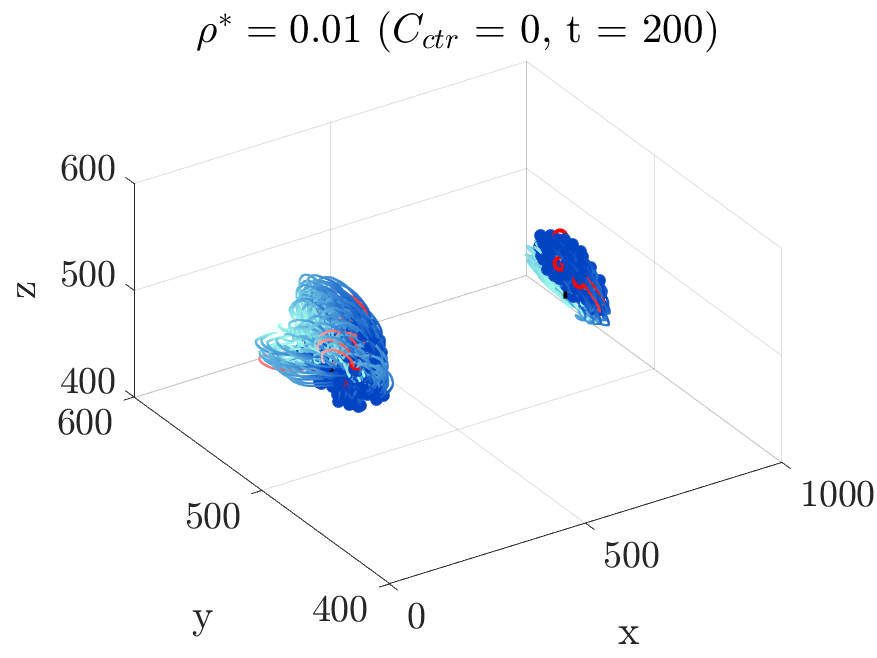}\\
	\includegraphics[width=0.327\linewidth]{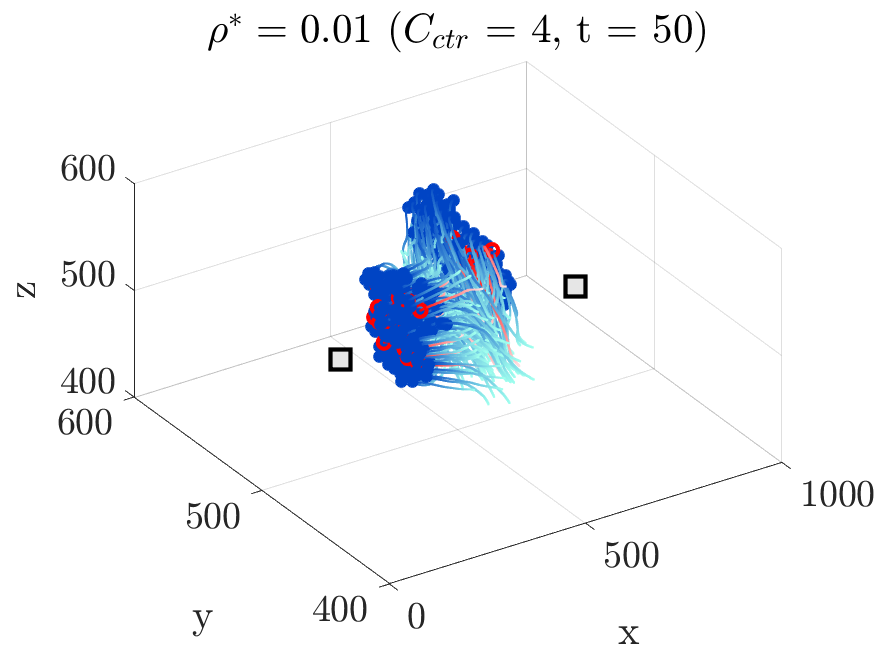}	\includegraphics[width=0.327\linewidth]{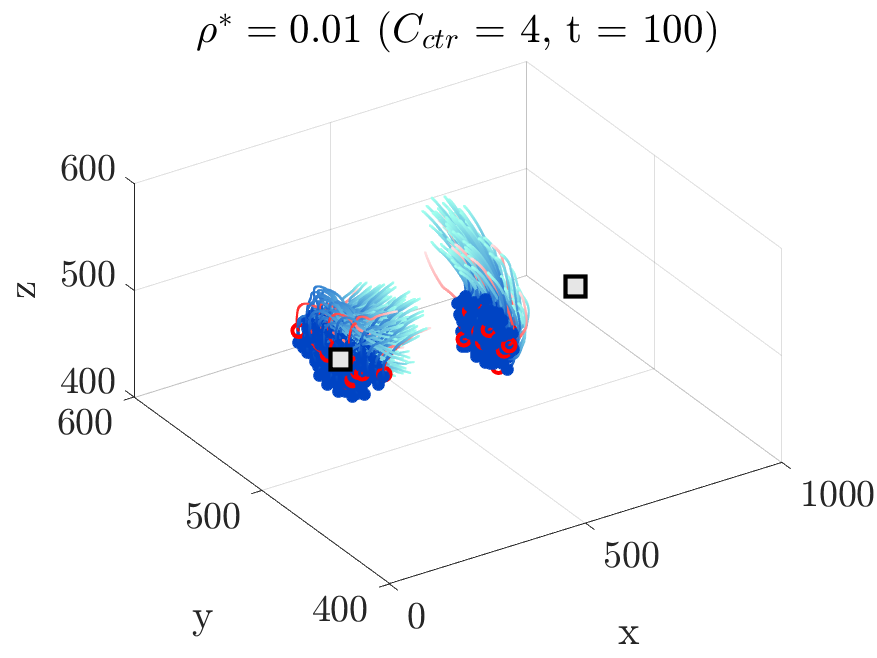}
	\includegraphics[width=0.327\linewidth]{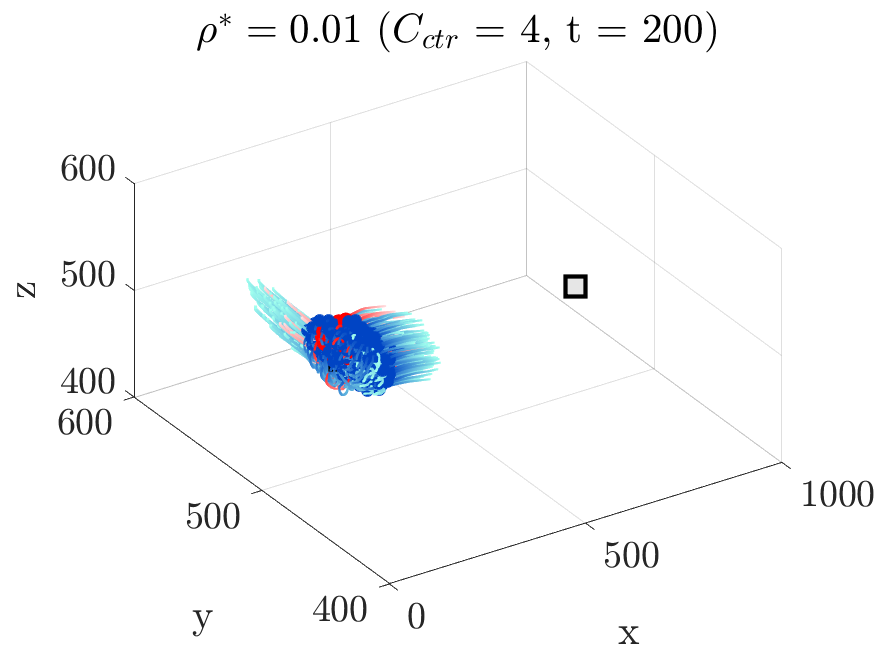}
	\caption{Three snapshots of the 3D dynamics at the microscopic level described in \eqref{eq:dynamics} with $\lambda$ evolving with rates \eqref{eq:rates_test} taken at time $t=50$, $t=100$ and $t=200$. First row, $C_{ctr} = 0$, $C_{src}=0.75$. Second row, $C_{ctr} = 4$, $C_{src}=0.75$. We represent in blue the agents in the followers status and in red the ones in the leaders status.}
	\label{fig:micro3D_dynamics}
\end{figure}

Assume the two food sources to be located in $x_1^{src} = (200,500,500)$ and $x_2^{src} = (800,500,500)$. Suppose that leaders emerge with density-dependent transition rates as defined in \eqref{eq:rates_test} and where we assume the constants to be $q_L = 4\times 10^{-3}$ and $q_F = 3\times 10^{-3}$.   
\paragraph{Microscopic case}
In Figure \ref{fig:micro3D_dynamics} three snapshots of the dynamics at time $t=50$, $t=100$, $t=200$. First row: $C_{ctr}= 0$, $C_{src} = 0.75$. Agents split in two groups moving toward the two food sources. Second row: $C_{ctr}= 4$, $C_{src} = 0.75$. At final time agents move toward one of the two food sources.
%In Figure \cref{fig:micro3D_food_percentages} the agents percentages for the 3D tests with food sources.  
%\begin{figure}[H]
%	\centering
%	\includegraphics[width=0.495\linewidth]{figure/micro3D/percentages_micro3D_SIfoodNOcentre}
%	\includegraphics[width=0.495\linewidth]{figure/micro3D/percentages_micro3D_SIfoodSIcentre}
%	\caption{Agents percentages: microscopic three dimensional model with no food sources. On the left $C_{centre} = 0$, on the right $C_{centre} = 4$. In black and in red the agents percentages computed by counting the effective number of followers and leaders per time steps.
%	}
%	\label{fig:micro3D_food_percentages}
%\end{figure}
%%%%%%%%%%%%%%%%%%%%%%%%%%%%%%

\begin{figure}[tbhp]
	\centering
	\includegraphics[width=0.327\linewidth]{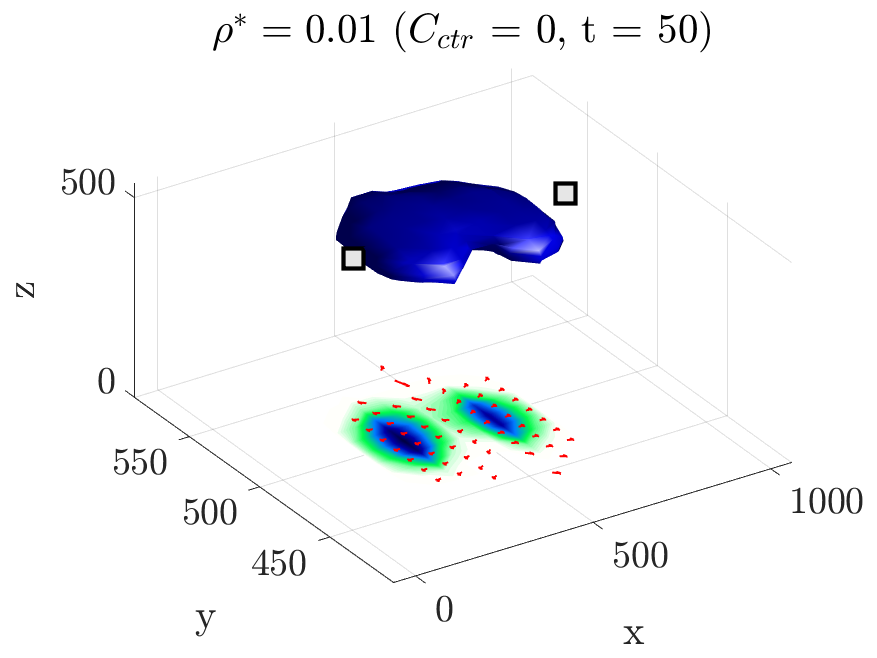}	\includegraphics[width=0.327\linewidth]{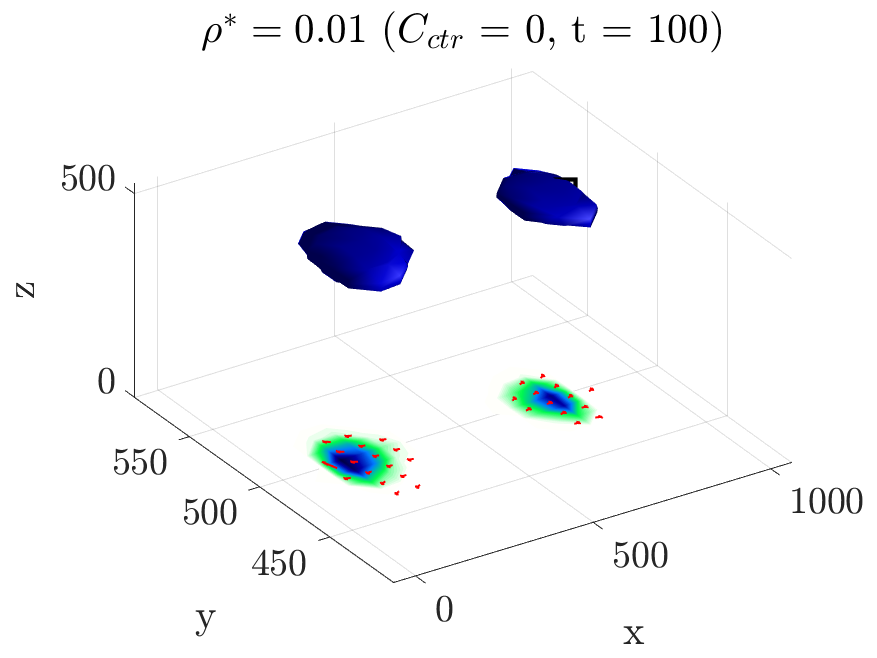}
	\includegraphics[width=0.327\linewidth]{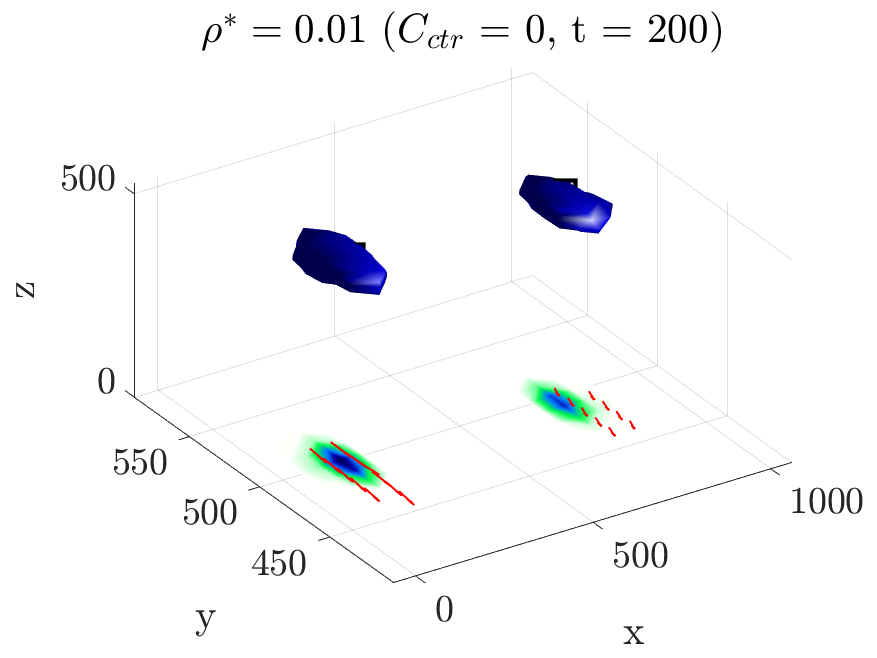}\\
	\includegraphics[width=0.327\linewidth]{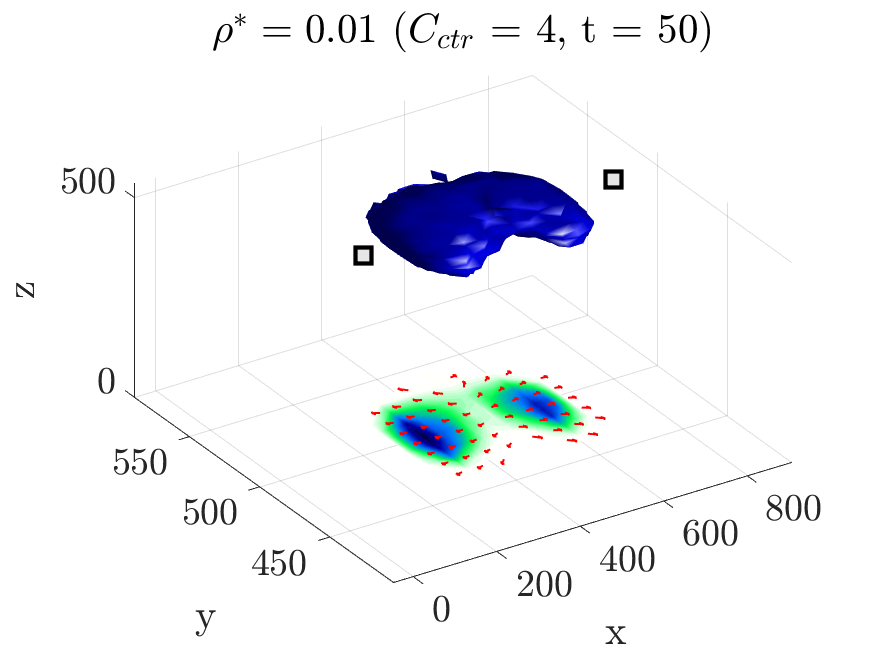}	\includegraphics[width=0.327\linewidth]{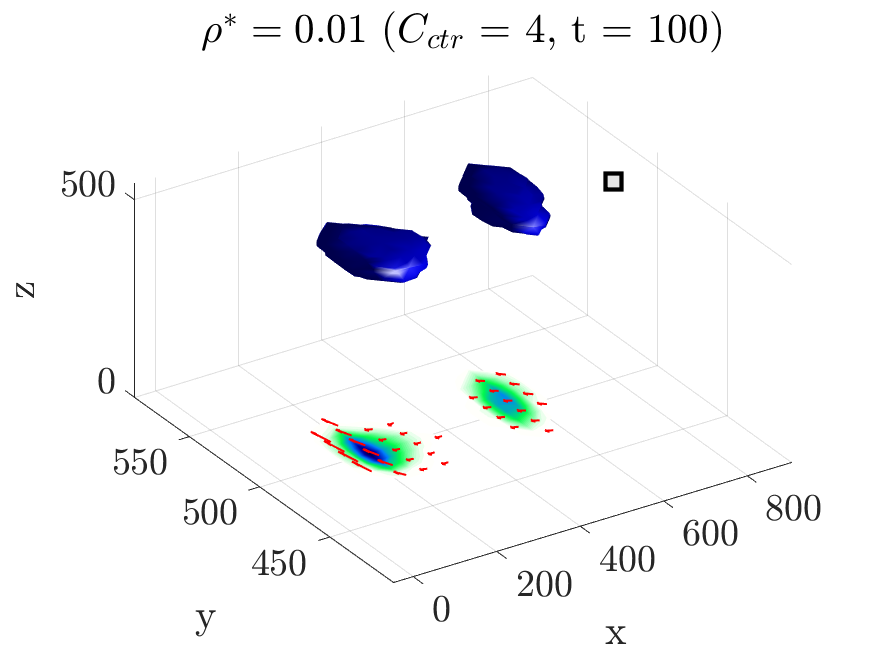}
	\includegraphics[width=0.327\linewidth]{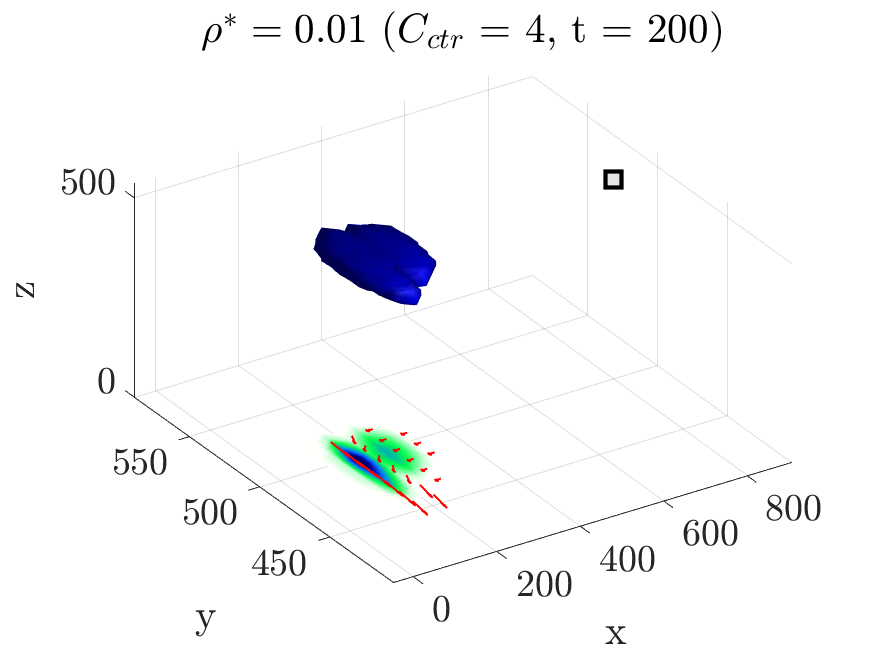}
	\caption{Three snapshots of the 3D dynamics at the mesoscopic level described in \eqref{eq:boltz_strong} and simulated by means of the Asymptotic Nanbu Algorithm \ref{alg_binary} with $\lambda$ evolving with rates \eqref{eq:rates_test} as in Algorithm \ref{alg_lambda}, taken at time $t=50$, $t=100$ and $t=200$. First row, $C_{ctr} = 0$, $C_{src}=0.75$. Second row, $C_{ctr} = 4$, $C_{src}=0.75$. In red the velocity vector field. }
	\label{fig:meso3D_dynamics}
\end{figure}

\paragraph{Mesoscopic case}
In Figure \ref{fig:meso3D_dynamics} we report  three snapshots of the dynamics at time $t=50$, $t=100$, $t=200$ and in red the velocity vector field. We add also the density distribution of the whole flock projected over the plane $(x,y)$ and in red the leaders velocity vector field. First row: $C_{crt}= 0$, $C_{src} = 0.75$. Second row: $C_{ctr}= 4$, $C_{src} = 0.75$. The behaviour is similar to the one in the microscopic case. 
%In Figure \cref{fig:3D_food_percentages} the agents percentages for the 3D tests with food sources in both the microscopic and mesoscopic case.  
%\begin{figure}[H]
%	\centering
%	\includegraphics[width=0.495\linewidth]{figure/percentages/percentages3D_SIfoodNOcentre}
%	\includegraphics[width=0.495\linewidth]{figure/percentages/percentages3D_SIfoodSIcentre}
%	\caption{Agents percentages: 3D model with two food sources. On the left $C_{centre} = 0$,  on the right $C_{centre} = 4$. Markers have been added just to indicate different lines. 
%	}
%	\label{fig:3D_food_percentages}
%\end{figure}
The videos of the simulations of this subsection are available at \href{https://drive.google.com/drive/folders/1RS0LyB18zSoyKqVmy99NBwsIyRvR-t_m?usp=share_link}{[VIDEO].}
\section{Conclusions}\label{sec:conclusion}
In this paper, we have studied collective behaviour of birds under a follower-leaders dynamics, starting from the model presented in \cite{cristiani2021all}.  Through the emergence of leaders, we recover the ability to split the initial configuration and initiate directional changes without the need of external influences. We derived a kinetic model to effectively depict the motion of a large swarm with transient leadership and topological interactions, and subsequently we simulated the dynamics introducing a novel stochastic particle method.
A significant emphasis was placed on studying topological interactions. We tackled the issue of the numerical evaluation of Nearest Neighbors reducing the computational costs of the search from quadratic to logarithmic by optimally organizing agents in a binary tree and performing a $k$-NN search. Moreover, we directed our attention to transient leadership, showcasing how labels can change over time, particularly for driving agents towards a common target. Various strategies for leaders' emergence were explored, including the extension to a continuous space of labels, as detailed  in Remark \ref{remark_continuos} and Remark \ref{remark_kincont}.  The numerical experiments show the computational feasibility of the topological sphere and simulation of the dynamics at mesoscopic level by means of the proposed novel stochastic algorithm \ref{alg_binary}. Moreover, these simulation asses the flexibility of the model of describing non-trivial phenomena such as flock splitting, and attraction towards food sources, or nesting areas, showing that these mechanisms are triggered by the presence of transient leadership and topological interactions. It would be intriguing to describe the original model from a kinetic viewpoint, reintroducing delay, which, as discussed in \cite{cristiani2021all}, appears to play a crucial role in achieving desired configurations.
Finally, several questions arise concerning the study of non-local terms in high dimensions. For instance, it could be beneficial to further enhance the numerical scheme implemented, focusing on other useful strategies for approximating topological interactions.
\section*{Acknowledgments}
GA and FF  are members of INdAM GNCS. GA was partially supported by the MIUR-PRIN Project 2022 No. 2022N9BM3N	``Efficient numerical schemes and optimal control methods for time-dependent PDEs”  and by MUR-PRIN Project 2022 PNRR No. P2022JC95T, ``Data-driven discovery and control of multi-scale interacting artificial agent systems", financed by the European Union - Next Generation EU.
	\bibliographystyle{abbrv}
	\bibliography{biblio}
\end{document}